\documentclass{compositio}
\usepackage[latin1]{inputenc}

\usepackage{amsmath}
\usepackage{amscd}
\usepackage{mathrsfs}
\usepackage[all]{xy}
\usepackage{eufrak}
\usepackage{graphicx}
\usepackage[dvips,colorlinks=true,linkcolor=blue]{hyperref}
\usepackage{makeidx}
\usepackage{index}
\usepackage{multicol}

\renewcommand{\a}{\mathcal{A}}

\newcommand{\Autc}{\mathrm{Aut}_K^{\mathrm{cont}}}

\newcommand{\B}{\mathrm{B}}
\newcommand{\bs}[1]{\boldsymbol{#1}}
\newcommand{\CVD}{$\Box$}

\newcommand{\C}{\mathrm{C}}
\newcommand{\Conf}{\mathrm{Conf}}

\newcommand{\D}{\mathrm{D}}
\newcommand{\Def}{\mathrm{Def}}

\newcommand{\e}{\mathbf{e}}

\newcommand{\Ed}{\mathcal{E}^{\dag}}

\newcommand{\Endc}{\mathrm{End}_K^{\mathrm{cont}}}
\newcommand{\et}{\mathrm{e}}
\newcommand{\F}{\mathrm{F}}
\newcommand{\Fb}{\bar{\mathrm{F}}}
\newcommand{\G}{\mathrm{G}}

\renewcommand{\H}{\mathcal{H}}
\newcommand{\Hd}{\mathcal{H}^{\dag}}
\newcommand{\Hom}{\mathrm{Hom}}

\renewcommand{\L}{\mathrm{L}}

\newcommand{\M}{\mathrm{M}}

\newcommand{\Mod}{\mathrm{Mod}}

\newcommand{\N}{\mathrm{N}}
\renewcommand{\O}{\mathcal{O}}

\newcommand{\ph}[1]{\langle #1\rangle}
\newcommand{\Q}{\mathcal{Q}}
\newcommand{\simto}{\stackrel{\sim}{\to}}

\newcommand{\R}{\mathcal{R}}
\newcommand{\Rep}{\underline{\mathrm{Rep}}}

\newcommand{\V}{\mathrm{V}}

\newcommand{\W}{\mathbf{W}}

\newtheorem{theorem}{Theorem}[section]
\newtheorem{proposition}[theorem]{\textsc{Proposition}}
\newtheorem{lemma}[theorem]{\textsc{Lemma}}
\newtheorem{corollary}[theorem]{\textsc{Corollary}}
\theoremstyle{definition}
\newtheorem{definition}[theorem]{\textsc{Definition}}
\newtheorem{example}[theorem]{\textsc{Example}}

\theoremstyle{remark}
\newtheorem{remark}[theorem]{\textsc{Remark}}
\newtheorem{hypothesis}[theorem]{\textsc{Hypothesis}}

\newtheorem{notation}[theorem]{\textsc{Notation}}
\numberwithin{equation}{subsection}

\numberwithin{equation}{subsection}

\makeindex

\begin{document}

\title[$p$-Adic Confluence of $q$-Difference Equations]{$\bs{p}$-Adic Confluence of $\bs{q}$-Difference Equations}

\author{Andrea Pulita}
\email{pulita@math.uni-bielefeld.de}

\address{SFB701 Spectral Structures and Topological
Methods in Mathematics, Fakultät für Mathematik Universität
Bielefeld P.O.Box 100 131 D-33501, Office: V4-223, Bielefeld,
Germany.}




\begin{abstract}
We develop the theory of $p$-adic confluence of $q$-difference
equations. The main result is the fact that, in the $p$-adic
framework, a function is (Taylor) solution of a differential
equation if and only if it is solution of a $q$-difference
equation. This fact implies an equivalence, called
\emph{Confluence}, between the category of differential equations
and those of $q$-difference equations. We develop this theory by
introducing a category of ``\emph{sheaves}'' on the disk
$\mathrm{D}^-(1,1)$, for which the stalk at $1$ is a differential
equation, the stalk at $q$ is a $q$-difference equation if $q$ is
not a root of unity, and the stalk at a root of unity $\xi$ is a
mixed object, formed by a differential equation and an action of
$\sigma_\xi$.
\end{abstract}

\maketitle

\makeatletter
\renewcommand\tableofcontents{%
    \subsection*{\contentsname}%
    \@starttoc{toc}%
    }
\makeatother

\begin{small}
\setcounter{tocdepth}{2}
\tableofcontents
\end{small}

\setcounter{section}{0}

\section*{\textsc{Introduction}}
\addcontentsline{toc}{section}{Introduction}

The main aim of this paper is to provide a \emph{theory of
confluence} for $q$-difference equations in the $p$-adic
framework.

\if{This is done in three steps. Firstly we introduce a new class
of objects, called \emph{$\sigma$-modules} and
\emph{$(\sigma,\delta)$-modules}, expressing the notion of
``\emph{family of $q$-difference equations, for $q$ varying in a
subgroup $\ph{S}$ of $K^\times$}''. Secondly we analyse the two
functors associating to such a ``family'' a single equation of the
family, and the associated differential equation
(\emph{$q$-tangent operator}) respectively. Thirdly we prove that
there is a full subcategory, formed by the so called
\emph{constant} $\sigma$-modules or $(\sigma,\delta)$-modules, for
which these two functors are locally equivalences. In other words
for a constant object the whole ``family of $q$-difference
equations'' is completely determined, up to restrict $\ph{S}$, by
one of its equations, or by the associated differential equation.
A constant object is in fact characterized by the the property
that every equation of the family has the same \emph{generic
Taylor solution}. In other words, in the $p$-adic framework,
Taylor solutions of differential equations are also $q$-Taylor
solution of $q$-difference equations, and conversely.}\fi

\subsection*{A motivation : the rough idea of the confluence}

Heuristically we say that a family of $q-$difference equations
$\{\;\sigma_q(Y_q)=A(q,T)\cdot Y_q\;\}_{q\in
D^-(1,\epsilon)-\{1\}}$ (where $\sigma_q$ is the automorphism
$f(T)\mapsto f(qT)$), is confluent to the differential equation
$\delta_1(Y_q)=G(1,T)\cdot Y_q$, with $\delta_1:=T\frac{d}{dT}$,
if one has $\lim_{q\to 1}\frac{A(q,T)-1}{q-1}=G(1,T)$ and, in some
suitable meaning
\begin{equation}\label{tends to Y_1,h}
\lim_{q\to 1} Y_q = Y_1\;.
\end{equation}
Roughly speaking, in this paper we show that in the $p$-adic
framework, if a differential equation is given, then, for
$\epsilon$ sufficiently small, one may choose the family
$\{G(q,T)\}_{q}$ in order to have $Y_q=Y_1$, for all
$q\in\mathrm{D}^-(1,\epsilon)$. Conversely if $q_{_0}$ is not a
root of unity, and if a single equation
$\sigma_{q_{_0}}(Y_{q_{_0}})=A(q_{_0},T)\cdot Y_{q_{_0}}$ is
given, then, under some assumptions on the radius of convergence
of its \emph{generic} Taylor solution $Y_{q_{_0}}$, one can find a
differential equation, and family as above with the property that
$Y_q=Y_{q_{_0}}=Y_1$, for all $q\in\mathrm{D}^+(1,|q_{_0}-1|)$. In
this sense, in the $p$-adic context, the solutions of
$q$-difference equations are not simply a
``\emph{discretization}'' of the solutions of differential
equations, but they are actually equal. We want now to state these
facts more precisely.

\subsection*{The work of Y.André and L.Di Vizio}
In \cite{An-DV} the authors initiated the study of the phenomena
of confluence in a  $p-$adic setting.
  For $K$  a complete discrete valuation field of mixed
characteristic, they found an equivalence between the category of
$q-$difference equations with Frobenius structure over the Robba
ring $\R_{K^{\mathrm{alg}}}$ (here called
$\sigma_q-\Mod(\R_{K^{\mathrm{alg}}})^{(\phi)}$), and the category
of differential equations with Frobenius structure over the Robba
ring $\R_{K^{\mathrm{alg}}}$ (here called
$\delta_1-\Mod(\R_{K^{\mathrm{alg}}})^{(\phi)}$). %

One of the restrictions of \cite{An-DV} is that the number $q$ is
required to satisfy $|q-1|<|p|^{\frac{1}{p-1}}$. Indeed, in the
annulus $|q-1|=|p|^{\frac{1}{p-1}}$ one encounters the $p-$th root
of unity and, if $\xi^p=1$, then the category
$\sigma_\xi-\Mod(\R_{K^{\mathrm{alg}}})^{(\phi)}$ is different in
nature from the category of differential equation, since it is not
$K^{\mathrm{alg}}-$linear.

The equivalence of \cite{An-DV} is obtained as follows. In
\cite{An} one proves that the Tannakian group of
$\delta_1-\Mod(\R_{K^{\mathrm{alg}}})^{(\phi)}$ is
$\mathcal{I}_{k^{\mathrm{alg}}(\!(t)\!)}\times\mathbb{G}_a$, where
$k$ is the perfect residue field of $K$, and
$\mathcal{I}_{k^{\mathrm{alg}}(\!(t)\!)}$ is the absolute Galois
group of $k^{\mathrm{alg}}(\!(t)\!)$. On the other hand in
\cite{An-DV} one shows that
$\sigma_q-\Mod(\R_{K^{\mathrm{alg}}})^{(\phi)}$ has the same
Tannakian group
$\mathcal{I}_{k^{\mathrm{alg}}(\!(t)\!)}\times\mathbb{G}_a$. By
composition with the respective Tannakian equivalences ($T_q$ and
$T_1$ below), one obtains then the so called the \emph{confluence
functor} `` $\mathrm{Conf}_q$ '' (in the notations of \cite{An-DV}
one has $T_1=\V_d^{(\phi)}$ and $T_q=\V_{\sigma_q}^{(\phi)}$):
\begin{equation}
\xymatrix{
\sigma_q-\Mod(\R_{K^{\mathrm{alg}}})^{(\phi)}\ar@{..>}[rr]^{\mathrm{Conf}_q}_{\cong}\ar[dr]_{T_q}^{\cong}&&
\delta_1-\Mod(\R_{K^{\mathrm{alg}}})^{(\phi)}\ar[dl]^{T_1}_{\cong}\\
&\Rep_{K^{\mathrm{alg}}}(\mathcal{I}_{k^{\mathrm{alg}}(\!(t)\!)}\times
\mathbb{G}_a)& }
\end{equation}

The strategy of \cite{An-DV} consists in showing that, as in the
case of differential equations (cf. \cite{An}), every object $\M$
in $\sigma_q-\Mod(\R_{K^{\mathrm{alg}}})^{(\phi)}$ is
quasi-unipotent, i.e. becomes unipotent after a special extension
of $\R_K$ (cf. section \ref{Special extensions}). Once a basis of
$\M$ is fixed, this means that $\M$ admits a complete basis of
solutions $\widetilde{Y} \in
\mathrm{GL}_n(\widetilde{\mathcal{R}_K}[\log(T)])$, where
$\widetilde{\R_K}$ is the union of all special extensions of
$\R_K$ (it is a sort of lifting of $k(\!(t)\!)^{\mathrm{alg}}$).
We will call ``\emph{étale}'' solutions the solutions of $\M$ in
$\widetilde{\R_K}[\log(T)]$. The proof of this relevant result
needs a substantial effort, and is actually not less complicated
than the classical $p-$adic local monodromy theorem for
differential equations itself (i.e. the fact that $T_1$ is an
equivalence).
Thanks to the fact that this important, but also very peculiar,
class of $q$-difference and differential equations are trivialized
by $\widetilde{\R_K}[\log(T)]$, one can define the functor $T_1$
(resp. $T_q$) associating to a differential (resp. $q$-difference)
equation $(\M,\delta_1^\M)$ (resp. $(\M,\sigma_q^{\M})$) the
$K^{\mathrm{alg}}$-vector space $T_1(\M,\delta_1^\M)$ (resp.
$T_q(\M,\sigma_q^{\M})$) of its ``étale'' solutions in
$\widetilde{\R_K}[\log(T)]$.\footnote{Following the definition
section \ref{section constant solutions},
$\V_{d}^{(\phi)}(\M):=(\M\otimes_{\R_K}\widetilde{\R_K}[\log(T)])^{\delta_1=0}$
is actually the dual of the space of solutions
$\mathrm{Hom}^{\delta_1}_{\R_K}(\M,\widetilde{\R_K}[\log(T)])$
(resp. same remark for
$\V_{\sigma_q}^{(\phi)}(\M):=(\M\otimes_{\R_K}\widetilde{\R_K}[\log(T)])^{\sigma_q=\mathrm{Id}}$
and
$\mathrm{Hom}^{\sigma_q}_{\R_K}(\M,\widetilde{\R_K}[\log(T)])$).}
The action of $\mathcal{I}_{k^{\mathrm{alg}}(\!(t)\!)} \times
\mathbb{G}_a$ on the space of the ``étale'' solutions arises from
its action on $\widetilde{\R_K}[\log(T)]$ by $\R_K$-linear
automorphisms commuting with $\delta_1$ and $\sigma_q$ on
$\widetilde{\R_K}[\log(T)]$.

Hence one sees for the first time in \cite{An-DV} the fact that
the ``étale'' solutions of a $q$-difference equation with
Frobenius structure, are also the ``étale'' solutions of a
differential equation. Moreover the functor $\Conf_q$ is nothing
but the functor sending a $q$-difference equation (with (strong)
Frobenius structure) into the differential equation having the
same solutions.

In the present paper we prove that this ``permanence'' of the
solutions holds also for \emph{Taylor solutions} (see below). We
develop then a $p$-adic theory of Confluence using, as a unique
tool, this fact, here called \emph{propagation principle}. We
prove indeed that this principle is sufficient to define the
Confluence and Deformation equivalences, over almost all $p$-adic
ring of functions, with very basic assumptions on the equations.
This theory requires only the definition and the formal properties
of the \emph{generic Taylor solution $Y(x,y)$}. For this reason it
is not a consequence of the heretofore developed theory.
Conversely we deduce, as a special case, the confluence of
\cite{An-DV} by comparing Taylor solutions and ``étale'' solutions
(cf. the end of the introduction).

\subsection*{The generic $q$-Taylor solution}

Let now $K$ be an arbitrary ultrametric complete valued field of
mixed characteristic $(0,p)$. Let
$X=\mathrm{D}^{+}(c_0,R_0)-\cup_{i=1,\ldots,n}\mathrm{D}^-(c_i,R_i)$
be an affinoid, where $\mathrm{D}^-(c,R)$ denotes the open disc
centered at $c$ of radius $R$. Let $\H_K(X)$ be the ring of
analytic elements on $X$. Consider a $q$-difference equation
\begin{equation}\label{equation 0.0.1}
\sigma_q(Y)=A(q,T)\cdot Y\;, \quad A(q,T)\in GL_n(\H_K(X))
\end{equation}
on $X$. Denote by $(\M,\sigma_q^{\M})$ the $q$-difference module
over $X$ defined by this equation.

A major difference between the complex and the $p$-adic settings
is that in the latter there are disks (not centered at $0$) which
are $q$-invariant. A disk $\mathrm{D}^-(c,R)\subset X(K)$ is
$q$-invariant (i.e. the map $x\mapsto qx$ is a bijection of
$\mathrm{D}^-(c,R)$) if and only if $|q-1||c|<R$, and $|q|=1$ (cf.
Lemma \ref{|q-1||c|<R-yt}). Starting from this consideration, in
\cite{DV-Dwork} the author define, for $q$-difference equations,
the $q$-analogue of the generic Taylor solution of a differential
equation (cf. Def. \ref{Definition of the q-Taylor solution}):
\begin{equation}
Y(x,y):=\sum_{n\geq 0}H_n(q,T)\frac{(x-y)_{q,n}}{[n]^!_q}\;,
\end{equation}
where $H_n(q,T)$ is obtained by iterating the equation
\eqref{equation 0.0.1}: $d_q^n(Y) = H_n(q,T)\cdot Y$, where
$d_q:=\frac{\sigma_q-1}{(q-1)T}$. For a large class of equations 
it happens that, for all $c\in X(K)$, the series $Y(x,c)$
represents a function which converges on a disk
$\mathrm{D}^-(c,R)$, with $|q-1||c|<R$. More precisely $Y(x,y)$
converges in a neighborhood of the diagonal of the type
$\mathcal{U}_R := \{(x,y)\in X \times X\;|\; |x-y| < R \}$, with
\begin{equation}
|q-1|\cdot\mathfrak{s}_X<R\;, %
\end{equation}%
where $\mathfrak{s}_X:=\sup_{c\in X}|c|$ as shown in the
following picture (one easily sees that
$\mathfrak{s}_X=\max(|c_0|,R_0)$):
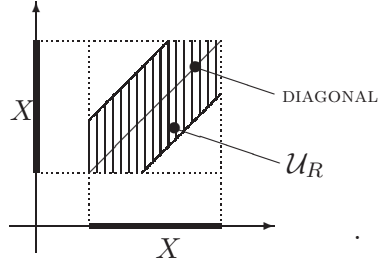
\begin{center}
\begin{picture}(100,90)
\put(10,0){\vector(0,1){95}} %
\put(0,10){\vector(1,0){100}} %
\put(10,30){\linethickness{2pt}\line(0,1){50}} %
\put(0,50){$X$} %
\put(55,-2){$X$} %
\put(30,10){\linethickness{2pt}\line(1,0){50}} %
\put(30,30){\line(1,1){50}}
\qbezier[35](10,30)(45,30)(80,30)%
\qbezier[35](10,80)(45,80)(80,80)%
\qbezier[35](30,10)(30,45)(30,80)%
\qbezier[35](80,10)(80,45)(80,80)%
\qbezier(30,50)(45,65)(60,80)%
\qbezier(50,30)(65,45)(80,60)%

\put(30,30){\line(0,1){20}}%
\put(33,30){\line(0,1){23}}%
\put(36,30){\line(0,1){26}}%
\put(39,30){\line(0,1){29}}%
\put(42,30){\line(0,1){32}}%
\put(45,30){\line(0,1){35}}%
\put(48,30){\line(0,1){38}}%
\put(51,31){\line(0,1){40}}%
\put(54,34){\line(0,1){40}}%
\put(57,37){\line(0,1){40}}%
\put(60,40){\line(0,1){40}}%
\put(63,43){\line(0,1){37}}%
\put(66,46){\line(0,1){34}}%
\put(69,49){\line(0,1){31}}%
\put(72,52){\line(0,1){28}}%
\put(75,55){\line(0,1){25}}%
\put(78,58){\line(0,1){22}}%
\put(59.5,44.5){$\bullet$}%
\put(67.5,67.5){$\bullet$}%
\put(105,30){$\mathcal{U}_R$}%
\put(103,57){\begin{scriptsize}\textsc{diagonal}\end{scriptsize}}%
\put(62,47){\line(3,-1){40}}%
\put(70,70){\line(3,-1){30}}%
\put(130,5){$.$}
\end{picture}
\end{center}

We call such equations \emph{Taylor admissible}. The matrix
function $Y(x,y):\mathcal{U}_R \to GL_{n}(K)$ is invertible and
satisfies the cocycle conditions: $Y(x,y)\cdot Y(y,z)=Y(x,z)$ and
$Y(x,y)^{-1}=Y(y,x)$, for all $(x,y),(y,z),(x,z)\in
\mathcal{U}_R$. Moreover $Y(qx,y)=A(q,x)Y(x,y)$ and, for all $c\in
X(K)$, the matrix $Y(x,c)\in GL_n(\a_K(c,R))$ is a fundamental
basis of solutions of the equation \eqref{equation 0.0.1}. In
particular the $q$-difference algebra $\a_K(c,R)$ of analytic
functions over
the disk $\mathrm{D}^-(c,R)$, trivializes $(\M,\sigma_q^{\M})$.

The following fact is the main point of this paper (cf. Theorem
\ref{main theorem second form}). If now $q'\neq q$ belongs to the
disk
$\mathrm{D}^-(q,R/\mathfrak{s}_X)=\mathrm{D}^-(1,R/\mathfrak{s}_X)$,
then the matrix
\begin{equation}
A(q',x):=Y(q'x,y)\cdot Y(x,y)^{-1}=Y(q'x,y)\cdot Y(y,x)=Y(q'x,x)
\end{equation}
\emph{is an analytic function of $x$ on all of $X$}. Indeed
$(q'x,x)\in\mathcal{U}_R$, for all $x\in X$, and hence the matrix
$A(q',x)$ maps $x\mapsto (q'x,x)\mapsto Y(q'x,x)=A(q',x)$. One
shows easily that $A(q',x)\in GL_n(\H_K(X))$, for all
$q'\in\mathrm{D}^-(1,R/\mathfrak{s}_X)$, since $Y(x,y)$ is
invertible. This fact implies that \emph{$Y(x,y)$ is
simultaneously the Taylor solution of every equation of the family
$\{\sigma_{q'}(Y)=A(q',T)Y\}_{q'}$, for all
$q'\in\mathrm{D}^-(1,R/\mathfrak{s}_X)$}. Equivalently, this means
that the $q$-difference module $(\M,\sigma_q^{\M})$ is canonically
endowed with an action of $\sigma_{q'}$, for all $q'\in
\mathrm{D}^-(1,R/\mathfrak{s}_X)$. This remarkable fact will be
called \emph{Propagation Principle}.
As one can see, this happens actually under the following weak 
assumptions on $(\M,\sigma_q^{\M})$:
\begin{eqnarray}
\textrm{i)}&\quad&\textrm{$q$ is not a root of unity;}\label{property i)}\\
\textrm{ii)}&\quad&\textrm{$Y(x,y)$ converges on some
$\mathcal{U}_R$ with $|q-1|\cdot\mathfrak{s}_X < R \leq
r_X$;\qquad\qquad\qquad\qquad\qquad\qquad\qquad\qquad}\label{property
ii)}
\end{eqnarray}
\if{\begin{enumerate}
\item $q$ is not a root of unity;%
\item $Y(x,y)$ converges on some $\mathcal{U}_R$ with
$|q-1|\cdot\mathfrak{s}_X < R \leq r_X$,
where $r_X=\min(R_0,R_1,\ldots,R_n)$ is a number depending on the
geometry of $X$.
\end{enumerate}}\fi
where $r_X=\min(R_0,R_1,\ldots,R_n)$ is a number depending on the
geometry of $X$. The category of $q$-difference modules
$(\M,\sigma_q^{\M})$ satisfying these two properties for a
suitable unspecified $R$ satisfying $|q-1|\mathfrak{s}_X<r\leq
R\leq r_X$ will be denoted by
$\sigma_q-\Mod(\H_K(X))^{[r]}$.

The assumption $|q-1|\mathfrak{s}_X<R$ assures that the image of
the map $x\mapsto(qx,x):X\mapsto X\times X$ is contained in
$\mathcal{U}_R$. The bound $R\leq r_X$ assures that the function
$Y(x,y)$ does not converge outside $X$. Indeed the properties of
$Y(x,y)$ outside $X$ are not invariant under $\H_K(X)$-base
changes in $\M$. Finally condition $\mathrm{ii)}$ also assures
that the map $x\mapsto qx$ is a bijection of $X$
 globally fixing each individual hole of $X$ (cf. section \ref{k_0}).
Since $r_X\leq\mathfrak{s}_X$, we are assuming implicitly that
$|q-1|<1$. But no restrictive assumptions on $X$ or on $K$ are
made.

Obviously this process works just as well if the initial function
$Y(x,y)$ is the generic Taylor solution of a differential
equation. The category of \emph{differential} equations whose
Taylor solution converges on $\mathcal{U}_R$, for an 
unspecified $R$ satisfying $r\leq R\leq r_{X}$,
will be denoted by $
\delta_1-\Mod(\H_K(X))^{[r]}$.

\subsection*{Discrete and analytic $\sigma$-modules}
Let $\Q(X)$ be the set of $q\in K$ for which $x\mapsto qx$ is a
bijection of $X$. Then $\Q(X)$ is a topological subgroup of
$K^{\times}$, and the disk $\mathrm{D}^-(1,R/\mathfrak{s}_X)$,
with $R\leq r_X$, is an open subgroup of $\Q(X)$. The group
$\Q(X)$ acts continuously on $\H_K(X)$ via $q\mapsto\sigma_q$. The
data of $\M$, together with the simultaneous
$\sigma_q$-semi-linear action of $\sigma_q^{\M}$, for all $q\in
\mathrm{D}^-(1,R/\mathfrak{s}_X)$, is then a \emph{semi-linear
representation of the sub-group
$\mathrm{D}^-(1,R/\mathfrak{s}_X)\subseteq\Q(X)$.} This
representation has three remarkable properties:
\begin{enumerate}
\item[(a)] The map $(q',x)\mapsto A(q',x)$ is analytic in
$(q',x)$. In
particular the representation is continuous; %
\item[(b)] The group $\mathrm{D}^-(1,R/\mathfrak{s}_X)$ depends on
$R$,
and hence on $\M$;  %
\item[(c)] The matrix $Y(x,y)$ is simultaneously the generic
Taylor solution of the $q$-difference module $(\M,\sigma_q^{\M})$,
for all $q\in \mathrm{D}^-( 1 , R/\mathfrak{s}_X )$.
\end{enumerate}

Inspired by the first two properties we define a new class of
objects called \emph{discrete or analytic $\sigma$-modules} as
follows. Consider a subset $S\subset\Q(X)$. A \emph{discrete
$\sigma$-module} on $S$ is nothing but a $\H_K(X)$ semi-linear
representation of the group $\ph{S}$ generated by $S$. If $S=U$ is
an open subset of  $\Q(X)$, we define \emph{analytic
$\sigma$-modules on $U$} to be a discrete $\sigma$-modules over
$U$ together with a certain condition of analyticity of
$\sigma^{\M}_q$ with respect to $q$. These categories are denoted
by $\sigma-\Mod(\H_K(X))_S^{\mathrm{disc}}$ and
$\sigma-\Mod(\H_K(X))_U^{\mathrm{an}}$ respectively. In this paper
the words ``\emph{discrete}'' or ``\emph{analytic}'' will be
referred to the discreteness or analyticity of $\sigma_q^{\M}$
with respect to $q$. We heuristically imagine the analytic
$\sigma$-modules as \emph{semi-linear representations of the
(co-variant) sheaf of groups} $U\mapsto \ph{U}$.

\begin{remark}\label{Morphisms between sigma modules}
It is important to notice that morphisms between analytic
$\sigma$-modules over $U$ are morphisms of representations. More
precisely once a basis of $\M$ (resp. $\N$) is fixed, we have a
family of operators $\{\sigma_q(Y)=A(q,T)Y\}_{q\in \ph{U}}$ (resp.
$\{\sigma_q(Y)=\widetilde{A}(q,T)Y\}_{q\in \ph{U}}$) such that
$A(q,T)$ (resp. $\widetilde{A}(q,T)$) depends analytically on
$(q,T)$.\footnote{The data of an analytic $\sigma$-module is
actually nothing but ``\emph{a family of $q$-difference equations
depending analytically on $q$}''.} A morphism $\alpha:\M\to\N$
then must simultaneously commute with $\sigma^{\M}_q$ and
$\sigma^{\N}_q$, for all $q\in\ph{U}$. In other words the matrix
$B$ of $\alpha$ must simultaneously verify
$A(q,T)B=\sigma_q(B)\widetilde{A}(q,T)$, for all $q\in\ph{U} $.
Actually there are \emph{non isomorphic} analytic $\sigma$-modules
over $U$ defining isomorphic $q$-difference equations at every
$q\in\ph{U}$ (see example \ref{Morphisms between sigma modules
-examples}). This is analogous to have \emph{non isomorphic}
sheaves having isomorphic stalks at every point.
\end{remark}

\subsection*{Taylor admissible $\sigma$-modules}
We now want to analyse property $(c)$: the constancy of the solutions. %
If $S\not\subseteq\bs{\mu}_{p^{\infty}}$, we call \emph{Taylor
admissible $\sigma$-modules over $S$} those $\sigma$-modules for
which the $q$-Taylor solution $Y(x,y)$ is the same for all
$q\in\ph{S}$, and satisfy the condition ii), for all $q\in S$ (cf.
\eqref{property ii)}). If $S=U$ is open, by the Propagation
Principle, Taylor admissible $\sigma$-modules are
\emph{automatically} analytic on $U$ (cf. Remark \ref{admissible
==> analytic}). This category is denoted by
$\sigma-\Mod(\H_K(X))_U^{\mathrm{adm}}\;\subseteq\;\sigma-\Mod(\H_K(X))_U^{\mathrm{an}}$.
We heuristically imagine Taylor admissible $\sigma$-modules as
\emph{semi-linear representations of the (co-variant) sheaf of
groups} $U\mapsto \ph{U}$, which are \emph{locally constant}.

Taylor admissibility is a particular case of a more classical
notion. If $\mathrm{C}/\H_K(X)$ is an algebra admitting an action
of the group $\ph{S}$ extending that on $\H_K(X)$, then a
semi-linear representation of $\ph{S}$ over $\H_K(X)$ is called
\emph{$\C$-admissible} if it is trivialized by $\C$.
%
%
%
%
For a discrete $\sigma$-module $\M$ over $S$ to be trivialized by
$\C$ means exactly that there exists $Y\in GL_n(\C)$ which is a
simultaneous solution of all operators defined by $\M$. If $\M$ is
trivialized by $\C$ we will say that $\M$ is \emph{$\C$-constant}.
We observe that if $S=q^{\mathbb{Z}}$, then $\C$ is nothing but a
$q$-difference algebra over $\H_K(X)$. So the constancy of the
solutions does not depend on the analyticity of $\M$, rather it is
a \emph{discrete} fact.

In section \ref{section - solution formal def} we define
\emph{discrete $\sigma$-algebras}, and we develop a basic
differential/difference Galois theory for discrete
$\sigma$-algebras. The analogue of the
Picard-Vessiot theorem providing the existence 
of a discrete $\sigma$-algebra trivializing a given discrete
$\sigma$-module \emph{is missing}. We are thus obliged to work
with the category of discrete $\sigma$-modules trivialized by a
fixed discrete $\sigma$-algebra $\C$.
In section \ref{theory of deformation} we develop formally the
theory of $\C$-Confluence and $\C$-Deformation, which will also
depend on the chosen discrete $\sigma$-algebra $\C$.

\begin{remark}\label{intro -constancy of solutions}
Notice that solutions will be defined formally as morphisms
$\M\to\C$ commuting simultaneously with the actions of $\sigma_q$
for all $q\in S$ (cf. Section \ref{section constant solutions}).
This fact, together with Remark \ref{Morphisms between sigma
modules}, 
explains why the notion of $\C$-\emph{constant}
$\sigma$-module implies the constancy of the solutions (with
respect to $q$).
\end{remark}

\subsection*{The Confluence functor}
Let $(\M,\sigma^{\M})$ be an analytic $\sigma$-module over $U$. By
analyticity we also have an action of the \emph{Lie algebra} of
$\ph{U}$ (here systematically identified with $K\cdot \delta_1$).
In other words the following limit converges to a connection
$\delta_1^{\M}:\M\to\M$ (cf. section \ref{analytic sigma delta
modules}):
\begin{equation}\label{intro - limit of confluence}
\delta_1^{\M}\;\;:=\;\;\lim_{\begin{smallmatrix} q\in\ph{U},q\to 1
\end{smallmatrix}}\frac{\sigma_q^{\M}-1}{q-1}\;\;\in
\;\;\mathrm{End}^{\mathrm{cont}}_K(\M)\;\;,
\end{equation}
where 
$q$ runs over the (open) group $\ph{U}$ generated by $U$. In terms
of matrices, the matrix $G(1,T)$ of $\delta_1^{\M}$ is
$G(1,T)=q\frac{\partial}{\partial q}\Bigl(A(q,T)\Bigr)_{|_{q=1}}$
(cf. equation \eqref{G(q,T)=qd/dqA(q,T)}). By continuity,
morphisms of analytic $\sigma$-modules also commute with the
connection (cf. remark \ref{sigma to sigma,delta is fully
faithful}).  Hence we obtain
 a functor called  $\mathrm{Conf}_U:\sigma-\Mod(\H_K(X))^{\mathrm{an}}_U\xrightarrow[]{\;\;\;\;}\delta_1-\Mod(\H_K(X))$, sending $(\M,\sigma^{\M})$ into
$(\M,\delta_1^{\M})$ (cf. Remark \ref{Definition of Conf_U}). This
functor is not an equivalence, but it does induce an equivalence:
\begin{equation}
\mathrm{Conf}_U^{\mathrm{Tay}}\;:\;\sigma-\Mod(\H_K(X))^{[r]}_U
\;\; \xrightarrow[]{\;\;\sim\;\;} \;\;
\delta_1-\Mod(\H_K(X))^{[r]}\;,
\end{equation}
where $\mathrm{Conf}_U^{\mathrm{Tay}}$ simply denotes the
restriction of $\mathrm{Conf}_U$ to the category
$\sigma-\Mod(\H_K(X))^{[r]}_U\subseteq
\sigma-\Mod(\H_K(X))^{\mathrm{adm}}_U$ of Taylor admissible
$\sigma$-modules verifying condition ii) with $r\leq R\leq r_{X}$
(cf. \eqref{property ii)}), where $r>0$ is large enough to have
$U\subset\mathrm{D}^-(1,r/\mathfrak{s}_X)$ (cf. Corollary
\ref{compare with}). The Propagation Principle gives a quasi
inverse functor (cf. Remark \ref{Definition of Conf_U} for a
formal presentation).

On the other hand let $q\in U-\bs{\mu}_{p^{\infty}}$. An analytic
$\sigma$-module over $U$ defines a $q$-difference module by
forgetting the action of $\sigma_{q'}^{\M}$, for all $q'\neq q$.
Again the Propagation Principle provides an equivalence
\begin{equation}
\mathrm{Res}^{U}_{q}\;:\;\sigma-\Mod(\H_K(X))^{[r]}_{U}\;\;
\xrightarrow[]{\;\;\sim\;\;}\;\; \sigma_q-\Mod(\H_K(X))^{[r]}\;,
\end{equation}
where $r\leq r_X$ is sufficiently large to have
$U\subseteq\mathrm{D}^-(1,r/\mathfrak{s}_X)$ (cf. Cor.
\ref{compare with}). We call the composite equivalence
$\mathrm{Conf}_q^{\mathrm{Tay}}$.  Thus we have
\begin{equation}\label{intro - Conf_q}
\mathrm{Conf}_q^{\mathrm{Tay}}:=\mathrm{Conf}_U^{\mathrm{Tay}}\circ(\mathrm{Res}^U_q)^{-1}\;:\;
\sigma_q-\Mod(\H_K(X))^{[r]}\xrightarrow[]{\;\;\sim\;\;}\delta_1-\Mod(\H_K(X))^{[r]}\;.
\end{equation}
The equivalence $\mathrm{Conf}_q^{\mathrm{Tay}}$ sends a
$q$-difference equation satisfying conditions i) and ii) (cf.
\eqref{property i)},\eqref{property ii)}), into the differential
equation having the same generic Taylor solution.

\subsection*{Roots of unity and $q$-tangent operators}
In this last equivalence the number $q$ must not belong to
$\bs{\mu}_{p^{\infty}}$. If $q'=\xi$, with $\xi^{p^n}=1$, the
category of $\sigma_{\xi}$-difference equations is not $K$-linear
and cannot be equivalent to the category of differential
equations. Nevertheless, if, for $q\notin\bs{\mu}_{p^{\infty}}$,
the radius $R$ of the $q$-Taylor solution is large, the
Propagation Principle gives an operator
$\sigma_{\xi}^{\M}:\M\to\M$ acting on $\M$. The idea is to replace
the category $\sigma_{\xi}-\Mod(\H_K(X))$ with another category.
The expected object ``at $\xi$'' should also be endowed with an
action of the Lie algebra, \emph{as we have just done in the case
$\xi=1$}. For all $q\in \ph{U}$ the action of the Lie algebra of
$\ph{U}$ is given by the limit
$\delta_q^{\M}:=\lim_{\begin{smallmatrix}q'\to
q\end{smallmatrix}}\frac{\sigma_{q'}^{\M}-\sigma_q^{\M}}{q'-q}\in\mathrm{End}^{\mathrm{cont}}_K(\M)
$, for $q,q'\in\ph{U}$, as shown in the picture:
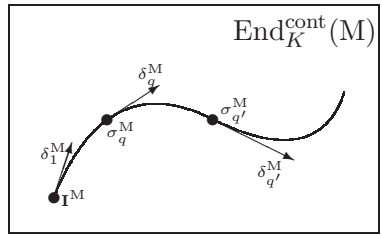
\begin{center}
\framebox{
\begin{picture}(130,80)
\put(78,70){$\mathrm{End}_K^{\mathrm{cont}}(\M)$}

\qbezier(10,10)(30,60)(70,40) %
\qbezier(70,40)(110,20)(120,50) %


\put(7.5,7.5){$\bullet$}%
\put(13,7){\begin{tiny}$\bs{\mathrm{I}}^{\mathrm{M}}$\end{tiny}}


\put(10,10){\vector(1,3){7}}%
\put(5,25){\begin{tiny}$\delta_1^{\mathrm{M}}$\end{tiny}} 




\put(27.5,37){$\bullet$}%
\put(30,33){\begin{tiny}$\sigma_q^{\mathrm{M}}$\end{tiny}}


\put(30,39.5){\vector(3,2){20}}%
\put(42,55){\begin{tiny}$\delta_q^{\mathrm{M}}$\end{tiny}} 


\put(67.5,37){$\bullet$} %


\put(73,42){\begin{tiny}$\sigma_{q'}^{\mathrm{M}}$\end{tiny}}


\put(71,39.5){\vector(2,-1){30}}%
\put(87,18){\begin{tiny}$\delta_{q'}^{\mathrm{M}}$\end{tiny}} %
\end{picture}}
\end{center}
Clearly $\delta_q^{\M}=\sigma_q^{\M}\circ\delta_1^{\M}$, so to
give $\delta_q^{\M}$ is equivalent to give $\delta_1^{\M}$. In a
root of unity the ``limit object'' is a mixed data
$(\M,\sigma_{\xi}^{\M},\delta_{1}^{\M})$, i.e. a connection
$\delta_1^{\M}$ on $\M$ together with an action of
$\sigma_{\xi}^{\M}$ on $\M$. We call these new objects
 \emph{$(\sigma_{\xi},\delta_{\xi})$-modules}. In the sequel every
 terminology is given simultaneously for $\sigma$-modules and
$(\sigma,\delta)$-modules. The additional data of
$\delta_{\xi}^{\M}$ makes the category of
$(\sigma_{\xi},\delta_{\xi})$-modules $K$-linear. Moreover
$\delta_{\xi}^{\M}$ preserve the ``\emph{information}'' in a
neighborhood of $\xi$, indeed we find equivalences
\begin{eqnarray}\label{intro - Conf_xi}
\mathrm{Conf}_{\xi}^{\mathrm{Tay}}:=\mathrm{Conf}_U^{\mathrm{Tay}}\circ(\mathrm{Res}^U_{\xi})^{-1}\;:\;
(\sigma_{\xi},\delta_{\xi})-\Mod(\H_K(X))^{[r]}
&\xrightarrow[]{\;\;\sim\;\;}& \delta_1-\Mod(\H_K(X))^{[r]}\;,\\
\mathrm{Def}_{\xi,q}^{\mathrm{Tay}}:=\mathrm{Res}^U_{q}\circ(\mathrm{Res}^U_{\xi})^{-1}\;:\;
(\sigma_{\xi},\delta_{\xi})-\Mod(\H_K(X))^{[r]}
&\xrightarrow[]{\;\;\sim\;\;}&
(\sigma_q,\delta_q)-\Mod(\H_K(X))^{[r]}\;.\quad\qquad
\end{eqnarray}
If $q$ is not a root of unity, then the data of $\delta_1^{\M}$ is
superfluous, indeed if the module is Taylor admissible the
Propagation Principle allows one to re-construct $\delta_1^{\M}$ from
$\sigma_q^{\M}$.

In the classical setting over the complex numbers $\mathbb{C}$,
understanding of the case $q=\xi\in\bs{\mu}_{p^\infty}$ remains an
open problem.

\subsection*{Quasi unipotence and comparison with Andr\'e-Di Vizio's Confluence}
Up to a correct definition for the notion of Taylor admissibility,
the previous theory can be generalized to more general rings of
functions. From section \ref{Propagation Theorem for other rings}
on we obtain the theory over $\R_K$. We prove that every
$q$-difference equations with Frobenius Structure over $\R_K$, is
quasi unipotent (i.e. is trivialized by
$\widetilde{\R_K}[\log(T)]$), for all
$q\in\mathrm{D}^-(1,1)-\bs{\mu}_{p^{\infty}}$, generalizing the
main result of \cite{An-DV}. We actually prove this theorem in the
more general context of $\sigma$-modules, and
$(\sigma,\delta)$-modules. We deduce it by the quasi unipotence of
$p$-adic \emph{differential} equations with Frobenius Structure
over $\R_K$, and by deformation. The idea is the following. As
already mentioned, we are obliged to work with $\sigma$-modules
trivialized by a fixed discrete $\sigma$-algebra $\C$, and the
$\C$-Confluence and $\C$-Deformations functors depend on $\C$. In
the ``quasi unipotent'' context  this algebra is
$\C:=\widetilde{\R_K}[\log(T)]$, while in the context of the
propagation theorem $\C:=\a_K(c,R)$, for an arbitrary point $c\in
X$, and suitable $R>0$. To compare Taylor solutions to the ``étale
solutions'' in $GL_n(\widetilde{\R_K}[\log(T)])$, the idea is to
find a discrete $\sigma$-algebra of \emph{functions over a disk}
containing $\widetilde{\R_K}[\log(T)]$. Actually such an algebra
does not exist.  Thus we use  a theorem of S.Matsuda (cf. Th.
\ref{canonical extension}) providing an equivalence between
$\delta_1-\Mod(\R_K)^{(\phi)}$ with the sub-category of
$\delta_1-\Mod(\Hd_K)^{(\phi)}$ formed by Special objects. Special
objects are trivialized by a special extension of $\Hd_K$ (cf.
Section \ref{Special extensions}). The ring $\a_K(1,1)$ is a
discrete $\sigma$-algebra over $\H_K^{\dag}$. We prove then that
the algebra $\widetilde{\C^{\textrm{ét}}_K}[\log(T)]$ generated
over $\Hd_K$ by all the ``étale solutions'' of Special objects
admits an embedding
$\widetilde{\C^{\textrm{ét}}_K}[\log(T)]\subset\a_{K^{\mathrm{alg}}}(1,1)$
commuting with $\delta_1$, with the Frobenius, and with
$\sigma_q^{\M}$, for all
$q\in\mathrm{D}^-(1,1)-\bs{\mu}_{p^{\infty}}$ (cf. Lemma
\ref{embedding in a_K(1,1)}). This will prove that the
$\C$-Confluence and the $\C$-Deformation functors defined by using
$\C=\a_K(1,1)$, or $\C = \widetilde{\R_K}[\log(T)]$ are actually
the same (cf. Cor. \ref{Conf Tay = Conf a_K(1,1) = Conf widetilde
R_K}). Moreover it proves also that the confluence of André-Di
Vizio coincides with our $\Conf_q^{\mathrm{Tay}}$ (cf. Section
\ref{Andre e Di Vizio}), thus it is independent on the Frobenius.

\subsection*{Structure of the paper}
Section \ref{section - notations} is devoted to notation. In
section \ref{section - Discrete or analytic sigma and sigma,delta
mod}, we give definitions and basic facts on
\emph{discrete/analytic $\sigma-$modules, and
$(\sigma,\delta)-$modules}. In section \ref{section - solution
formal def} we define \emph{discrete $\sigma$-algebras and
$(\sigma,\delta)$-algebras}, and we give the abstract definition
of \emph{solutions}. In section
\ref{theory of deformation} we give the formal notion of 
 \emph{confluence}. In section \ref{Taylor solutions--} we
introduce \emph{generic Taylor solutions} and \emph{generic radius
of convergence}. In section \ref{section - main theorem} we define
\emph{Taylor admissible objects} and obtain the main
\emph{Propagation Theorem} \ref{main theorem second form}. In the
last section \ref{section - quasi unipotence and plmt} we apply
the previous theory to the Robba ring, and to the \emph{$p-$adic
local monodromy theorem}.
\vspace{-0.3cm}
\begin{acknowledgements}
We thank L. Di Vizio who has always been willing to talk and to
explain the technical difficulties of her papers. The author
wishes to express also his gratitude to G. Christol for useful
advice and constant encouragements. A particular thanks goes to
B.Chiarellotto for his support and for generic corrections, and
also to F.Sullivan for English corrections. We thank moreover Y.
Andr\'e, J.P. Ramis, and J.Sauloy for useful discussions. The last
drafting of the paper have been highly improved thanks to some
remarks and advises of Y.André, L.Di Vizio, and the referee. Any
remaining inaccuracies are entirely my fault.

\end{acknowledgements}

\vspace{-0.5cm}
\begin{center}
\textbf{Index of categories}
\end{center}
\addcontentsline{toc}{section}{Index of categories}
\vspace{-0.3cm}
\begin{multicols}{3}{
\begin{picture}(10,0)\put(-16,0){$\sigma-\mathrm{Mod}(\B)_{S}^{\mathrm{disc}}$}\end{picture}
\hfill \pageref{definition of discrete sigma} \\
$\sigma_q-\mathrm{Mod}(\B)$\hfill\pageref{sigma_q-Mod =
sigma-Mod^disc_q} \\
$(\sigma,\delta)-\mathrm{Mod}(\B)_{S}^{\mathrm{disc}}$
\hfill\pageref{(sigma,delta)-Mod(B)_S^disc} \\
$(\sigma_q,\delta_q)-\mathrm{Mod}(\B)$\hfill\pageref{sigma_q,delta_q-Mod
= sigma,delta-Mod^disc_q} \\
$\delta_1-\Mod(\B)$\hfill\pageref{sigma_q,delta_q-Mod =
sigma,delta-Mod^disc_q} \\
$\sigma-\mathrm{Mod}(\B)_U^{\mathrm{an}}$ \hfill\pageref{fgrt}\\
$\sigma-\mathrm{Mod}(\R_K)_U^{\mathrm{an}}$\hfill\pageref{def of
an sigma mod over R_K and Hd_K}\\
$\sigma-\mathrm{Mod}(\R_K)_S^{\mathrm{disc}}$\hfill\pageref{def of
an sigma mod over R_K and Hd_K}\\
$\sigma-\mathrm{Mod}(\Hd_K)_U^{\mathrm{an}}$\hfill\pageref{def of
an sigma mod over R_K and Hd_K}\\
$\sigma-\mathrm{Mod}(\Hd_K)_S^{\mathrm{disc}}$\hfill\pageref{def
of an sigma mod over R_K and Hd_K}\\
$(\sigma,\delta)-\mathrm{Mod}(\B)_U^{\mathrm{an}}$\hfill\pageref{definition
of analytic sigma delta modules;}\\
$\sigma-\mathrm{Mod}(\B,\mathrm{C})^{\mathrm{const}}_S$\hfill\pageref{hihi}\\
$(\sigma,\delta)-\mathrm{Mod}(\B,\mathrm{C})^{\mathrm{const}}_S$\hfill\pageref{hihi}\\
$\sigma-\mathrm{Mod}(\B,\mathrm{C})^{\mathrm{an,const}}_U$\hfill\pageref{hihi}\\
$(\sigma,\delta)-\mathrm{Mod}(\B,\mathrm{C})^{\mathrm{an,const}}_U$\hfill\pageref{hihi}\\
$\sigma_q-\mathrm{Mod}(\B,\mathrm{C})_{S}$
\hfill\pageref{siiisiiisi}\\
$\sigma_q-\Mod(B,\C)^{\mathrm{an}}_U$
\hfill\pageref{siiisiiisi}\\
$\sigma-\Mod(\R_K)^{[r]}_S$ \hfill\pageref{sigma-Mod(R)^[r]_S}\\
$\sigma-\Mod(\Hd_K)^{[r]}_S$
\hfill\pageref{sigma-Mod(R)^[r]_S}\\
$(\sigma,\delta)-\Mod(\R_K)^{[r]}_S$\hfill\pageref{truy}\\
$(\sigma,\delta)-\Mod(\Hd_K)^{[r]}_S$\hfill\pageref{truy}\\
$\sigma-\Mod(\H_K(X))_S^{[r]}$\hfill\pageref{Remark - ffff}\\
$\sigma-\Mod(\H_K(X))_S^{\mathrm{adm}}$\hfill\pageref{Remark -
ffff}\\
$(\sigma,\delta)-\Mod(\H_K(X))_S^{[r]}$\hfill\pageref{Remark -
ffff}\\
$(\sigma,\delta)-\Mod(\H_K(X))_S^{\mathrm{adm}}$\hfill\pageref{Remark
- ffff}\\
$\sigma-\Mod(\R_K)_S^{\mathrm{adm}}$\hfill\pageref{Taylor adm over
R_K and Hd_K .....}\\
$\sigma-\Mod(\Hd_K)_S^{\mathrm{adm}}$\hfill\pageref{Taylor adm
over R_K and Hd_K .....}\\
$\sigma-\mathrm{Mod}(\R_K)_S^{(\phi)}$\hfill\pageref{def of frob
structure - order}\\
$\sigma-\mathrm{Mod}(\Hd_K)_S^{(\phi)}$\hfill\pageref{def of frob
structure - order}\\
$\delta_1-\Mod(\Hd_K)^{\mathrm{Sp}}$\hfill\pageref{Can}\\
$\sigma-\Mod(\Hd_K)_S^{\mathrm{Sp}}$\hfill\pageref{canonical
extension for sigma and (sigma,delta) modules}\\
$(\sigma,\delta)-\Mod(\Hd_K)_S^{\mathrm{Sp}}$\hfill\pageref{canonical
extension for sigma and (sigma,delta) modules}\\ 
$\sigma_q-\Mod(\R_K)^{\mathrm{conf}(\phi)}$\hfill\pageref{sigma_q-Mod(R_K)^(conf(phi))}\\
$\sigma_q-\Mod(\R_K)^{\mathrm{conf}}$\hfill\pageref{sigma_q-Mod(R_K)^conf}
}
\end{multicols}

\section{Notation}\label{section - notations}

We refer to \cite{De-Mi} for the definitions concerning Tannakian
categories. In the sequel when we say that a given category
$\mathcal{C}$ is (or is not) $K$-linear, we mean that the ring of
endomorphisms of the unit object is (or is not) exactly
\emph{equal} to $K$. We set $\mathbb{R}_\geq :=
\{r\in\mathbb{R}\;|\; r\geq 0\}$, and $\delta_1:=T\frac{d}{dT}$.

\subsection{Rings of functions}\label{rings}Let $R>0$ and $c\in
K$. The ring of analytic functions on the disk $\mathrm{D}^-(c,R)$
is
\begin{equation}\label{A_K(c,R)}\index{A_K(c,R)@$\a_{K}(c,R)$}
 \a_{K}(c,R):=\{\sum_{n\geq
0}a_n(T-c)^n\;|\;a_n\in K, \liminf_{n}|a_n|^{-1/n}\geq R \}\;.
\end{equation}
Its topology is given by the family of norms $|\sum
a_i(T-c)^i|_{(c,\rho)}:=\sup |a_i|\rho^i$, for all $\rho<R$. Let
$\emptyset\neq I\subseteq\mathbb{R}_{\geq 0}$ be some interval. We
denote the annulus relative to $I$ by \label{C(I)}
$\mathcal{C}_K(I):=\{x\in K\;|\;|x|\in I\}$. By $\mathcal{C}(I)$,
without the index $K$, we mean the annulus itself and not its
$K-$valued points. The ring of analytic functions on
$\mathcal{C}(I)$  is
\begin{equation}\label{A_K(I)}
\a_K(I):=\{\sum_{i\in\mathbb{Z}} a_iT^i\;|\; a_i\in K, \lim_{i\to
\pm\infty}|a_i|\rho^i=0, \textrm{ for all } \rho\in I\}\;.
\end{equation}
We set $|\sum_ia_iT^i|_\rho:=\sup_i|a_i|\rho^i<+\infty$, for all
$\rho\in I$. The ring $\a_K(I)$ is complete for the topology given
by the family of norms $\{|.|_\rho\}_{\rho\in I}$. Set
$I_\varepsilon:=]1-\varepsilon,1[$, $0<\varepsilon<1$. The Robba
ring is  defined as
\begin{equation}
\R_K:=\bigcup_{\varepsilon>0}\a_K(I_\varepsilon)\;,
\end{equation}
and is complete with respect to the limit Frechet topology.

\subsection{Affinoids}\label{affinoid - section -tfrg jj}

\begin{definition}\label{A=P^1- U D^-(c_i,r_i)}
A $K-$\emph{affinoid} is an analytic subset of $\mathbb{P}^1$
defined by
\begin{equation}\index{X@$X:=\mathrm{D}^+(c_0,R_0)-\bigcup_{i=1}^n\mathrm{D}^-(c_i,R_i)$}
X:=\mathrm{D}^+(c_0,R_0)-\bigcup_{i=1}^n\mathrm{D}^-(c_i,R_i)\;,
\end{equation} for some $R_0,\ldots,R_n>0$, $c_0,\ldots,c_n\in K$,
$c_1,\ldots,c_n\in\mathrm{D}^+_K(c_0,R)$. We denote by $X$ the
$K-$affinoid itself, and for all ultrametric valued $K-$algebras
$(L,|.|)$, we denote by  $X(L)$ its $L-$rational points.
\end{definition}

Let $H^{\mathrm{rat}}_K(X)$ be the ring of rational fractions
$f(T)$ in $K(T)$, without poles in $X(K^{\mathrm{alg}})$, and let
$\|.\|_X$ be the norm on $H^{\mathrm{rat}}_K(X)$ given by
$\|f(T)\|_{X}:=\sup_{x\in X(K^\mathrm{alg})}|f(x)|$. We denote by
\begin{equation}\index{H_K(X)@$\H_K(X)$}\label{H_K(X)}
\mathcal{H}_K(X)
\end{equation}
the completion of $(H^{\mathrm{rat}}_K(X),\|\cdot\|_X)$. It is
known that if $\rho_1,\rho_2\in |K^{\mathrm{alg}}|$, and if
$X=\mathrm{D}^+(0,\rho_2)-\mathrm{D}^-(0,\rho_1)$, then
$\H_K(X)=\a_K([\rho_1,\rho_2])$. Let now $\varepsilon>0$. If
$X=\mathrm{D}^+(c_0,R_0)-\bigcup_{i=1}^n\mathrm{D}^-(c_i,R_i)$, we
set
$X_\varepsilon:=\mathrm{D}^+(c_0,R_0+\varepsilon)-\bigcup_{i=1}^n\mathrm{D}^-(c_i,R_i-\varepsilon)$.
We then set
\begin{equation}\label{HHHH}
\Hd_K(X):=\bigcup_{\varepsilon>0}\mathcal{H}_K(X_\varepsilon)\;.
\end{equation}
The ring $\Hd_K(X)$ is complete with respect to the limit
topology. Let $X_1:=\{x\;|\;|x|=1\}$, we set
\begin{equation}\index{Hd_K@$\Hd_K$}\label{Hd_K}
\mathcal{H}_K:=\mathcal{H}_K(X_1)\;,\quad \Hd_K:=\Hd_K(X_1)\;.
\end{equation}

\subsection{Norms}\label{norms}

Every semi-norm $|.|_{\B}$ on a ring $\B$ will be extended to a
semi-norm on $M_{n\times n}(\B)=M_{n}(\B)$, by setting
$|(b_{i,j})_{i,j}|_{\B}:=\max_{i,j}|b_{i,j}|_{\B}$.

\begin{definition}\label{multiplicative seminorm}
Let $X$ be an affinoid. A \emph{bounded multiplicative
semi-norm} on $\H_K(X)$ is a function
$|.|_*:\H_K(X)\to\mathbb{R}_{\geq 0}$, such that $|0|_*=0$,
$|1|_*=1$, $|f-g|_*\leq \max(|f|_*,|g|_*)$, $|fg|_*=|f|_*|g|_*$,
and $|.|_*\leq C \|.\|_X$, for some constant $C>0$.
\end{definition}

\subsubsection{}\label{exemple of norms} Let
$(L,|.|)/(K,|.|)$ be an extension of valued fields. Let $c\in
X(L)$, then $|.|_c:f\mapsto |f(c)|_{L}$ is a bounded
multiplicative semi-norm on $\H_K(X)$. If
$\mathrm{D}^+(c,R)\subseteq X$, then
$|f|_{(c,R)}:=\sup_{x\in\mathrm{D}_{L^{\mathrm{alg}}}^+(c,R)}|f(x)|$
is a bounded multiplicative semi-norm on $\H_K(X)$. Moreover if
$f=\sum_{i\geq 0} a_i(T-c)^i$, $a_i\in L$ is the Taylor expansion
of $f$ at $c\in X(L)$, then $|f|_{(c,R)}=\sup_i|a_i|R^{i}$.

\begin{definition}
Let $f(T)=\sum_{i\in\mathrm{Z}}a_{i}(T-c)^i$, $a_{i}\in K$, be a
formal power series. We set
$|f|_{(c,\rho)}:=\sup_{i}|a_{i}|\rho^i$, this number can be equal
to $+\infty$.
\end{definition}

\begin{definition}\label{log graphic ..}
Let $r\mapsto N(r):\mathbb{R}_{\geq 0}\to\mathbb{R}_{\geq 0}$ be a
function. The $\log$-function attached to $N$ is defined by
$\widetilde{N}(t):=\log(N(\exp(t)))$ (i.e.
$\widetilde{N}\;:\;\mathbb{R}\cup\{-\infty\}\;\xrightarrow[\sim]{\exp}\;
\mathbb{R}_{\geq 0}\;\xrightarrow[]{N}\; \mathbb{R}_{\geq
0}\;\xrightarrow[\sim]{\log}\;\mathbb{R}\cup\{-\infty\}$). We will
say that $N$ has  a given property logarithmically if
$\widetilde{N}$ has that property.
\end{definition}

\begin{definition}\index{Ray(f(T),c)@$Ray(f(T),c)$}
\label{radii} Let $f(T)=\sum_{i\geq 0} a_{i}(T-c)^i$, $a_{i}\in K$
be a formal power series. The radius of convergence of $f(T)$ at
$c$ is $Ray(f(T),c):=\liminf_{i\geq 0}|a_{i}|^{-1/i}$. If
$F(T)=(f_{h,k}(T))_{h,k}$, is a matrix, then we set
$Ray(F(T),c):=\min_{h,k}\;Ray(f_{h,k}(T),c)$.
\end{definition}

\begin{lemma}[(\protect{\cite[ch.II]{Ch-Ro}})]\label{log-convex property}
Let $f(T)\in K[[T-c]]$. Suppose that $|f|_{(c,\rho_0)}<\infty$,
for some $\rho_0>0$. Then:
\begin{enumerate}
\item For all $\rho < \rho_0$ one has $Ray(f(T),c) \geq \rho$, and $|f|_{(c,\rho)} < \infty$; %
\item the function
$\rho\mapsto|f|_{(c,\rho)}:[0,\rho_0]\longrightarrow\mathbb{R}_{\geq
0}$ is log-convex, piecewise log-affine, and log-increasing:
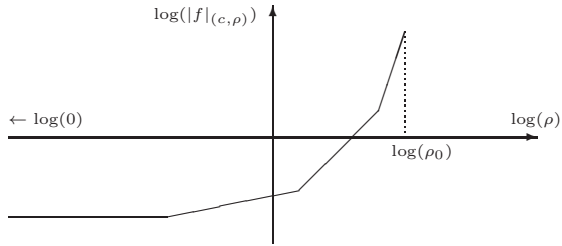
\begin{center}
\begin{picture}(100,90)
\put(50,0){\vector(0,1){90}} %
\put(-50,40){\vector(1,0){200}} %
\put(-50,10){\line(1,0){60}}%
\put(10,10){\line(5,1){50}}%
\put(60,20){\line(1,1){30}}%
\put(90,50){\line(1,3){10}}%
\qbezier[20](100,40)(100,60)(100,80)%
\put(95,33){\begin{tiny}$\log(\rho_0)$\end{tiny}} %
\put(140,45){\begin{tiny}$\log(\rho)$\end{tiny}} %
\put(5,85){\begin{tiny}$\log(|f|_{(c,\rho)})$\end{tiny}} %
\put(-50,45){\begin{tiny}$\leftarrow\log(0)$\end{tiny}} %
\end{picture}
\end{center}
\item One has $|f(T)|_{(c,\rho)}=\sup_{|x-c|\leq\rho, x\in
K^{\mathrm{alg}}}|f(x)|_{K^{\mathrm{alg}}}=\lim_{r\to\rho^{-}}\sup_{|x-c|=r,
x\in K^{\mathrm{alg}}}|f(x)|_{K^{\mathrm{alg}}}$; %
\item All zeros of $f(T)$ are algebraic. Moreover $f(T)$ has a
zero $\zeta\in K^{\mathrm{alg}}$, with $|\zeta-c|=\rho<\rho_{0}$,
if and only if the previous graph has a break at
$\log(\rho)$.\hfill\CVD
\end{enumerate}
\end{lemma}

\subsection{Generic points}
\label{Omega} Let $(\Omega,|.|)/(K,|.|)$ be a complete field such
that $|\Omega|=\mathbb{R}_{\geq 0}$, and that $k_{\Omega}/k$ is
not algebraic.

\begin{proposition}[(\protect{\cite[9.1.2]{Ch-Ro}})]
\label{definition of generic point}
\index{t_c,rho@$t_{c,\rho}$}
For every disk $\mathrm{D}^+(c,\rho)$, $c\in K$, there exists a
point $t_{c,\rho}\in \Omega$, called \emph{generic point of
$\mathrm{D}^+(c,\rho)$} such that
$|\,t_{c,\rho}-c\,|_{\Omega}=\rho$, and that
$\mathrm{D}^-_\Omega(\;t_{c,\rho}\;,\;\rho\;)\cap
K^{\mathrm{alg}}=\emptyset$. \hfill\CVD
\end{proposition}

\subsubsection{}\label{remremju}  A generic point
defines a bounded multiplicative semi-norm on $\H_K(X)$, and hence
defines a Berkovich point (cf. \cite{Ber}). The reader knowing the
language of Berkovich will not find difficulties in translating
the contents of this paper into the language of Berkovich.

For all $f(T)\in \H_K(\mathrm{D}^+(c,\rho))$, one has
\begin{equation}\label{Fund. Property of Gen Pts}
|f(t_{c,\rho})|_{_\Omega} \;=\; |f(T)|_{(c,\rho)}\; =\; \sup_{
\begin{smallmatrix}|x-c|\leq\rho\\
x\in
K^{\mathrm{alg}}
\end{smallmatrix}}|f(x)| \;=\;
\lim_{r\to\rho^{-}}\sup_{
\begin{smallmatrix}
|x-c|=r\\
x\in K^{\mathrm{alg}}
\end{smallmatrix}}|f(x)|\;.
\end{equation}
Hence, although the point $t_{c,\rho}$ is not uniquely determined
by the fact that $\mathrm{D}^-_\Omega(t_{c,\rho},\rho)\cap
K^{\mathrm{alg}}=\emptyset$, the norm $|\cdot|_{(c,\rho)}$ (i.e.
the Berkovich point $|\cdot|_{(c,\rho)}$) does not depend on the
choice of $t_{c,\rho}$.

By point iii) of Lemma \ref{log-convex property}, if $\rho\in |K|$
(resp. $\rho\in |K^{\mathrm{alg}}|$;\; $\rho\notin
|K^{\mathrm{alg}}|$), then one also has
$|f(t_{c,\rho})|=\max_{\begin{smallmatrix}|x|=\rho\\x\in
K\end{smallmatrix}}|f(x)|$ (resp.
$|f(t_{c,\rho})|=\max_{\begin{smallmatrix}|x|=\rho\\x\in
K^{\mathrm{alg}}\end{smallmatrix}}|f(x)|$;\; $|f(t_{c,\rho})|
=\lim_{r\to\rho^{-}} \max_{\begin{smallmatrix}
|x-c|=r\in|K^{\mathrm{alg}}|\\
x\in K^{\mathrm{alg}}\end{smallmatrix}}|f(x)|$).

\begin{proposition}[(\cite{Ber})]
Let
$X=\mathrm{D}^+(c_0,R_0)-\cup_{i=1,\ldots,n}\mathrm{D}^-(c_i,R_i)$
be an affinoid. Let $t_{c_i,R_i}\in X(\Omega)$ be the generic
point of $\mathrm{D}^{+}(c_i,R_i)$. Then, for all $f\in\H_K(X)$,
one has
\begin{equation}\label{Shilov Boundary -tekkkkkk}
\|f(T)\|_X = \max( |f(t_{c_{0},R_0})|_{\Omega} , \ldots ,
|f(t_{c_{n},R_n})|_{\Omega} )\;.
\end{equation}
\end{proposition}

\begin{lemma}\label{|f'| leq |f| / r_A}
Let
$X=\mathrm{D}^+(c_0,R_0)-\cup_{i=1,\ldots,n}\mathrm{D}^-(c_i,R_i)$
be an affinoid. Let $r_X:=\min(R_0,\ldots,R_n)$. Then
$\|\frac{d}{dT}f(T)\|_X\leq r_X^{-1}\|f(T)\|_X$.
\end{lemma}
\begin{proof} This follows easily from the Mittag-Leffler decomposition of
$f(T)$ together with the observations that
$\|f(T)\|_X=\max_{i=0,\ldots,n}(|f(t_{c_i,R_i})|)$ (cf.
\eqref{Shilov Boundary -tekkkkkk}), and $|f'(t_{c_i,R_i})|\leq
R_i^{-1}|f(t_{c_i,R_i})|$,  $\forall\;i$.
\end{proof}

\section{Discrete or analytic $\sigma-$modules and $(\sigma,\delta)-$modules}
\label{section - Discrete or analytic sigma and sigma,delta mod}

\begin{definition}
Let $\B$ be one of the rings of section \ref{rings}. We denote by
\begin{eqnarray}\index{Q(B)@$\Q(\B)$, $\Q_1(\B)$}
\Q(\B)&=&\{q\in K\;|\; \sigma_q:f(T)\mapsto f(qT)
\textrm{ is an automorphism of }\B\}\\
\Q_1(\B)&=&\Q(\B)\cap \mathrm{D}^-(1,1)\;.\label{Q and Q_1}
\end{eqnarray}
We will write $\Q$ and $\Q_1$ when no confusion is possible.
\end{definition}

Notice that $\Q(\B)\subset (K^{\times},|.|)$ is a topological
group and always contains  a disk $\mathrm{D}^-(1,\tau_0)$, for
some $\tau_0>0$. One has $\Q(\a_K(I))=\Q(\R_K)=\Q(\Hd_K)= \{q\in
K\;|\;|q|=1\}$. One sees easily that $\Q(\H_K(X))\subset\{q\in
K\;|\;|q|=1\}$ (cf. section \ref{k_0}, and Lemma
\ref{|q-1||c|<R-yt}).

\begin{definition}\index{S@$S^{\circ}=S-\bs{\mu}(\Q)$, $\ph{S}$}\index{m(Q)@$\bs{\mu}(\mathcal{Q})$}
Let $S\subseteq\Q$ be a subset. We denote by $\ph{S}$ the subgroup
of $\Q$ generated by $S$. Let $\bs{\mu}(\mathcal{Q})$ be the set
of all roots of unity belonging to $\Q$. Then we set
\begin{equation}\label{S'}
S^{\circ}:=S-\bs{\mu}(Q)\;.
\end{equation}
\end{definition}

\subsection{Discrete $\sigma-$modules}
By assumption, every finite dimensional free $\mathrm{B}-$module
$\M$ has the product topology.
\begin{definition}[(discrete $\sigma-$modules)]\label{definition of discrete sigma}
\index{sigma-Mod(B)_S^disc@$\sigma-\mathrm{Mod}(\B)_{S}^{\mathrm{disc}}$}
Let $S\subset\Q$ be an arbitrary subset. An object of
\begin{equation}
\sigma\!\!-\!\!\mathrm{Mod}(\B)_S^{\mathrm{disc}}
\end{equation}
is a finite dimensional free $\B-$module $\mathrm{M}$, together
with a group morphism
\begin{equation}\index{sigma^M@$\sigma^{\mathrm{M}}\protect{:}\ph{S}\to\mathrm{Aut}_K^{\mathrm{cont}}(\mathrm{M})$}
\sigma^{\mathrm{M}}:\ph{S}\xrightarrow[]{}\mathrm{Aut}_K^{\mathrm{cont}}(\mathrm{M})
\;,
\end{equation}
sending $q\mapsto \sigma_q^{\mathrm{M}}$, such that, for all $q\in
S$, the operator $\sigma_q^{\mathrm{M}}$ is
$\sigma_q-$semi-linear, that is
\begin{equation}\label{sigma_q-semilinear}
\sigma_q^{\mathrm{M}}(fm)=\sigma_q(f)\cdot\sigma_q^{\mathrm{M}}(m)\;,
\end{equation}
for all $f\in\B$, and all $m\in\mathrm{M}$. Objects
$(M,\sigma^{\mathrm{M}})$ in
$\sigma-\mathrm{Mod}(\B)_S^{\mathrm{disc}}$ will be called
\emph{discrete $\sigma-$modules over $S$}. A morphism between
$(\mathrm{M},\sigma^{\mathrm{M}})$ and
$(\mathrm{N},\sigma^{\mathrm{N}})$ is a $\mathrm{B}-$linear map
$\alpha:\mathrm{M}\to\mathrm{N}$ such that
\begin{equation}\label{morphism ==> commute simultaneously}
\alpha\circ\sigma_q^{\mathrm{M}}=\sigma_q^{\mathrm{N}}\circ\alpha\;,
\end{equation}
for all $q\in S$. We will denote the $K-$vector space of morphisms
by $\Hom^{\sigma}_S(\M,\N)$.\label{Hom^sigma_S(M,N)}
\index{Hom^sigma_S(M,N)@$\Hom^{\sigma}_S(\M,\N)$}
\end{definition}

\begin{notation}
If $S=\{q\}$ is reduced to a point, then the category of discrete
$\sigma-$modules over $\{q\}$ is the usual category of
$q-$difference modules. We will therefore use a simplified notation:
\begin{equation}\label{sigma_q-Mod = sigma-Mod^disc_q}\index{sigma_q-Mod(B)@$\sigma_q-\mathrm{Mod}(\B)$}
\sigma_q-\mathrm{Mod}(\B)\;\;:=\;\;\sigma-\mathrm{Mod}(\B)^{\mathrm{disc}}_{\{q\}}\;\;.
\end{equation}
\end{notation}
\begin{remark}\label{first remark}
1.--- Conditions \eqref{sigma_q-semilinear} and \eqref{morphism
==> commute simultaneously} for $q\in S$ imply the same conditions
for every $q\in \ph{S}$.

2.--- 
If $\M\neq 0$, the map
$\sigma^{\M}:\ph{S}\to\mathrm{Aut}_K^{\mathrm{cont}}(\mathrm{M})$
is injective. Indeed, since $\B$ is a domain and $\M$ is free, the
equality $\sigma_q^{\M}(fm)=\sigma_{q'}^{\M}(fm)$, $\forall \;
f\in\B$, $\forall\; m\in\M$, implies
$\sigma_q(f)\sigma_q^{\M}(m)=\sigma_{q'}(f)\sigma_{q'}^{\M}(m)$,
and hence the contradiction: $\sigma_q(f)=\sigma_{q'}(f),\forall
f\in \B$.

3.--- The morphism $\sigma^{\M}$ on $\ph{S}$ is determined by its
restriction to the set $S$.  Conversely, if a map $S\to\Autc(\M)$
is given, then this map extends to a group morphism
$\ph{S}\to\Autc(\M)$ if and only if the following conditions are
verified:
\begin{equation*}
\begin{array}{rl} i.&
\sigma_q^{\M}\circ\sigma_{q'}^{\M}=\sigma_{q'}^{\M}\circ\sigma_{q}^{\M}\;,
\;\;\textrm{ for all }q,q'\in S\;;\\
ii.& \textrm{If }\exists\; 
n,m\in\mathbb{Z}\;,\textrm{and}\;\;q_1,q_2\in S\;,\;\textrm{such that
}\;q_1^n=q_2^m\;,
\textrm{ then } (\sigma_{q_1}^{\M})^n=(\sigma_{q_2}^\M)^m\;;\\
iii.&\textrm{If }1\in S\;,\textrm{ then
}\sigma_1^{\M}=\mathrm{Id}\;.
\end{array}
\end{equation*}
\end{remark}

\subsubsection{Matrices of $\sigma^{\M}$.}\label{matrices A(q,T)}

Let $\e=\{\mathrm{e}_1,\ldots,\mathrm{e}_n\}\subset\M$ be a basis
over $\B$. If
$\sigma_q^{\M}(\mathrm{e}_i)=\sum_{j}a_{i,j}(q,T)\cdot\mathrm{e}_j$,
then in this basis $\sigma^{\M}_q$ acts as
\begin{equation}\label{sigma(f)=A(q,T)sigma(f)}
\sigma_q^{\M}(f_1,\ldots,f_n)=
(\sigma_q(f_1),\ldots,\sigma_q(f_n))\cdot A(q,T)\;,
\end{equation}
where $A(q,T):=(a_{i,j}(q,T))_{i,j}$. By definition
$A(1,T)=\mathrm{Id}$, and one has
\begin{equation}\label{A(q q',T)=A(q',qT) A(q,T)}
A(q q',T)=A(q',qT)\cdot A(q,T)\;.
\end{equation}
In particular $A(q^n,T)=A(q,q^{n-1}T)\cdot A(q,q^{n-2}T)\cdots
A(q,T)$.

\subsubsection{Internal $\mathrm{Hom}$ and $\otimes$.}

Let $(\M,\sigma^{\M})$, $(\mathrm{N},\sigma^{\mathrm{N}})$ be two
discrete $\sigma-$modules over $S$. We define a structure of
discrete $\sigma-$module on $\Hom_{\B}(\M,\N)$ by setting
\index{sigma_q^Hom@$\sigma_q^{\mathrm{Hom}(\M,\N)}$}
$\sigma_q^{\mathrm{Hom}(\M,\N)}(\alpha):=\sigma_q^{\N}\circ\alpha\circ(\sigma_q^{\M})^{-1}$,
for all $q\in S$, and all $\alpha\in\mathrm{Hom}_{\B}(\M,\N)$. We
define on $\M\otimes_{\B}\N$ a structure of discrete
$\sigma-$module over $S$ by setting
\index{sigma_q^tens@$\sigma_q^{\M\otimes\N}$}
$\sigma_q^{\M\otimes\N}(m\otimes
n):=\sigma_q^{\M}(m)\otimes\sigma_q^{\N}(n)$, for all $q\in S$,
and all $m\in\M$, $n\in\N$.

\subsubsection{}\label{Tannakian if S' neq empty} If $S^{\circ}\neq
\emptyset$ (cf. \eqref{S'}), then the category
$\sigma-\mathrm{Mod}(\B)^{\mathrm{disc}}_{S}$ is $K-$linear. If
$\B$ is a Bezout ring (i.e. every finitely generated ideal of $\B$
is principal), then $\sigma-\mathrm{Mod}(\B)^{\mathrm{disc}}_{S}$
is Tannakian (cf. \cite[12.3]{An-DV}). The ring $\H_K(X)$ is
always principal. If $K$ is spherically closed, then $\a_K(I)$,
$\R_K$, $\Hd_K$ are Bezout rings.

\subsubsection{}
As already mentioned in the introduction, the following is an
example of two \emph{non isomorphic} analytic $\sigma$-modules
over $X$, having isomorphic ``stalks'' at every $q\in U\subset
\Q(X)$. This is analogous to have \emph{non isomorphic} sheaves
having isomorphic stalks at every point.

\begin{example}\label{Morphisms between sigma modules
-examples} Let $X=\{|x|=1\}$, then $\Q(X)=\{x\in K\;|\;|x|=1\}$.
Let $U:=\mathrm{D}^-(1,1)$, and let $\pi\in K$ satisfy
$|\pi|=|p|^{\frac{1}{p-1}}$. Put then $A(q,x):=\exp(\pi (q-1) x)$,
and $\widetilde{A}(q,x):=\exp(\pi q(q-1) x)$. Let $\M$ (resp.
$\N$) be the discrete $\sigma$-module over $U$ defined by the
family $\{\sigma_q(Y)=A(q,x)\cdot Y\}_{q\in U}$ (resp.
$\{\sigma_q(Y)=\widetilde{A}(q,x)\cdot Y\}_{q\in U}$).  In this
fixed basis of $\M$ and $\N$, the matrices of every isomorphism
between $(\M,\sigma_q^{\M})$ and $(\N,\sigma_q^{\N})$ are of the
form $B(q,x)=\lambda \cdot\exp(\pi(1-q)x)\in\H_K(X)^{\times}$,
with $\lambda\in K^{\times}$. Hence for all $q\in U$ the equation
$\sigma_q(Y)=A(q,x)Y$ is isomorphic to
$\sigma_q(Y)=\widetilde{A}(q,x)Y$. But since $B(q,x)$ depends on
$q$, \emph{$\M$ and $\N$ are not isomorphic as analytic
$\sigma$-modules over $U$}.
\end{example}

\subsection{Discrete $(\sigma,\delta)-$modules}
Let $S\subset\Q(\B)$ be an arbitrary subset.
\begin{definition}[(discrete $(\sigma,\delta)-$modules)] An object of
\begin{equation}
\label{(sigma,delta)-Mod(B)_S^disc}
\index{sigma,delta-Mod(B)_S^disc@$(\sigma,\delta)-\mathrm{Mod}(\B)_{S}^{\mathrm{disc}}$}
(\sigma,\delta)\!\!-\!\!\mathrm{Mod}(\B)_S^{\mathrm{disc}}
\end{equation}
is a discrete $\sigma-$module over $S$, together with a
connection\footnote{i.e. $\delta_1^{\M}$ verifies
$\delta_1^{\M}(fm)=\delta_1(f)\cdot m+f\cdot\delta_1^{\M}(m)$,
$\forall\;f\in\B$, $\forall\;m\in\M$. Recall that
$\delta_1:=T\frac{d}{dT}$.}
 $\delta_1^{\M}:\M\to\M$.\label{delta_1^M:M-->M}\index{delta_1^M@$\delta_1^{\M}:\M\to\M$}
Objects $(\M,\sigma^{\mathrm{M}},\delta_1^{\M})$ of
$(\sigma,\delta)-\mathrm{Mod}(\B)_S^{\mathrm{disc}}$ will be
called \emph{discrete $(\sigma,\delta)-$modules over $S$}. A
morphism between $(\M,\sigma^{\M},\delta_1^{\M})$ and
$(\mathrm{N},\sigma^{\N},\delta_1^{\N})$ is a morphism
$\alpha:(\M,\sigma^{\M})\to(\N,\sigma^{\N})$ of discrete
$\sigma-$modules satisfying
\begin{equation}
\alpha\circ\delta_1^{\mathrm{M}}=\delta_1^{\mathrm{N}}\circ\alpha\;.
\end{equation}
We will denote the $K-$vector space of morphisms by
$\Hom^{(\sigma,\delta)}_S(\M,\N)$.\label{Hom^(sigma,delta)(M,N)}
\index{Hom^(sigma,delta)_S(M,N)@$\Hom^{(\sigma,\delta)}_S(\M,\N)$}
\end{definition}

\begin{notation}By analogy with \eqref{sigma_q-Mod = sigma-Mod^disc_q}, if
$S=\{q\}$, then we set:
\begin{equation}\label{sigma_q,delta_q-Mod = sigma,delta-Mod^disc_q}
(\sigma_q,\delta_q)-\mathrm{Mod}(\B)\;\;:=\;\;(\sigma,\delta)-\mathrm{Mod}(\B)^{\mathrm{disc}}_{\{q\}}\;\;.
\end{equation}
If $q=1$ we denote it by $\delta_1-\Mod(\B)$.
\end{notation}

As already mentioned in the introduction, we introduce the
operator
\begin{equation}
\delta_q^{\M}:=\sigma_q^{\M}\circ\delta_1^\M\;.
\end{equation}
For all $f\in\B$, all $m\in\M$, and all $q\in \ph{S}$, one has
\begin{equation}\label{delta(fm)=delta(f)sigma(m)+sigma(f)delta(m)}
\delta_q^{\M}(f\cdot
m)=\sigma_q(f)\cdot\delta_q^{\M}(m)+\delta_q(f)\cdot\sigma_q^{\M}(m)\;.
\end{equation}
Moreover, for all $\alpha\in\Hom^{(\sigma,\delta)}(\M,\N)$, and
all $q\in
\ph{S}$, one has 
$\alpha\circ\delta_q^{\mathrm{M}}=\delta_q^{\mathrm{N}}\circ\alpha$. 
Heuristically we imagine $\M$ as endowed with the map
$q\mapsto\delta_q^{\M}:\ph{S}\to\Endc(\M)$.  This justifies
notations \eqref{(sigma,delta)-Mod(B)_S^disc} and
\eqref{sigma_q,delta_q-Mod = sigma,delta-Mod^disc_q}.

\subsubsection{Matrices of $\delta_q^{\M}$.}

Let $\e=\{\mathrm{e}_1,\ldots,\mathrm{e}_n\}\subset\M$ be a basis
over $\B$. Let $A(q,T)\in GL_n(\B)$ be the matrix of
$\sigma_q^{\M}$ in the basis $\e$ (cf.
\eqref{sigma(f)=A(q,T)sigma(f)}). If
$\delta_q^{\M}(\mathrm{e}_i)=\sum_{j}g_{i,j}(q,T)\cdot\mathrm{e}_j$,
and if $G(q,T)=(g_{i,j}(q,T))_{i,j}$, then $\delta^{\M}_q$ acts in
the basis $\e$ as:
\begin{equation}\label{matrix of delta}
\delta_q^{\M}(f_1,\ldots,f_n)=(\delta_q(f_1),\ldots,\delta_q(f_n))\cdot
A(q,T)+(\sigma_q(f_1),\ldots,\sigma_q(f_n))\cdot G(q,T)\;.
\end{equation}
One has moreover
\begin{equation}\label{G(qq',T) = G(q',qT) cdot A(q,T)}
G(q'\cdot q,T) = G(q',qT)\cdot A(q,T)\;.
\end{equation}

\subsubsection{Internal $\mathrm{Hom}$ and $\otimes$.}

Let $(\M,\sigma^{\M},\delta^{\M})$,
$(\mathrm{N},\sigma^{\mathrm{N}},\delta^{\N})$ be two discrete
$(\sigma,\delta)-$modules over $S$.
We define a structure of discrete $(\sigma,\delta)-$module on
$\mathrm{Hom}_{\B}(\M,\mathrm{N})$ by setting
\begin{equation}\label{delta^Hom(M,N)}\index{delta_q^Hom@$\delta_q^{\mathrm{Hom}(\M,\mathrm{N})}$}
\delta_q^{\mathrm{Hom}(\M,\mathrm{N})}(\alpha):=
\Bigl(\delta_q^{\N}\circ\alpha-\sigma_q^{\mathrm{Hom}(\M,\N)}(\alpha)\circ\delta_q^{\M}\Bigr)\circ(\sigma_q^{\M})^{-1}\;.
\end{equation}
This definition gives the relation $\delta_q^{\N}(\alpha\circ m)=
 \sigma_q^{\mathrm{H}}(\alpha)\circ\delta_q^{\M}(m) +
 \delta_q^{\mathrm{H}}(\alpha)\circ\sigma_q^{\M}(m)$,
 for all
$\alpha\in\mathrm{Hom}_{\B}(\M,\N)$, and all $m\in\M$, where
$\mathrm{H}:=\mathrm{Hom}_{\B}(\M,\N)$.
We define on $\M\otimes_{\B}\N$ a structure of discrete
$(\sigma,\delta)-$module over $S$ by setting
\begin{equation}\label{delta_q^M otimes N}\index{delta_q^M tens N@$\delta_q^{\M\otimes\N}$}
\delta_q^{\M\otimes\N}(m\otimes
n):=\delta_q^{\M}(m)\otimes\sigma_q^{\N}(n)+\sigma_q^{\M}(m)\otimes\delta_q^{\N}(n)\;,
\end{equation}
for all $q\in S$, and all $m\in\M$, $n\in\N$.

\subsubsection{}
If $\B$ is Bezout, then
$(\sigma,\delta)-\Mod(\B)^{\mathrm{disc}}_S$ is $K-$linear and
Tannakian.

\subsection{Analytic $\sigma-$modules} Analytic
$\sigma-$modules are defined only if the ring $\B$ is equal to one
of the following rings: $\a_K(I)$, $\H_K(A)$, $\Hd_K(X)$, $\H_K$,
$\Hd_K$, $\R_K$. Notice that if $U\subset \Q(\B)$ is an open
subset, then the subgroup $\ph{U}\subseteq\Q(\B)$ generated by $U$
is open, i.e. $\ph{U}$ contains a disk $\mathrm{D}^-_K(1,\tau)$,
for some $\tau>0$.

\begin{definition}\label{definition of analytic sigma
module}\label{A(Q,T)} Let $\B:=\H_K(X)$. Let $(\M,\sigma^{\M})$ be
a discrete $\sigma-$module over $U$. Let $A(q,T)\in GL_n(\B)$ be
the matrix of $\sigma_q^{\M}$ in a fixed basis. We will say that
$(\M,\sigma^{\M})$ is an \emph{analytic $\sigma-$module} if, for
all $q\in U$, there exist a disk
$\D^-(q,\tau_q)=\{q'\;\;|\;\;|q'-q|<\tau_q\}$, with $\tau_q>0$,
and a matrix $A_q(Q,T)$ such that:
\begin{enumerate}
\item $A_q(Q,T)$ is an analytic element on the domain
$(Q,T)\in\mathrm{D}^-(q,\tau_q)\times X$\;; %
\item For all $q'\in\mathrm{D}^-_K(q,\tau_q)$, one has
$A_q(Q,T)_{|Q=q'}=A(q',T)$.
\end{enumerate}
This definition does not depend on the choice of basis $\e$. We
define
\begin{equation}\label{fgrt}\index{sigma-Mod(B)_U^an@$\sigma-\mathrm{Mod}(\B)_{U}^{\mathrm{an}}$}
\sigma-\mathrm{Mod}(\B)_U^{\mathrm{an}}
\end{equation} %
as the full sub-category of
$\sigma-\mathrm{Mod}(\B)_U^{\mathrm{disc}}$, whose objects are
analytic $\sigma-$modules. Let $I\subset\mathbb{R}_{\geq 0}$ be an
interval. We give the same definition over the ring $\B:=\a_K(I)$,
namely, if $\mathcal{C}(I):=\{|T|\in I\}$, the point i) is
replaced by
\begin{enumerate}
\item[i$'$)] $A_q(Q,T)$ is an analytic function on the domain
$(Q,T)\in\mathrm{D}^-(q,\tau_q)\times \mathcal{C}(I)$\;.
\end{enumerate}
\end{definition}

\begin{example}
The discrete $\sigma$-modules appearing in the Example
\ref{Morphisms between sigma modules -examples}, are actually
analytic.
\end{example}

\subsubsection{Analyticity of $\Hom(\M,\N)$ and $\M\otimes\N$.} If $(\M,\sigma^{\M})$ and
$(\N,\sigma^{\N})$ are two analytic $\sigma-$modules over $U$,
then $(\Hom(\M,\N),\sigma^{\Hom(\M,\N)})$ and
$(\M\otimes\N,\sigma^{\M\otimes\N})$ are analytic. This follows
from the explicit dependence of the matrices of
$\sigma^{\Hom(\M,\N)}$ and $\sigma^{\M\otimes\N}$ on the matrices
of $\sigma^{\M}$ and $\sigma^{\N}$.

\subsubsection{Discrete and analytic $\sigma-$modules over $\a_K(I)$, $\R_K$
and $\Hd_K(X)$.}

If $I_1\subset I_2$, then the restriction functor
$\sigma-\mathrm{Mod}(\a_K(I_2))_U^{\mathrm{an}}\to
\sigma-\mathrm{Mod}(\a_K(I_1))_U^{\mathrm{an}}$ is fully faithful.
Indeed the equality $f_{|_{I_1}}=g_{|_{I_1}}$ implies $f=g$, for
all $f,g\in\a_{K}(I_2)$ (analytic continuation
\cite[5.5.8]{Ch-Ro}).
\begin{definition}\label{def of an sigma mod over R_K and Hd_K}
Let $S\subseteq\Q$ be a subset, and let $U\subseteq\Q$ be an open
subset. We set
\begin{eqnarray}
\index{sigma-Mod(R_K)_U^an@$\sigma-\mathrm{Mod}(\R_K)_U^{\mathrm{an}}$,
$\sigma-\mathrm{Mod}(\R_K)_S^{\mathrm{disc}}$,
$\sigma-\mathrm{Mod}(\Hd_K)_U^{\mathrm{an}}$,
$\sigma-\mathrm{Mod}(\Hd_K)_S^{\mathrm{disc}}$}
\sigma-\mathrm{Mod}(\R_K)_U^{\mathrm{an}}&:=&\bigcup_{\varepsilon
>0}\sigma-\mathrm{Mod}(\a_K(]1-\varepsilon,1[))_U^{\mathrm{an}}\;;\\
\sigma-\mathrm{Mod}(\R_K)_S^{\mathrm{disc}}&:=&\bigcup_{\varepsilon
>0}\sigma-\mathrm{Mod}(\a_K(]1-\varepsilon,1[))_S^{\mathrm{disc}}\;.
\end{eqnarray}
Similarly, one can define
$\sigma-\mathrm{Mod}(\Hd_K(X))_U^{\mathrm{an}}$ and
$\sigma-\mathrm{Mod}(\Hd_K(X))_S^{\mathrm{disc}}$.
\end{definition}

\begin{remark}
Since $U$ is open, one has $U^{\circ}\neq \emptyset$ (cf.
\eqref{S'}). By Section \ref{Tannakian if S' neq empty}, if $\B$
is one of the previous rings (and if it is a Bezout ring), then
$\sigma-\Mod(\B)^{\mathrm{an}}_U$ is $K-$linear and Tannakian.
\end{remark}

\subsection{Analytic $(\sigma,\delta)-$modules}\label{analytic
sigma delta modules} We maintain the previous notations. In
section \ref{construction of delta} below we define a \emph{fully
faithful} functor
\begin{equation}
(\textrm{Forget
}\delta)^{-1}:\;\;\sigma-\mathrm{Mod}(\B)_U^{\mathrm{an}}\xrightarrow[]{\;\;\quad\;\;}
(\sigma,\delta)-\mathrm{Mod}(\B)_U^{\mathrm{disc}}\;,
\end{equation}
which is a ``local'' section of the functor $\textrm{Forget
}\delta:(\sigma,\delta)-\mathrm{Mod}(\B)_U^{\mathrm{disc}}
\xrightarrow[]{} \sigma-\mathrm{Mod}(\B)_U^{\mathrm{disc}}$.
\label{definition of analytic sigma delta modules;} The essential
image of the functor $(\textrm{Forget }\delta)^{-1}$ will be
denoted by
\begin{equation}\index{sigma,delta-Mod(B)_U^an@$(\sigma,\delta)-\mathrm{Mod}(\B)_{U}^{\mathrm{an}}$}
(\sigma,\delta)-\mathrm{Mod}(\B)_U^{\mathrm{an}}\;.
\end{equation}
By definition, the functor which ``forgets'' the action of
$\delta$ is therefore an equivalence
\begin{equation}\label{sigma-an=(sigma,delta)-an}
(\sigma,\delta)-\mathrm{Mod}(\B)_U^{\mathrm{an}}\xrightarrow[\sim]{\textrm{Forget
}\delta}\sigma-\mathrm{Mod}(\B)_U^{\mathrm{an}}\;.
\end{equation}

Notice that a morphism between analytic $(\sigma,\delta)-$modules
is, by definition, a morphism of \emph{discrete}
$(\sigma,\delta)-$modules.

\subsubsection{Construction of $\delta$.} \label{construction of
delta} \label{sigma to sigma,delta is fully faithful} Let
$(\M,\sigma^{\M})$ be an analytic $\sigma-$module. We shall define
a $(\sigma,\delta)$-module structure on $\M$. It follows from
definitions \ref{A(Q,T)} and \ref{def of an sigma mod over R_K and
Hd_K} that the map
$q\mapsto\sigma^{\M}_q:\ph{U}\to\mathrm{Aut}_K(\M)$ is
\emph{derivable}, in the sense that, for all $q\in \ph{U}$, the
limit
\begin{equation}\label{delta:=q lim q'-->q sigma_q'-sigma_q/q'-q}
\delta_q^{\M}:=q\cdot\lim_{q'\to
q}\frac{\sigma_{q'}^{\mathrm{M}}-\sigma_q^{\mathrm{M}}}{q'-q}=\textrm{``
}(q\frac{d}{dq}\sigma^M)(q)\textrm{ ''}
\end{equation}
exists in $\mathrm{End}_K^{\mathrm{cont}}(\M)$, with respect to
the simple convergence topology (cf. \eqref{G(q,T)=qd/dqA(q,T)}).
Moreover, for all $q\in \ph{U}$, the rule
\eqref{delta(fm)=delta(f)sigma(m)+sigma(f)delta(m)} holds, and
$\delta^{\M}_q=\sigma_q^{\M}\circ\delta_1^{\M}$.

Let $\alpha:(\M,\sigma^{\M})\to(\mathrm{N},\sigma^{\mathrm{N}})$
be a morphism of analytic $\sigma-$modules, that is
$\alpha\circ\sigma_q^{\M} = \sigma_q^{\N}\circ\alpha$, for all
$q\in U$. Passing to the limit in the definition \eqref{delta:=q
lim q'-->q sigma_q'-sigma_q/q'-q}, one shows that $\alpha$
commutes with $\delta^{\M}_q$, for all $q\in U$. Hence the
inclusion $\mathrm{Hom}^{(\sigma,\delta)}_{U}(\M,\N)\subseteq
\mathrm{Hom}^{\sigma}_U(\M,\mathrm{N})$ is an equality. If
$\e=\{\mathrm{e}_1,\ldots,\mathrm{e}_n\}\subset\M$ is a basis in
which the matrix of $\sigma_q^{\M}$ is $A(q,T)$, 
then the matrix of $\delta_q^{\M}$ is (cf. \eqref{matrix of
delta}, Def. \ref{A(Q,T)} and \ref{def of an sigma mod over R_K
and Hd_K})
\begin{equation}\label{G(q,T)=qd/dqA(q,T)}
G(q,T)\;\;:=\;\;q\cdot\lim_{q'\to
q}\frac{A(q',T)-A(q,T)}{q'-q}\;\;=\;\;\Bigl(\partial_Q\bigl(
A_q(Q,T)\bigr)\Bigr)_{|_{Q=q}}\;,
\end{equation}
where $\partial_Q$ is the derivation $Q\frac{d}{dQ}$, and
$A_q(Q,T)$ is the matrix of Definition \ref{A(Q,T)}.

\begin{remark}\label{Definition of Conf_U}
By the above definitions, there is an obvious functor
\begin{equation}\index{Conf_U@$\Conf_U$}
\Conf_U:\sigma-\Mod(\B)_U^{\mathrm{an}}\xrightarrow[]{\qquad}\delta_1-\Mod(\B)\;,
\end{equation}
obtained by composing $(\textrm{Forget }\delta)^{-1}$ (cf.
\eqref{sigma-an=(sigma,delta)-an}) with $\textrm{Forget
}\sigma:(\sigma,\delta)-\Mod(\B)^{\mathrm{an}}_U\xrightarrow[]{\quad}\delta_1-\Mod(\B)$.
\end{remark}

\section{Solutions (formal definition)}
\label{section - solution formal def}

\subsection{Discrete $\sigma-$algebras and $(\sigma,\delta)-$algebras}

Let $S\subseteq\Q(\B)$ be a 
subset.
\begin{definition}[(Discrete $\sigma-$algebra over $S$)]\label{discrete sigma algebra}
\index{C^sigma_S@$\mathrm{C}^{\sigma}_S$} A
$\mathrm{B}$-\emph{discrete $\sigma-$algebra over $S$}, or simply
a \emph{discrete $\sigma-$algebra over $S$}  is a $\B-$algebra
$\mathrm{C}$ such that:
\begin{enumerate}
\item $\C$ is an \emph{integral domain}, %
\item there exists a group morphism
$\sigma^{\C}:\ph{S}\to\mathrm{Aut}_K(\C)$ such that
$\sigma^{\C}_q$ is a ring automorphism extending $\sigma_q^{\B}$,
for
all $q\in \ph{S}$; %
\item one has $\mathrm{C}^\sigma_S=K$, where
$\mathrm{C}^{\sigma}_S:= \{c\in \mathrm{C}\;|\;
\sigma_q(c)=c\;,\textrm{ for all }q\in S\}$.
\end{enumerate}
We will call $\mathrm{C}^\sigma_S$ \emph{the sub-ring of
$\sigma-$constants of }$\mathrm{C}$.
 We will write
$\sigma_q$ instead of $\sigma_q^{\C}$, when no confusion is
possible.
\end{definition}

Observe that no topology is required on $\mathrm{C}$. The word
\emph{discrete} is employed, here and later on, to emphasize that
we do not ask ``continuity'' with respect to $q$.
Notice also that if a discrete $\sigma-$algebra $\mathrm{C}$ is
free and of finite rank as $\B-$module, then it is a discrete
$\sigma-$module.

\subsubsection{}\label{no discrete sigma-alg over root of unity}
If $S^{\circ}\neq \emptyset$ (cf. \eqref{S'}), then
$\B^{\sigma}_S=K$, and $\B$ itself is a discrete $\sigma-$algebra
over $S$.  On the other hand, If $S=\{\xi\}$ is reduced to a root
of unity $\xi\in\bs{\mu}(\Q)$, since
$\B^{\sigma}_{S}=\B^{\sigma_\xi}\neq K$, it follows that $\B$
itself is not a discrete $\sigma-$algebra over $S$. Hence there is
no discrete $\sigma-$algebra over $S=\{\xi\}$. To deal with this
problem we introduce the following

\begin{definition}[(Discrete $(\sigma,\delta)-$algebra over
$S$)]\label{discrete sigma,delta algebra}
\index{C^sigma,delta_S@$\C^{(\sigma,\delta)}_{S}$} A
\emph{discrete $(\sigma,\delta)-$algebra $\mathrm{C}$ over $S$} is
a $\B-$algebra such that:
\begin{enumerate}
\item $\C$ satisfies properties i) and ii) of Definition
\ref{discrete sigma algebra}, %
\item there exists a derivation $\delta_1^{\C}$, extending $\delta_1=T\frac{d}{dT}$ on $\B$, and commuting with $\sigma_q^{\mathrm{C}}$, for all $q\in\ph{S}$, %
\item one has $\C^{(\sigma,\delta)}_{S}=K$, where
$\C^{(\sigma,\delta)}_{S}:=\{f\in\C
\;|\;f\in\mathrm{C}^{\sigma}_S\;,\textrm{ and }\;\delta_1(f)=0\}$.
\end{enumerate}
We will call $\C^{(\sigma,\delta)}_S$ the \emph{sub-ring of
$(\sigma,\delta)-$constants of $\mathrm{C}$}. We will write
$\delta_1$ instead of $\delta_1^{\C}$, if no confusion is
possible.
\end{definition}

The operator $\delta_q^{\C}:=\sigma_q^{\C}\circ\delta_1^{\C}$
satisfies property
\eqref{delta(fm)=delta(f)sigma(m)+sigma(f)delta(m)}.  Since
$\B^{(\sigma,\delta)}_S=K$, it follows that $\B$ is always a
$(\sigma,\delta)-$algebra over $S$, for an arbitrary sub-set
$S\subseteq\Q(\B)$, even for $S=\{\xi\}$, with
$\xi\in\bs{\mu}(\Q(\B))$.

\subsection{Constant Solutions }\label{section constant solutions}

\begin{definition}[(Constant solutions on $S$)]
\label{discrete solution} Let $(\mathrm{M},\sigma^{\mathrm{M}})$
(resp.$(\M,\sigma^{\M},\delta^{\M})$) be a \emph{discrete}
$\sigma-$module (resp. $(\sigma,\delta)-$module) over $S$, and let
$\mathrm{C}$ be a discrete $\sigma-$algebra (resp.
$(\sigma,\delta)-$algebra) over $S$. A \emph{constant solution} of
$\M$, with values in $\mathrm{C}$, is a $B-$linear morphism
$$\alpha:\mathrm{M}\longrightarrow \mathrm{C}$$ such that
$\alpha\circ\sigma_q^{\M}=\sigma_q^{\C}\circ\alpha$, for all $q\in
S$ (resp. $\alpha$  simultaneously satisfies
$\alpha\circ\delta_1^{\M}=\delta_1^{\C}\circ\alpha$, and
$\alpha\circ\sigma_q^{\M}=\sigma_q^{\C}\circ\alpha$, for all $q\in
S$). We denote by $\mathrm{Hom}_S^{\sigma}(\M,\mathrm{C})$ (resp.
$\mathrm{Hom}_S^{(\sigma,\delta)}(\M,\mathrm{C})$) the $K-$vector
space of the solutions of $\M$ in $\mathrm{C}$.
\end{definition}

\subsubsection{Matrices of solutions.}
\label{morphisms as solutions} Let $\M$ be a discrete
$\sigma-$module (resp. $(\sigma,\delta)-$module). Let $\C$ be a
discrete $\sigma-$algebra (resp. $(\sigma,\delta)-$algebra) over
$S$. Recall that, if $S=\{\xi\}$, with $\xi^n=1$, then there is no
discrete $\sigma-$algebra, over $S$ (cf. Section \ref{no discrete
sigma-alg over root of unity}).

Let $\e=\{\mathrm{e}_1,\ldots,\mathrm{e}_n\}$ be a basis of $\M$,
and let $A(q,T)$ (resp. $G(q,T)$) be the matrix of $\sigma_q^{\M}$
(resp. $\delta_q^{\M}$) in this basis (cf. \eqref{matrix of
delta}). We identify a morphism $\alpha:\M\to\C$ with the vector
$(y_i)_{i}\in\C^n$, given by $y_i:=\alpha(\mathrm{e}_i)$. In this
way constant solutions become solutions in the usual \emph{vector
form}. Indeed
\begin{eqnarray}
\left(\begin{smallmatrix}
\sigma_q(y_1)\\
\begin{picture}(1,12)
\put(-3,1){$\vdots$}
\end{picture}\\
\sigma_q(y_n)
\end{smallmatrix}\right)&=&A(q,T)\cdot
\left(
\begin{smallmatrix}
y_1\\
\begin{picture}(1,12)
\put(-3,1){$\vdots$}
\end{picture}\\
y_n
\end{smallmatrix}
\right)
\;,\;\;\textrm{ for all }q\in S\;,\\
\Bigl(\textrm{resp.
}\left(
\begin{smallmatrix}
\delta_q(y_1)\\
\begin{picture}(1,12)
\put(-3,1){$\vdots$}
\end{picture}\\
\delta_q(y_n)
\end{smallmatrix}
\right)&=&G(q,T)\cdot \left(
\begin{smallmatrix}
y_1\\
\begin{picture}(1,12)
\put(-3,1){$\vdots$}
\end{picture}\\
y_n
\end{smallmatrix}
\right)\;,\;\;\textrm{ for all }q\in S\;\;\Bigr)\;.
\end{eqnarray}

\begin{definition} By a \emph{fundamental matrix of solutions}
 of $\M$ (in the basis
$\e$) we mean a matrix $Y\in GL_n(\mathrm{C})$ satisfying
\emph{simultaneously}
\begin{equation}\label{sigma(Y)=YA}
\sigma_q(Y)=A(q,T)\cdot Y\;,\quad\textrm{ for all }q\in S\;,
\end{equation}
(resp. satisfying \emph{simultaneously}
\begin{equation}\label{delta(Y)=YG}
\left\{\begin{array}{rcl}
\sigma_q(Y)&=& A(q,T)\cdot Y\;,\quad\textrm{ for all }q\in S\;,\\
\delta_1(Y)&=& G(1,T)\cdot Y\;.\quad\phantom{\textrm{ for all
}q\in S}\quad)\;.
\end{array}\right.
\end{equation}
\end{definition}

\subsubsection{Unit object and $\sigma$-constants.}
Let $\mathbb{I}=\B$ be the unit object. By the description given
above, every solution $\alpha\in\Hom^{\sigma}_S(\mathbb{I},\C)$
(resp. $\alpha\in\Hom^{(\sigma,\delta)}_S(\mathbb{I},\C)$) can be
identified with $y:=\alpha(1)\in \C^{\sigma}_S$ (resp.
$y:=\alpha(1)\in\C^{(\sigma,\delta)}_S$). We obtain
$\C^\sigma_S\cong\Hom^{\sigma}_S(\mathbb{I},\C)$ (resp.
$\C^{(\sigma,\delta)}_S\cong\Hom^{(\sigma,\delta)}_S(\mathbb{I},\C)$).
In particular $\B^{\sigma}_S$ (resp. $\B^{(\sigma,\delta)}_S$) is
identified with $\mathrm{End}^{\sigma}_S(\mathbb{I})$ (resp.
$\mathrm{End}_S^{(\sigma,\delta)}(\mathbb{I})$), and the category
is $K$-linear if and only if $\B^{\sigma}_S=K$ (resp.
$\B^{(\sigma,\delta)}_S=K$).

\subsubsection{Dimension of the space of solutions.}

Let $F:=\mathrm{Frac}(\mathrm{C})$ be the fraction field of
$\mathrm{C}$, then both $\sigma_q$ and $\delta_1$ extend to $F$
(cf. \cite[Ex.1.5]{VS}).
\begin{lemma}[(Wronskian Lemma)]
Let $\M$ be a $(\sigma,\delta)-$module (resp. $\sigma-$module)
over $S$, and let $\C$ be a discrete $(\sigma,\delta)-$algebra
(resp. $\sigma-$algebra) over $S$. One has
\begin{equation}\label{dim_K Hom(M,C) <= rk M}
\mathrm{dim}_K \mathrm{Hom}_S^{(\sigma,\delta)}(\M,\mathrm{C})\leq
\mathrm{rk}_{\B}(\M)\;.
\end{equation}
(resp. if $S^{\circ}\neq \emptyset$ (cf. \eqref{S'}), then
$\mathrm{dim}_K \mathrm{Hom}_S^{\sigma}(\M,\mathrm{C})\leq
\mathrm{rk}_{\B}(\M)$.)
\end{lemma}
\begin{proof} One has $\mathrm{dim}_K
\mathrm{Hom}_S^{(\sigma,\delta)}(\M,\mathrm{C})\leq\mathrm{dim}_K
\mathrm{Hom}^{\delta_1}(\M,\mathrm{C}) \leq
\mathrm{rk}_{\B}(\M)\;.$
On the other hand, if $q\in S^{\circ}$, then 
$\mathrm{Hom}^{\sigma_q}(\M,\mathrm{C}) \leq \mathrm{rk}_{\B}(\M)$
(cf. \cite[Lemma 1.1.11]{DV-Inventiones}). Hence $\mathrm{dim}_K
\mathrm{Hom}_S^{\sigma}(\M,\mathrm{C})\leq\mathrm{dim}_K
\mathrm{Hom}^{\sigma_q}(\M,\mathrm{C}) \leq
\mathrm{rk}_{\B}(\M)$.\end{proof}

\section{$\C$-Constant Confluence}
In this section we state the formal results regarding confluence.
We introduce the notion of $\C$-constant modules. As explained in
the introduction, this notion is an adaptation of the notion of
$\C$-admissibility in the sense of representation theory. On the
other hand it can be interpreted as a generalization of the Galois
theory for differential and $q$-difference equations. According to
this point of view, in our context we have the problem that the
analogue of the Picard-Vessiot algebra trivializing a given object
$\M$ does not exist for arbitrary objects $\M$. Also the
uniqueness of the Picard-Vessiot algebra remains an open problem.
We avoid these problems by working with the category of modules
trivialized by a given algebra $\C$ which is fixed once and for
all. We hope that this problem will be overcome in the future.

\label{theory of deformation}

\subsection{$\C$-Constant modules}
Let $\B$ be one of the rings of Sections \ref{rings} and
\ref{affinoid - section -tfrg jj}, let $S\subset\Q(\B)$ be a
subset, and let $U\subset\Q(\B)$ be an open subset.
\begin{definition}[($\C$-Constant modules)]\label{definition of constant}
Let $M$ be a discrete $\sigma-$module over $S$. We will say that
$M$ is $\C$-\emph{constant} on $S$, or equivalently that
\emph{$\M$ is trivialized by $\C$}, if there exists a discrete
$\sigma-$algebra $\mathrm{C}$ over $S$ such that
\begin{eqnarray}
\mathrm{dim}_K\mathrm{Hom}_S^{\sigma\phantom{,\delta)}}(\M,\mathrm{C})&=&\mathrm{rk}_{\B}\M\;.
\end{eqnarray}
We give the analogous definition for $(\sigma,\delta)$-modules.
The full sub-category of
$\sigma-\mathrm{Mod}(\B)^{\mathrm{disc}}_S$ (resp.
$(\sigma,\delta)-\mathrm{Mod}(\B)^{\mathrm{disc}}_S$), whose
objects are trivialized by $\mathrm{C}$, will be denoted by
\begin{equation}\label{hihi}
\index{sigma-Mod(B,C)^const_S@$\sigma-\mathrm{Mod}(\B,\mathrm{C})^{\mathrm{const}}_S$}
\index{sigma,delta-Mod(B,C)^const_S@$(\sigma,\delta)-\mathrm{Mod}(\B,\mathrm{C})^{\mathrm{const}}_S$}
\sigma-\mathrm{Mod}(\B,\mathrm{C})^{\mathrm{const}}_S\qquad
(\textrm{resp. }\;
(\sigma,\delta)-\mathrm{Mod}(\B,\mathrm{C})^{\mathrm{const}}_S\;\;)\;.
\end{equation}
The full subcategory of
$\sigma-\mathrm{Mod}(\B,\mathrm{C})^{\mathrm{const}}_U$ (resp.
$(\sigma,\delta)-\mathrm{Mod}(\B,\mathrm{C})^{\mathrm{const}}_U$)
whose objects are analytic will be denoted by
\begin{equation}
\index{sigma,delta-Mod(B,C)an,const_U@$(\sigma,\delta)-\mathrm{Mod}(\B,\mathrm{C})^{\mathrm{an,const}}_U$}
\sigma-\mathrm{Mod}(\B,\mathrm{C})^{\mathrm{an,const}}_U\qquad(\textrm{resp.
}\;(\sigma,\delta)-\mathrm{Mod}(\B,\mathrm{C})^{\mathrm{an,const}}_U\;\;)\;.
\end{equation}
\end{definition}

 Notice that $\M$ is trivialized by
$\mathrm{C}$ if there exists $Y\in GL_n(\mathrm{C})$,
$n:=\mathrm{rk}_{\B}\M$, such that $Y$ is \emph{simultaneously} a
solution, for all $q\in S$, of the family of equations
\eqref{sigma(Y)=YA} (resp. both the conditions of
\eqref{delta(Y)=YG}). Roughly speaking, $\M$ is $\C$-constant on
$S$ if it admits a basis of $q-$solutions in $GL_n(\C)$ which
``does not depend on $q\in S$''.

\begin{lemma}\label{tensor product is again strongly conf:key lemma}
Let $\M$, $\N$ be two discrete $\sigma-$modules (resp.
$(\sigma,\delta)-$modules). If $\M$, $\N$ are both trivialized by
$\mathrm{C}$, then $\M\otimes \N$, $\Hom(\M,\N)$, $M^{\vee}$,
$N^{\vee}$ are trivialized by $\mathrm{C}$.
\end{lemma}
\begin{proof} The fundamental matrix solution
of $\M\otimes\N$ (resp. $\mathrm{Hom}(\M,\N)$) is obtained 
by taking products of entries of the two matrices of solutions of
$\M$ and $\N$ respectively. Hence ``it does not depend on $q\in
S$''. The assertion on $\M^{\vee}$, $\N^{\vee}$ is a particular
case of the previous one.\end{proof}

\begin{lemma}\label{The forget from S to S' is fully faithful}
\index{Res^S_S'@$\mathrm{Res}^{S}_{S'}$} Let $S'\subseteq S$ be a
non empty subset. Let $\mathrm{C}$ be a discrete
$(\sigma,\delta)-$algebra over $S$. Then the restriction functor
$\mathrm{Res}^{S}_{S'}$, sending $(\M,\sigma^{\M},\delta^{\M}_1)$
into $(\M,\sigma^{\M}_{|_{\ph{S'}}},\delta^{\M}_{1})$:
\begin{equation}
\mathrm{Res}^{S}_{S'}:(\sigma,\delta)-\mathrm{Mod}(\B,\mathrm{C})^{\mathrm{const}}_{S}\longrightarrow
(\sigma,\delta)-\mathrm{Mod}(\B)^{\mathrm{disc}}_{S'}\;\;
\end{equation}
is fully faithful and its image is contained in the category 
$(\sigma,\delta)-\mathrm{Mod}(\B,\mathrm{C})^{\mathrm{const}}_{S'}$.
The same fact is true for discrete $\sigma-$modules under the
assumption: $(S')^{\circ}\neq \emptyset$.
\end{lemma}
\emph{Proof. } The proof is the same in both cases: here we give
the proof in the case of $(\sigma,\delta)-$modules. We must show
that the inclusion $\mathrm{Hom}_S^{(\sigma,\delta)}(\M,\N)\to
\mathrm{Hom}_{S'}^{(\sigma,\delta)}(\M,\N)$
 is an isomorphism, for all $\M,\N$ in
$(\sigma,\delta)-\mathrm{Mod}(\B,\mathrm{C})^{\mathrm{const}}_{S}$.
In other words, we have to show that if $\alpha:\M\to\N$ commutes
with $\sigma_{q'}$, for all $q'\in S'$, then it commutes also with
$\sigma_q$, for all $q\in S$. One has
\begin{eqnarray}
\mathrm{Hom}_S^{(\sigma,\delta)} (\M,\N)&=&\Hom^{(\sigma,\delta)}_S(\M\otimes\N^{\vee},\B)\;;\\
\mathrm{Hom}_{S'}^{(\sigma,\delta)}(\M,\N)&=&\Hom^{(\sigma,\delta)}_{S'}(\M\otimes\N^{\vee},\B)\;.
\nonumber
\end{eqnarray}
Observe that $\M\otimes\N^{\vee}$ is the dual of the ``internal
hom'' $\mathrm{Hom}(\M,\N)$. By Lemma \ref{tensor product is again
strongly conf:key lemma}, $\M\otimes\N^{\vee}$ is trivialized by
$\mathrm{C}$. The restriction of $\M\otimes\N^{\vee}$ to $S'$ is
obviously $\C$-constant on $S'$, since it is trivialized by
$\mathrm{C}$. This implies that
\begin{equation}
\Hom^{(\sigma,\delta)}_S(\M\otimes\N^{\vee},\mathrm{C})=
\Hom^{(\sigma,\delta)}_{S'}(\M\otimes\N^{\vee},\mathrm{C})\;.
\end{equation}
This shows that a morphism with values in $\B\subseteq\mathrm{C}$
commutes with all $\sigma_q$ and $\delta_q$, for all $q\in S$, if
and only if it commutes with all $\sigma_q$ and $\delta_q$, for
all $q\in S'$. Hence
\begin{equation}
\Hom^{(\sigma,\delta)}_S(\M\otimes\N^{\vee},\B)=\Hom^{(\sigma,\delta)}_{S'}(\M\otimes\N^{\vee},\B)\;.\qquad\Box
\end{equation}

\subsubsection{Restriction to a roots of unity.} By the previous lemma, if $\xi\in
S\cap\bs{\mu}(\Q)$, then
\begin{equation}
\mathrm{Res}^{\;S}_{\{\xi\}}:(\sigma,\delta)-\mathrm{Mod}(\B,\mathrm{C})^{\mathrm{const}}_{S}\longrightarrow
(\sigma_\xi,\delta_\xi)-\mathrm{Mod}(\B)
\end{equation}
is again fully faithful. On the other hand, if
$S^\circ\neq\emptyset$, then the restriction
\begin{equation}
\mathrm{Res}^{\;S}_{\{\xi\}}:\sigma-\mathrm{Mod}(\B,\mathrm{C})^{\mathrm{const}}_{S}\longrightarrow
\sigma_\xi-\mathrm{Mod}(\B)
\end{equation}
is \emph{not} fully faithful, since
$\sigma-\mathrm{Mod}(\B,\mathrm{C})^{\mathrm{const}}_{S}$ is
$K-$linear, while $\sigma_\xi-\mathrm{Mod}(\B)$ is not $K-$linear
(i.e. $K\subset\mathrm{End}(\mathbb{I})$, but $K\neq
\mathrm{End}(\mathbb{I})$, cf. Section \ref{section - notations}).

\subsubsection{The case of an open subset.}\label{fresfresfres}
We observe that if $U$ is open, then the condition $U^{\circ}\neq
\emptyset$ is automatically verified. Hence, by Lemma \ref{The
forget from S to S' is fully faithful}, if $S\subset U$ is a (non
empty) subset, the restriction
\begin{equation}\mathrm{Res}^{U}_{S}:
(\sigma,\delta)-\mathrm{Mod}(\B,\mathrm{C})^{\mathrm{an,const}}_U\xrightarrow[]{\quad\quad}
(\sigma,\delta)-\mathrm{Mod}(\B,\mathrm{C})^{\mathrm{const}}_{S}
\end{equation}
is fully faithful. The same is true for $\sigma-$modules, under
the assumption $S^\circ\neq \emptyset$. In particular, if
$U'\subset U$ is an open subset, then the restriction functor is
fully faithful:
\begin{equation}\mathrm{Res}^{U}_{U'}:
(\sigma,\delta)-\mathrm{Mod}(\B,\mathrm{C})^{\mathrm{an,const}}_U\xrightarrow[]{\quad\quad}
(\sigma,\delta)-\mathrm{Mod}(\B,\mathrm{C})^{\mathrm{an,const}}_{U'}\;.
\end{equation}

\subsection{$\C$-Constant deformation and $\C$-constant confluence}\label{constant deformation section}
In this section we give the formal definition of the confluence
and deformation functors. As usual $S\subseteq\Q(\B)$ is an
arbitrary subset, and $U\subseteq\Q(\B)$ is an open subset.
\begin{definition}[(Extensible objects)]
\label{definition of strongly confluent at q} Let $q\in S$. Let
$\mathrm{C}$ be a discrete $\sigma-$algebra over $S$. A
$q-$difference module $\M$ is said to be \emph{$\C$-extensible to
$S$} if it belongs to the essential image of the restriction
functor
\begin{eqnarray*}\label{forget functor from S to q}
\mathrm{Res}^{\;S}_{\{q\}}:
\sigma-\mathrm{Mod}(\B,\mathrm{C})^{\mathrm{const}}_{S}&\longrightarrow&
\sigma_q-\mathrm{Mod}(\B)\;.
\end{eqnarray*}
\index{sigma_q-Mod(B,C)_S@$\sigma_q-\mathrm{Mod}(\B,\mathrm{C})_{S}$}
\index{sigma_q,delta_q-Mod(B,C)_S@$(\sigma_q,\delta_q)-\mathrm{Mod}(\B,\mathrm{C})_{S}$}
The \emph{full sub-category} of $\sigma_q-\mathrm{Mod}(\B)$ whose
objects are $\C$-extensible to $S$, will be denoted by
$\sigma_q-\mathrm{Mod}(\B,\mathrm{C})_{S}$. If $U$ is open, and if
$q\in U$, we will denote by
\begin{equation}\label{siiisiiisi}
\index{sigma_q,delta_q-Mod(B,C)^an_U@$(\sigma_q,\delta_q)-\mathrm{Mod}(\B,\mathrm{C})^{\mathrm{an}}_{U}$}
\index{sigma_q-Mod(B,C)^an_U@$\sigma_q-\Mod(B,\C)^{\mathrm{an}}_U$}
\sigma_q-\Mod(B,\C)^{\mathrm{an}}_U
\end{equation} the
full sub-category of $\sigma_q-\Mod(\B)_U$ whose objects belong to
the essential image of
$\sigma-\mathrm{Mod}(\B,\C)^{\mathrm{an,const}}_{U}$. We give
analogous definitions for $(\sigma,\delta)$-modules.
\end{definition}

Lemma \ref{The forget from S to S' is fully faithful} and
Definition \ref{definition of strongly confluent at q} easily  give
the following formal statement:
\begin{corollary}\label{forget functor from S to q strongss}
With the notation of Lemma \ref{The forget from S to S' is fully
faithful}, one has an equivalence
\begin{equation}\label{forget functor from S to q strong}
\mathrm{Res}^{\;S}_{\{q\}}\;:\;(\sigma,\delta)-\mathrm{Mod}(\B,\mathrm{C})^{\mathrm{const}}_{S}
\xrightarrow[]{\;\;\sim\;\;}
(\sigma_q,\delta_q)-\mathrm{Mod}(\B,\mathrm{C})_{S}\;.
\end{equation}
The same fact is true for $\sigma-$modules, under the additional
hypothesis that $q\in S^{\circ}$. \hfill\CVD
\end{corollary}

\begin{definition}\label{Conf_q,q'}
1.-- Let $S\subseteq\Q(\B)$ be a subset and let $q,q'\in \ph{S}$.
We will call the \emph{$\C$-constant deformation functor}, denoted
by
\begin{equation}\label{Def_q,q'}\index{Def_q,q'@$\Def_{q,q'}^{\C}$}
\Def_{q,q'}^{\C}:(\sigma_q,\delta_q)-\mathrm{Mod}(\B,\mathrm{C})_{S}\xrightarrow[\qquad]{\sim}
(\sigma_{q'},\delta_{q'})-\mathrm{Mod}(\B,\mathrm{C})_{S}\;,
\end{equation}
the equivalence obtained by composition of the restriction
functor \eqref{forget functor from S to q strong}:
\begin{equation}
\Def_{q,q'}^{\C}:=\mathrm{Res}^{\;S}_{\{q'\}}\circ(\mathrm{Res}^{\;S}_{\{q\}})^{-1}\;.
\end{equation}

2.-- We will call the \emph{$\C$-constant confluence functor}, the
equivalence
\begin{equation}\label{Conf_q} \index{Conf_q@$\Conf_{q}^{\C}:=\Def_{q,1}^{\C}$}
\Conf_{q}^{\C}:=\Def_{q,1}^{\C}
:(\sigma_q,\delta_q)-\mathrm{Mod}(\B,\mathrm{C})_{S}\xrightarrow[\quad]{\sim}
(\sigma_{1},\delta_{1})-\mathrm{Mod}(\B,\mathrm{C})_{S}\;.
\end{equation}

3.-- Suppose that $q\in S^{\circ}$ and $q'\in S$, then we will
call again the \emph{$\C$-constant deformation functor}, denoted
again by
\begin{equation}
\Def_{q,q'}^{\C}:\sigma_q-\mathrm{Mod}(\B,\mathrm{C})_{S}\xrightarrow[\qquad]{}
\sigma_{q'}-\mathrm{Mod}(\B,\mathrm{C})_{S}\;,
\end{equation}
the functor obtained by composition of the restriction functor
\eqref{forget functor from S to q strong}:
$\Def_{q,q'}^{\C}:=\mathrm{Res}^{\;S}_{\{q'\}}\circ(\mathrm{Res}^{\;S}_{\{q\}})^{-1}$.
If $q'\in S^{\circ}$, then $\Def_{q,q'}^{\C}$ is an equivalence.
\end{definition}

It follows from Corollary \ref{forget functor from S to q
strongss}, that if $q,q'\in U$, one has an \emph{equivalence},
again called $\mathrm{Def}_{q,q'}^{\C}$
\begin{equation}
\Def_{q,q'}^{\C}\;:\;(\sigma_q,\delta_q)-\mathrm{Mod}(\B,\mathrm{C})^{\mathrm{an}}_{U}
\xrightarrow[\sim]{}
(\sigma_{q'},\delta_{q'})-\mathrm{Mod}(\B,\mathrm{C})^{\mathrm{an}}_{U}\;.
\end{equation}
The same fact is true for analytic $\sigma-$modules under the
condition $q,q'\notin\bs{\mu}(\Q)$.

\subsubsection{}\label{Confluence depends on C} Notice that the functor
$\mathrm{Res}^{\;S}_{\{q\}}$ does not depend on $\C$, but
$(\mathrm{Res}^{\;S}_{\{q\}})^{-1}$ is a particular section of
$\mathrm{Res}^{\;S}_{\{q\}}$ with values in the category of
objects trivialized by $\C$ (cf. Corollary \eqref{forget functor
from S to q strongss}). Hence $(\mathrm{Res}^{\;S}_{\{q\}})^{-1}$,
$\Conf^{\C}_q$ and $\mathrm{Def}_{q,q'}^{\C}$ actually depend on
$\C$.

\subsubsection{}\label{shhhshhsssuu} According to Definition
\ref{definition of strongly confluent at q} (cf. also
\ref{sigma_q-Mod = sigma-Mod^disc_q} and \ref{sigma_q,delta_q-Mod
= sigma,delta-Mod^disc_q}), if $q\in U\subset U'$, then, by Lemma
\ref{The forget from S to S' is fully faithful} (cf. Section
\ref{fresfresfres}), the following restriction functors are fully
faithful immersions:
\begin{equation}
\begin{array}{llcl}
&\mathrm{Res}^{U'}_U:\sigma-\Mod(\B,\C)_{U'}\qquad %
&\;\; \longrightarrow \;\;&
\sigma-\Mod(\B,\C)_U\\
&\mathrm{Res}^{U'}_U:\sigma-\Mod(\B,\C)_{U'}^{\mathrm{an,const}}%
&\;\; \longrightarrow \;\; &
\sigma-\Mod(\B,\C)_U^{\mathrm{an,const}}\\
\textrm{$\Bigl($resp. }&
\mathrm{Res}^{U'}_U:(\sigma,\delta)-\Mod(\B,\C)_{U'} %
&\;\; \longrightarrow \;\;&
(\sigma,\delta)-\Mod(\B,\C)_U\\
&\mathrm{Res}^{U'}_U:(\sigma,\delta)-\Mod(\B,\C)_{U'}^{\mathrm{an,const}}
&\;\; \longrightarrow \;\;& %
(\sigma,\delta)-\Mod(\B,\C)_U^{\mathrm{an,const}}\qquad\Bigr)\;.
\end{array}
\end{equation}
We can then consider the following diagram in which we
heuristically imagine categories appearing in the first two lines
as the stalks at $q$ of suitable corresponding stacks over
$\Q(X)$:
\begin{equation}
\xymatrix{\bigcup_{U}\;\sigma-\Mod(\B,\C)^{\mathrm{an,const}}_U
\ar@{=}[r]^-{\mathrm{Def.}\;\ref{sigma-an=(sigma,delta)-an}}
\ar[d]_{\bigcup_{U}\;\mathrm{Res}^{\;U}_{\{q\}}}\ar@{}[dr]|{\odot}\quad&
\quad\bigcup_{U}\;(\sigma,\delta)-\Mod(\B,\C)^{\mathrm{an,const}}_U\ar[d]_{\wr}^{\bigcup_{U}\;\mathrm{Res}^{\;U}_{\{q\}}}\\
\bigcup_{U}\;\sigma_q-\Mod(\B,\C)_U\quad\ar[d]_{i_\sigma}\ar@{}[dr]|{\odot}&\quad
\bigcup_{U}\;(\sigma_q,\delta_q)-\Mod(\B,\C)_U\ar[l]^{\textrm{Forget
}\delta_q}\ar[d]^{i_{(\sigma,\delta)}}\\
\sigma_q-\Mod(\B)\quad&\quad(\sigma_q,\delta_q)-\Mod(\B)
\ar[l]^{\textrm{Forget }\delta_q}}
\end{equation}
where $U$ runs over the set of open neighborhoods of $q$, and where
$i_\sigma$ and $i_{(\sigma,\delta)}$ are the trivial inclusions of
full sub-categories. In the sequel we will study the full
subcategory of $\sigma_q-\Mod(\B)$ (resp.
$(\sigma_q,\delta_q)-\Mod(\B)$) formed by \emph{Taylor admissible
objects}, this category is contained in the essential image of
$i_\sigma$ (resp. $i_{(\sigma,\delta)}$) (see Th. \ref{main
theorem}). In this case we will obtain an analogous diagram (see
Cor. \ref{compare with})) in which $i_{(\sigma,\delta)}$ is an
equivalence (for all $q\in U$), and $i_\sigma$ is an equivalence
only if $q$ is not a root of unity.

If $q$ is not a root of unity, then all the arrows of this diagram
will be equivalences, hence giving $\delta_q$ is superfluous.
If $q$ is a root of unity, then the right hand side vertical
arrows will be equivalences, while the arrow on the left hand side will not.
 In this last case the $q-$tangent operator is necessary to
``\emph{preserve the information in the neighborhood of $q$}''.
In this case  the good notion of ``\emph{stalk at $q$}'' of an analytic
$\sigma$-module is  the notion of
$(\sigma_q,\delta_q)$-module and not simply that  of
$\sigma_q$-module.

One may have the feeling that the functor
``$\mathrm{Forget}\;\delta_q$'' contains  ``information'' if $q$
is a root of unity, but we will see (Prop. \ref{frobenius are
trivial in the root of 1}) that, if $\B=\R_K$ or if $\B=\Hd_K$,
then this functor sends every $(\sigma,\delta)-$module with
Frobenius structure into a direct sum of copies of the unit object.

\subsubsection{Dependence on $\C$.}\label{C_1 subset C_2 ==> Then
Def^C_1=Def^C_2} Let $\C_1\subseteq \C_2$ be two algebras as
above.  Then clearly $\mathrm{Def}_{q,q'}^{\C_2}$ extends
$\mathrm{Def}_{q,q'}^{\C_1}$ to the larger category of modules
trivialized by $\C_2$. One of the main problems of the theory is
that, if there are no inclusions between $\C_1$ and $\C_2$, then
it is not clear whether there exists a discrete $\sigma$-algebra
(resp. $(\sigma,\delta)$-algebra) $\C_3$ containing both $\C_1$
and $\C_2$. For this reason, if the same object is trivialized by
$\C_1$, and also by $\C_2$, it is not clear whether its
deformations with respect to $\C_1$ and $\C_2$ are equal. We will
encounter this problem in section \ref{quasi unipotence of q-diff
- section}.

\section{Taylor solutions}

\label{Taylor solutions--}

In this section $\B=\H_K(X)$, for some affinoid
$X=\mathrm{D}^+(c_0,R_0)-\cup_{i=1}^n\mathrm{D}^-(c_i,R_i)$, and
$S=\{q\}\in\Q(\H_K(X))\subseteq\{q\in K\;|\;|q|=1\}$ is reduced to
a point. Let $(\Omega,|.|)/(K,|.|)$ be an arbitrary extension of
complete valued fields. Let $c\in X(\Omega)$ and let
$\rho_{c,X}>0$ be the largest real number such that
$\mathrm{D}_{\Omega'}^-(c,\rho_{c,X})\subseteq X(\Omega')$, for
all complete valued field extensions $(\Omega',|.|)/(\Omega,|.|)$.
One has
\begin{equation}\label{rho_c,A-y-t-}\index{rho_c,X@$\rho_{c,X}$}
\rho_{c,X}=\min(R_0, |c-c_1|, |c-c_2|, \cdots, |c-c_n|).
\end{equation}
Notice that $c$ can be equal to a generic point (cf. Definition
\ref{definition of generic point}). We want to find solutions of
$q-$difference equations converging in a disc centered at $c$,
i.e. matrix solutions in the form \eqref{sigma(Y)=YA}, with values
in the $\sigma_q-$algebra $\mathrm{C}:=\a_K(c,R)$, for some
$0<R\leq\rho_{c,X}$.

\subsection{The $q-$algebras $\Omega\{T-c\}_{q,R}$ and $\Omega[\![T-c]\!]_q$}
Unless we explicitly state the contrary, we will not assume that
$q\notin\bs{\mu}(\Q)$. The following results generalize the
analogous constructions of \cite{DV-Dwork} to the case of a root of
unity.

\begin{lemma}\label{|q-1||c|<R-yt}
Let $0 < R \leq \rho_{c,X}$. The algebra $\a_\Omega(c,R)$ is an
$\H_\Omega(X)$-discrete $\sigma-$algebra over $S=\{q\}$, if and
only if both of the following conditions hold:
\begin{equation}\label{.FF.hyp}
|q-1||c|<R\;,\qquad\mathrm{and}\qquad |q|=1\;.
\end{equation}
\end{lemma}

\begin{definition}\label{def of
(T-c)_q,n ..g.G.g.g} Let $q\in K^{\times}$ be an arbitrary number.
Following \cite{DV-Dwork} and \cite{An-DV} we set
\begin{eqnarray}\label{basis:;}\index{T-c_q,n@$(T-c)_{q,n}$}
(T-c)_{q,n}&:=&(T-c)(T-qc)(T-q^2c)\cdots(T-q^{n-1}c)\;, \\
\protect{[n]}_q&:=&1+q+q^2+\cdots+q^{n-1}\;,\\
\protect{[n]}_q^{!}&:=&\frac{(q-1)(q^2-1)(q^3-1)\cdots(q^n-1)}{(q-1)^n}\;.
\end{eqnarray}
\end{definition}

\subsubsection{$q$-binomial.}For all $q\in K^{\times}$, we define the $q$-binomial
$\binom{n}{i}_q$ by the relation
\begin{equation}
(1-T)(1-qT)\cdots(1-q^{n-1}T)=\sum_{i=0}^{n}(-1)^i\binom{n}{i}_qq^{\frac{i(i-1)}{2}}T^i\;,
\end{equation}
where, if $i=0$, the symbol $q^{\frac{i(i-1)}{2}}$ is by
definition equal to $1$. This extends the definition given in
\cite{DV-Dwork} (cf. eq. \eqref{Di vizio definition of
q-factorials} below) to the case of root of unity. If $1\leq i\leq
n-1$, by induction one has
\begin{equation}
\binom{n}{i}_q=\binom{n-1}{i-1}_q+q^i\binom{n-i}{i}_q=
q^{n-i}\binom{n-1}{i-1}_q+\binom{n-i}{i}_q\;.
\end{equation}
If $q$ is not a root of unity, then one can write
\begin{equation}\label{Di vizio definition of q-factorials}
\binom{n}{i}_q=\frac{[n]_q\cdot[n-1]_q\cdots[n-i+1]_q}{[i]^!_q}\;.
\end{equation}

If $q$ is an $m-$th root of $1$, then $[n]_q^!=0$, for all $n\geq
m$. The family $\{(T-c)_{q,n}\}_{n\geq 0}$ is adapted to the
$q-$derivation
\begin{equation}\label{d_q}\index{d_q@$d_q$}
d_q:=\frac{\sigma_q-1}{(q-1)T}=\frac{\Delta_q}{T}
\end{equation} in the sense
that for all $n\geq 1$ one has
$d_q((T-c)_{q,n})=[n]_q\cdot(T-c)_{q,n-1}$. One has always the
relation $d_q(fg)=\sigma_q(f)d_q(g)+d_q(f)g$. More generally our
definition of $q$-binomials allow us to generalize the proof of
\cite[Lemma 1.2, (1.2.2)]{DV-Dwork} to the case of a root of
unity. We obtain the formula
\begin{equation}\label{d_q^n(fg)(T)=...}
d_q^n(fg)(T)=\sum_{i=0}^n\binom{n}{i}_{\!\!q}d_q^{n-i}(f)(q^iT)d_q^i(g)(T)\;.
\end{equation}

\subsubsection{}The following Lemma extends \cite[Section
1.3]{DV-Dwork} to the case of a root of unity.
\begin{lemma}\label{Radius twisted ...}
Let $(\Omega,|.|)/(K,|.|)$ be a complete extension of valued
fields. Let $|q-1||c|<R$, $|q|=1$, and let $f(T)=\sum_{n\geq
0}a_n(T-c)^n\in\a_\Omega(c,R)$. Then:
\begin{enumerate}
\item $f(T)$ can be written uniquely as the following series of
functions:
\begin{equation}
f(T)=\sum_{n\geq 0}\widetilde{a}_n(T-c)_{q,n}\in\a_\Omega(c,R)\;,
\end{equation}
with $\widetilde{a}_n\in \Omega$ satisfying
$\sup_n|\widetilde{a}_n|\rho^n<\infty$, for all $\rho<R$;

\item for all $|q-1||c|<\rho<R$ one has $|f(T)|_{(c,\rho)}=
\sup_{n\geq 0}|a_n|\rho^n=\sup_{n\geq 0}|\widetilde{a}_n|\rho^n$;%

\item  one has
$Ray(f(T),c)=\liminf_n|a_n|^{-1/n}=\liminf_n|\widetilde{a}_n|^{-1/n}$;%

\item if moreover $q\notin\bs{\mu}(\Q)$, then one has the so
called $q$-Taylor expansion (cf. \cite{DV-Dwork}):
\begin{equation}\label{Twisted Taylor formula}
f(T)=\sum_{n\geq 0}d_q^n(f)(c)\frac{(T-c)_{q,n}}{[n]^!_q}\;.
\end{equation}
\end{enumerate}
\end{lemma}
\begin{proof} Since $\a_\Omega(c,R)=\varprojlim_{r\to
R^-}\H_\Omega(\mathrm{D}^+(c,r))$,  we need only prove the
proposition for $\H_\Omega(\mathrm{D}^+(c,r))$, with
$|q-1||c|<r<R$. We recall that a series of functions $\sum_{n\geq
0} f_n$, $f_n\in\H_K(\mathrm{D}^+(c,r))$ converges to a function
$f\in\H_K(\mathrm{D}^+(c,r))$ if and only if
$\lim_n|f_n|_{(c,r)}=0$. Writing $(T-q^ic)=(1-q^i)c+(T-c)$, one
sees easily that $(T-c)_{q,n}=\sum_{i=0}^n
\widetilde{b}_{n,i}(T-c)^i$, with $\widetilde{b}_{i,j}$ satisfying
i) $\widetilde{b}_{0,0}=1$, ii) $\widetilde{b}_{0,i}=0$ $\forall
\; i\geq 1$, iii) $\widetilde{b}_{n,n}=1$ $\forall \; n\geq 0$,
iv) $\widetilde{b}_{n,i}=0$ $\forall \; i > n$, and v) for all
$0\leq i<n$:
\begin{equation}
\widetilde{b}_{n,i}=c^{n-i}\cdot\sum_{0\leq k_1<\cdots<k_{n-i}\leq
n-1}(1-q^{k_1}) (1-q^{k_2})\cdots (1-q^{k_{n-i}})\;.
\end{equation}
In other words
$[1,(T-c)_{q,1},(T-c)_{q,2},\ldots,(T-c)_{q,n}]^{t}=\widetilde{B}\cdot[1,(T-c),(T-c)^{2},\ldots,(T-c)^{n}]^{t}$
where $\widetilde{B}=(\widetilde{b}_{n,i})_{n,i=0,\ldots,n}$ is an
$(n+1) \times (n+1)$ lower triangular matrix satisfying
i),ii),iii),iv),v). Since $|q^i-1|\leq |q-1|$, one has also the
property vi) $|\widetilde{b}_{n,i}|\leq (|q-1||c|)^{n-i}<r^{n-i}$,
for all $0\leq i<n$. Hence for all $n\geq 0$, one has
$(T-c)_{q,n}=(T-c)^n+g_n(T)$, with $|g_n(T)|_{(c,r)}<r^n$, so
$|(T-c)_{q,n}|_{(c,r)}= |(T-c)^{n}|_{(c,r)}=r^{n}$. It is easy to
prove that also the matrix
$B:=\widetilde{B}^{-1}=(b_{n,i})_{n,i=0,\ldots,n}$ satisfies the
properties i),ii),iii),iv),v),vi).
Consider now $f(T)=\sum_{n\geq 0}a_n(T-c)^n$. Writing
$f_m(T):=\sum_{n=0}^{m}a_n(T-c)^{n}=\sum_{n=0}^ma_n\sum_{i=0}^n
b_{n,i}(T-c)_{q,i}$ and rearranging terms one finds $f_m(T) = \sum_{n =
0}^m \widetilde{a}_{n,m}(T-c)_{q,n}$, with
$\widetilde{a}_{n,m}=\sum_{k=0}^{m-n}a_{n+k}b_{n+k,n}$. By
property vi) and by the assumption that $\lim_{n}|a_n|r^n=0$ the sum
$\widetilde{a}_{n}:=\sum_{k\geq 0}a_{n+k}b_{n+k,n}$ converges in
$\Omega$. Moreover
\begin{equation}\label{inneq -tret}
|\widetilde{a}_n|r^n \;\; \leq \;\; \max_{k\geq
0}|a_{n+k}||b_{n+k,n}| \cdot r^{n} \;\; \leq \;\; \max_{k\geq
0}|a_{n+k}|r^{n+k}\;.
\end{equation}
This proves that $\lim_{n} |\widetilde{a}_n| r^n = 0 $, and hence
that the series of functions $\sum_{n\geq
0}\widetilde{a}_n(T-c)_{q,n}$ is convergent in
$\H_\Omega(\mathrm{D}^+(c,r))$. If
$f_m^0(T):=\sum_{n=0}^m\widetilde{a}_n(T-c)_{q,n}$, one sees that
$|f_m^0-f_m|_{(c,r)}\leq \sup_{k\geq 0}|a_{m+k}|r^{m+k}$ which
tends to $0$, so $\lim_mf_m^0(T)=\lim_mf_m(T)=f(T)$ in
$\H_\Omega(\mathrm{D}^+(c,r))$.
Now the inequality \eqref{inneq -tret} shows that $\max_{n\geq
0}|\widetilde{a}_n|r^n \leq \max_{n\geq 0}|a_n|r^n$, and a
symmetric argument using the matrix $\widetilde{B}$ instead of $B$
proves the opposite inequality so $\max_{n\geq
0}|\widetilde{a}_n|r^n = \max_{n\geq 0}|a_n|r^n=|f(T)|_{(c,r)}$.
This last equality shows the uniqueness of the coefficients
$\{\widetilde{a}_n\}_n$ since if $\sum_{n\geq
0}\widetilde{a}_n(T-c)_{q,n}=\sum_{n\geq
0}\widetilde{a}_n'(T-c)_{q,n}$, then $\sum_{n\geq
0}(\widetilde{a}_n-\widetilde{a}_n')(T-c)_{q,n}=0$, and hence
$\sup_n(|\widetilde{a}_n-\widetilde{a}_n'|r^n)=0$, so that
$\widetilde{a}_n=\widetilde{a}_n'$, for all $n\geq 0$. Clearly the
radius of convergence of $f(T)$ is equal to both $\sup_{n\geq
0}\{r\geq 0\;|\;|a_n|r^n\textrm{ is bounded}\}$ and $\sup_{n\geq
0}\{r\geq 0\;|\;|\widetilde{a}_n|r^n\textrm{ is bounded}\}$.
Hence, by classical arguments on the radius of convergence, one
has
$Ray(f(T),c)=\liminf_n|a_n|^{-1/n}=\liminf_n|\widetilde{a}_n|^{-1/n}$.
The assertion iv) is proved in \cite{DV-Dwork}.
\end{proof}

\begin{remark}\label{multipliation rule}
If $f(T)=\sum_{n\geq 0}f_n(T-c)_{q,n}$, and if $g(T)=\sum_{n\geq
0}g_n(T-c)_{q,n}$, then $f(T)g(T)=\sum_{n\geq 0}h_n(T-c)_{q,n}$,
where $h_n=h_n(q;c;f_0,\ldots,f_n;g_0,\ldots,g_n)$ is a polynomial
in $\{q,c,f_0,\ldots,f_n,g_0,\ldots,g_n\}$. Indeed one has
$(T-c)_{q,n}\cdot(T-c)_{q,m}=\sum_{k=\max(n,m)}^{n+m}\alpha_k^{(n,m)}(T-c)_{q,k}$,
with $\alpha_k^{(n,m)}=\alpha_k^{(n,m)}(q,c)\in\Omega$. This also shows
 that if $v_{q,c}(f):=\min\{n\;|\; f_n\neq 0\}$, then one has
\begin{equation}
v_{q,c}(fg)\geq \max(v_{q,c}(f),v_{q,c}(g))\;.
\end{equation}

If moreover $q\notin\bs{\mu}(\Q)$, then, by using equation
\eqref{d_q^n(fg)(T)=...} and \eqref{Twisted Taylor formula}, one
has
\begin{equation}
h_n=\sum_{j=0}^n\sum_{s=0}^j\frac{[n]^!_q[j]^!_q[s+n-j]^!_q}{([s]^!_q)^{2}[n-j]^!_q}\cdot
q^{\frac{s(s-1)}{2}}(q-1)^{s}c^sf_{s+n-j}g_j\;.
\end{equation}
\end{remark}

\subsubsection{The algebras $\Omega[\![T-c]\!]_{q}$, and $\Omega\{T-c\}_{q,R}$.}

\begin{definition}\index{OmegaT-c_q,R@$\Omega[[T-c]]_{q}$, $\Omega\{T-c\}_{q,R}$}
For all $q\in\Q(X)$ we set
\begin{eqnarray}\label{Omega((T-c))_q}
\Omega[\![T-c]\!]_{q}&:=&\{\sum_{n\geq
0}f_n(T-c)_{q,n}\;|\;f_n\in\Omega\}\;.\\
\qquad\Omega\{T-c\}_{q,R}&:=&\{\;\sum_{n\geq
0}f_n(T-c)_{q,n}\;|\;f_n\in\Omega,\;\liminf_n|f_n|^{-1/n}\geq
R\;\}\;.
\end{eqnarray}
We define a multiplication on $\Omega[\![T-c]\!]_{q}$ and
$\Omega\{T-c\}_{q,R}$ by the rule given in Remark
\ref{multipliation rule}.
\end{definition}

\begin{lemma}
$\Omega[\![T-c]\!]_{q}$ and $\Omega\{T-c\}_{q,R}$ are commutative
$\Omega$-algebras, $\forall\;q\in\Q$.
\end{lemma}
\begin{proof} We prove only the associativity, the others verifications are similar. We have
to prove that $(fg)h=f(gh)$. By Lemma \ref{Radius twisted ...} the
assertion is proved if $f,g,h\in\Omega\{T-c\}_{q,R}$, with
$|q-1||c|<R$, since in this case
$\Omega\{T-c\}_{q,R}\cong\a_{\Omega}(c,R)$. On the other hand one
can assume that $f,g,h$ are polynomials since, by Remark
\ref{multipliation rule}, the $n$-th coefficient of $(fg)h$ and of
$f(gh)$ is a polynomial in $q$ and in the first $n$ coefficients of
$f,g,h$.
\end{proof}

\begin{remark}
If there exists a (smallest) integer $k_0$ such that
$|\overline{q}^{k_0}-1||c|<R$, then one shows that
$\Omega\{T-c\}_{q,R}=\prod_{i=0}^{k_0\!-1}\a_{\Omega}(q^ic,\widetilde{R})$,
where $\widetilde{R}$ depends explicitly on $R$, $c$, and $q$ (cf.
\cite[15.3]{DV-Dwork}). In this case $\Omega\{T-c\}_{q,R}$ is not
a domain and hence is not a $\H_\Omega(X)-$discrete
$\sigma-$algebra over $S=\{q\}$.
\end{remark}

\begin{remark}
If $x,y$ are variables, then $\Omega[\![x-y]\!]_q$ is not an
algebra, but merely a vector space. Indeed the multiplication law
involves $y$ in the coefficients ``$h_n$'' of Remark
\ref{multipliation rule}. This minor mistake occurs occasionally
in \cite{DV-Dwork}, but it is an irrelevant inaccuracy and does
not jeopardize any proposition of \cite{DV-Dwork}. The matrix
$Y(x,y)$ always seems to be used there under the assumption
\eqref{GOOD condition for Y(x,y)}.
\end{remark}
\subsection{q-invariant Affinoid}

 \label{k_0} Let $|q|=1$, $q\in K$. Let $X := \mathrm{D}^+(c_0,R_0)-
\cup_{i=1}^n \mathrm{D}^-(c_i,R_i)$,
$c_1,\ldots,c_n\in\mathrm{D}_K^+(c_0,R_0)$, $c_0\in K$, be a
$K$-affinoid. Then $X$ is $q$-invariant if and only if
$|q-1||c_0|<R_0$, and the map $x\mapsto q x$ permutes the family
of disks $\{\;\mathrm{D}^-(c_i,R_i)\;\}_{i=1,\ldots,n}$. This
happens if and only if for all $i=1,\ldots, n$ there exists (a
smallest) $k_i\geq 1$, such that $|q^{k_i}-1||c_i|<R_i$, and
moreover the family of disks
$\{\;\mathrm{D}^-(q^kc_i,R_i)\;\}_{k=1,\ldots,k_i}$ is finite and
contained in $\{\;\mathrm{D}^-(c_i,R_i)\;\}_{i=1,\ldots,n}$. If
$k_0$ is the minimum common multiple of the $k_i$'s, then
$x\mapsto q^{k_0}x$ leaves  every disk globally fixed and, by
Lemmas \ref{|q-1||c|<R-yt} and \ref{Radius twisted ...},  one has
\begin{equation}\label{d_q(f)leq r_A^-1 f, fifi}
\|d_{q^{k_0}}(f)\|_X\;\leq \;r_X^{-1}\|f\|_X\;,
\end{equation}
for all $f\in\H(X)$ (cf. Lemma \ref{|f'| leq |f| / r_A}). Indeed
the by Mittag-Lefler decomposition \cite{Ch-Ro}, we  reduce to
showing that every series $f=\sum_{j\leq -1}a_j(T-c_i)^{j}$, such
that $|a_j|R_i^j$ tends to zero, satisfies
$|d_{q^{k_0}}(f)|_{(c_i,R_i)}\leq R_i^{-1}\cdot|f|_{(c_i,R_i)}$,
and this is true by Lemma \ref{Radius twisted ...}.

Such a bound does not exist for $d_q$ itself. One can easily
construct counterexamples via the Mittag-Leffler decomposition.

\subsection{The generic Taylor solution}
We recall the definition of the classical Taylor solution of a
differential equation
\begin{definition}
Let $\delta_1-G(1,T)$, be a differential equation. Let
$G_{[n]}(T)$ be the matrix of $(d/dx)^{n}$. We set
\begin{equation}\label{Y_G-+-} \index{Y_G(1,T)(x,y)@$Y_{G(1,T)}(x,y)$}
Y_{G(1,T)}(x,y):=\sum_{n\geq 0}G_{[n]}(y)\frac{(x-y)^n}{n!}\;.
\end{equation}
\end{definition}
By induction on the rule $G_{[n+1]}=G_{[n]}'+G_{[n]}G_{[1]}$, one
finds $\|G_{[n]}\|_X\leq\max(\|G_{[1]}\|_X,r_X^{-1})^n$, hence
\begin{equation}\label{R_c in the differential case}
Ray(Y_G(T,c),c)=\liminf_{n}\Bigl(\frac{|G_{[n]}(c)|_\Omega}{|n!|}\Bigr)^{-1/n}\geq
\frac{|p|^{\frac{1}{p-1}}}{\max(r_X^{-1},\|G_{[1]}\|_X)}\;.
\end{equation}
In other words $Y_G(x,y)$ is an analytic element over a
neighborhood $\mathcal{U}_R$ of the diagonal of the type
\begin{equation}\label{U_R}\index{U_R@$\mathcal{U}_R:=\{(x,y)\;|\;|x-y|<R\;\}$}
\mathcal{U}_R:=\{(x,y)\in X\times X\;|\;|x-y|<R\}\;,
\end{equation}
for some $R>0$.
\begin{lemma} One has
$Y_{G}(x,x)=\mathrm{Id}$, and for all $(x,y)\in\mathcal{U}_R$:
\begin{eqnarray}
d/dy\; (Y_G(x,y))&=&-Y_G(x,y)\cdot G_{[1]}(y)\;,\\
Y_G(x,y)^{-1}&=&Y_G(y,x)\;,\\
Y_G(x,y)\cdot Y_G(y,z)&=&Y_G(x,z)\;,\\
d/dx\;(Y_G(x,y))&=&G_{[1]}(x)\cdot Y_G(x,y)\;.
\end{eqnarray}
\end{lemma}
\begin{proof} See \cite[p.137]{Astx} (cf. Lemma
\ref{Y(x,y)Y(y,z)=Y(x,z)}). The proof is analogous to that of
Lemma \ref{Y(x,y)Y(y,z)=Y(x,z)}.\end{proof}

\begin{definition}\label{Definition of the q-Taylor solution}
Let $q\in\Q-\bs{\mu}(\Q)$. Consider the $q-$difference equation
\begin{equation}\label{Y(qT)=AY--}
\sigma_q(Y)=A(q,T)\cdot Y\;,\qquad A(q,T)\in GL_n(\H_K(X))\;.
\end{equation}
Let $H_n$ be defined by $d_q^n(Y)=H_n\cdot Y$. We formally set
\begin{equation}\label{Y_A(T,c)}\index{Y_A(x,y)@$Y_{A(q,T)}(x,y)$}
Y_{A(q,T)}(x,y)\;=\;\sum_{n\geq
0}H_n(y)\frac{(x-y)_{q,n}}{[n]^!_q}\;.
\end{equation}
We will omit the index $A(q,T)$ if no confusion is possible.
Observe that $Y_{A(q,T)}(x,y)$ is a symbol and does not
necessarily define a convergent function.
\end{definition}

\begin{example}
With the notations of example \ref{Morphisms between sigma modules
-examples}, the generic Taylor solution of the equations
$\sigma_q(Y)=A(q,x)Y$, $\sigma_q(Y)=\widetilde{A}(q,x)Y$, are
$Y_{A(q,x)}(x,y)=\exp(\pi (x-y))$ and
$Y_{\widetilde{A}(q,x)}(x,y)=\exp(\pi q (x-y))$ respectively.
Notice that $Y_{A(q,x)}(x,y)$ is \emph{constant} with $q$.
\end{example}

\begin{definition}
For all (not necessarily bounded nor multiplicative) semi-norms
$|.|_*$ on $\H_K(X)$ extending the absolute value of $K$ we set
\begin{equation}\label{R_c--}
Ray(Y_{A(q,T)}(x,y),
|.|_*)\;:=\;\liminf_n(|H_n(y)|_*/|[n]^!_q|)^{-1/n}\;.
\end{equation}
If $Y_{A(q,T)}(x,y)$ is a convergent function on some neighborhood
of the diagonal of $X\times X$, then, for
$|f(T)|_*:=|f(c)|_\Omega$, $c\in X(\Omega)$, one finds Definition
\ref{radii}, namely $Ray(Y(x,y),|.|_c)=Ray(Y(x,c),c)$. In this
case we will write $Ray(Y(x,y), c):=Ray(Y(x,y), |.|_c)$ (cf.
Section \ref{exemple of norms}). If $X'\subseteq X$ is a
sub-affinoid we simply write $Ray(Y(x,y), X'):=Ray(Y(x,y),
\|.\|_{X'})$.
\end{definition}

\subsection{Transfer Principle}

\label{transfert}
 As in the differential setting, if
$X':=\mathrm{D}^+(c'_0,R'_0)-\cup_{i=1}^n\mathrm{D}^-(c_i',R_i')\subseteq
X$ is a $q-$invariant sub-affinoid, such that every disk
$\mathrm{D}^-(c_i',R_i')$ is also $q-$invariant, then the estimate
 \eqref{d_q(f)leq r_A^-1 f, fifi} holds (cf. Remark
\ref{k_0-example}).  Then, by induction on the rule
$H_{n+1}=d_q(H_{n})+\sigma_q(H_{n})H_{1}$, one shows that
$\|H_{n}\|_{X'}\leq\max(\|H_{1}\|_{X'},r_{X'}^{-1})^n$, hence
\begin{eqnarray}\label{R_A estimation}
Ray(Y(x,y),X')&:=&\liminf_n(\|H_n\|_{X'}/[n]^!_q)^{-1/n}
=\min_{c\in X'(\Omega)}Ray(Y(x,y),c)\nonumber\\
&\geq&\frac{\liminf_n([n]^!_q)^{1/n}}{\max(r_{X'}^{-1},\|H_1\|_{X'})}\;,
\end{eqnarray}
where $(\Omega,|.|)/(K,|.|)$ is an algebraically closed extension
of complete valued field such that $|\Omega|=\mathbb{R}_{\geq 0}$.
Observe that the second equality follows by the fact that
$\|.\|_X=\sup_{c\in X(\Omega)}|.|_c$. In particular if
$X'=\mathrm{D}^+(c,\rho)\subseteq X$, with
$|q-1||c|<\rho\leq\rho_{c,X}$, is a $q-$invariant disk, then
$Ray(Y(x,y),c)$ is greater than or equal to
\begin{equation}\label{R_c geq R_D...}
Ray(Y(x,y),\mathrm{D}^+(c,\rho))=\min_{c'\in\;
\mathrm{D}^+_{\Omega}(c,\rho)}Ray(Y(x,y),c')\;\geq\;
\frac{\liminf_n([n]^!_q)^{1/n}}{\max(\rho^{-1},|H_1|_{(c,\rho)})}\;.
\end{equation}

Notice that if $|q-1||c|<R_c:=Ray(Y(x,y),c)$, then $Y(x,c)\in
M_n(\a_\Omega(c,R_c))$, but $Y(x,c)$ is invertible only in
$GL_n(\a_{\Omega}(c,\widetilde{R}))$, with
$\widetilde{R}:=\min(\rho_{c,X},Ray(Y(x,y),c))$ (cf. Lemmas
\ref{det neq 0} and \ref{Y(x,y)Y(y,z)=Y(x,z)}).

\subsection{Properties of the generic Taylor solution}

\label{Y defines an analytic function on U_R} The formal matrix
solution $Y_A(x,y)$ is not always a function in a neighborhood of
type $\mathcal{U}_R$ of the diagonal of $X\times X$. But if for
all $c\in X(K^{\mathrm{alg}})$ one has $|q-1||c|<R\leq
\min(\rho_{c,X},Ray(Y(x,y),c))$, then, by Lemma \ref{Radius
twisted ...}, and by the  Transfer Principle (cf. Equation
\eqref{R_c geq R_D...}), $Y_A(x,y)$ actually defines an
\emph{invertible} function on $\mathcal{U}_R$ (cf. Lemmas \ref{det
neq 0} and \ref{Y(x,y)Y(y,z)=Y(x,z)}). If
$X=\mathrm{D}^+(c_0,R_0)-\cup_{i=1}^n\mathrm{D}^-(c_i,R_i)$, the
condition $|q-1||c|<R\leq \min(\rho_{c,X},Ray(Y(x,y),c))$, for all
$c\in X(K^{\mathrm{alg}})$, implies
\begin{equation}\label{FRREISTUGHY}
|q-1|\sup(R_0,|c_0|)\;=\;|q-1|\max_{c\in X}|c|\; <\; R\; \leq\;
\min_{c\in X}\rho_{c,X}\;=\;\min(R_0,\ldots,R_n)\;=\;r_X\;.
\end{equation}
In particular, since $r_X=\min(R_0,\ldots,R_n)\leq
\sup(|c_0|,R_0)$, this is possible only if
\begin{equation}
|q-1|<1\;,\quad\textrm{i.e.\; if\; } q\in\Q_1(X)\;.
\end{equation}

\begin{hypothesis}
From now on, without explicit mention to the contrary, we will
assume that
\begin{equation}
q\in\Q_1(X)\;.
\end{equation}
\end{hypothesis}

\begin{lemma}\label{det neq 0} Let
$q\in\Q_1(X)-\bs{\mu}(\Q_1(X))$. Let $f(x,y)$ be an analytic
function in a neighborhood of type $\mathcal{U}_R\subset X\times
X$ of the diagonal of $X\times X$.  Assume that\footnote{i.e.
assume that $|q-1|\max(|c_0|, R_0) <R\leq\rho_{c,X}$ for all $c\in
X(K^{\mathrm{alg}})$.}
\begin{equation}
|q-1|\max(|c_0|, R_0) < R \leq r_X\;.
\end{equation}
If moreover $f(x,y)$ satisfies $f(x,qy)=a(y)\cdot f(x,y)$, with
$a(y)\in \H_K(X)^{\times}$, then $f(x,y)$ is invertible.
\end{lemma}
\begin{proof} Since $f$ is an analytic function, it is sufficient to prove
that $f$ has no zeros in $\mathcal{U}_R$. We need only show
that for all $c\in X(\Omega)$, the function $g_c(y):= f(c,y)$ has
no zeros in $\mathrm{D}^-(c,R)$. One has $d_q(g_c(y))=h(y)\cdot
g_c(y)$, with $h(y)=\frac{a(y)-1}{(q-1)y}$. Assume that
$g_c(\widetilde{c})=0$, for some
$\widetilde{c}\in\mathrm{D}^-(c,R)=\mathrm{D}^-(\widetilde{c},R)$,
then, by Lemma \ref{Radius twisted ...}, $g_c(y)=\sum_{n\geq
0}a_k(y-\widetilde{c})_{q,n}$, with $a_0=0$. Since
$q\notin\bs{\mu}(\Q)$, we have $[n]_qa_n=0$ if and only if $a_n=0$.
Hence, by Remark \ref{multipliation rule} one has
$v_{q,\widetilde{c}}(d_q(g_c))=v_{q,\widetilde{c}}(g_c)-1$.
 On the other hand, $v_{q,\widetilde{c}}(hg_c)\geq
v_{q,\widetilde{c}}(g_c)$, which  contradicts
$d_q(g_c)=hg_c$. \end{proof}

\begin{lemma}\label{Y(x,y)Y(y,z)=Y(x,z)}
Let $q\in\Q_1(X)-\bs{\mu}(Q_1(X))$, and let
\begin{equation}
\begin{array}{rclcrcl}
\sigma_q^x&:&f(x,y)\mapsto f(qx,y)\;,&&\sigma_q^y&:&f(x,y)\mapsto f(x,qy)\;,\\
&&&&&&\\
d_q^x&:=&\frac{\sigma_q^x-1}{(q-1)x}\;,&&d_q^y&:=&\frac{\sigma_q^y-1}{(q-1)y}\;.
\end{array}
\end{equation}
Suppose that $Y_A(x,y)$ converges on $\mathcal{U}_R$, with (cf.
Section \ref{Y defines an analytic function on U_R})
\begin{equation}\label{GOOD condition for Y(x,y)}
|q-1|\max(|c_0|, R_0) \;< \;R\; \leq\; r_X\;.
\end{equation}
Then $Y_A(x,y)$ is invertible on $\mathcal{U}_R$ and satisfies
$Y_A(x,x)=\mathrm{Id}$ and:
\begin{eqnarray}
d_q^y\; Y_A(x,y)&=&-\sigma_q^y(Y_A(x,y))\cdot H_1(y)\;,\label{eq-2-}\\
\sigma_q^y\; Y_A(x,y)&=&Y_A(x,y)\cdot A(q,y)^{-1}\;,\label{eq-3-}\\
Y_A(x,y)\cdot Y_A(y,z)&=&Y_A(x,z)\;,\label{eq-4-}\\
Y_A(x,y)^{-1}&=&Y_A(y,x)\;,\label{eq-4.5-}\\
d_q^x\;Y_A(x,y)&=&H_1(x)\cdot Y_A(x,y)\;.\label{eq-5-}\\
\sigma_q^x\;Y_A(x,y)&=&A(q,x)\cdot Y_A(x,y)\;.\label{eq-6-}
\end{eqnarray}
\end{lemma}
\begin{proof} The relation $Y(x,x)=\mathrm{Id}$ is evident, while \eqref{eq-2-} is
easy to compute explicitly, and is equivalent to \eqref{eq-3-}.
Since $Y(x,y)$ converges on $\mathcal{U}_R$, equation
\eqref{eq-3-} implies that the determinant $d(x,y)$ of $Y(x,y)$
satisfies $d(x,qy)=a(y)d(x,y)$, with
$a(y)=\det(A(q,y)^{-1})\in\H_K(X)^{\times}$. By Lemma \ref{det neq
0}, $d(x,y)$ is then invertible on $\mathcal{U}_R$, and hence also
$Y(x,y)$ is invertible. By equation \eqref{d_q^n(fg)(T)=...}, and
since $q\notin\bs{\mu}(\Q)$, the relation
$d_q^{y}(Y(x,y)Y(x,y)^{-1})=0$ gives
\begin{equation}\label{frtgmlk b}
d_q^{y}(Y(x,y)^{-1})=-\sigma_q^{y}(Y(x,y)^{-1})\cdot
d_q^{y}(Y(x,y))\cdot Y(x,y)^{-1}\;.
\end{equation}
Hence, for all $x,y,z$ such that $|x-y|,|z-y|<R$, the relation
\eqref{frtgmlk b}, together with  relation \eqref{eq-2-}, give
$d_q^{y}(Y(x,y)\cdot Y(z,y)^{-1})=0$. 
Since $q\notin\bs{\mu}(\Q)$, this implies, by Lemma \ref{Radius
twisted ...}, that the function $Y(x,y)Y(z,y)^{-1}$ does not
depend on $y$. Specializing for $y=x$, and $y=z$, one finds
$Y(x,z)=Y(z,x)^{-1}$, and $Y(x,y)\cdot Y(y,z)=Y(x,z)$. Then, by
the above expression for $d^x_q(Y(y,x)^{-1})=d^x_q(Y(x,y))$, the
relations \eqref{eq-5-} and \eqref{eq-6-} follow from
\eqref{eq-4.5-} and \eqref{eq-2-}. \end{proof}

\subsubsection{The case $|q-1|=1$, $|q|=1$.}
If for a $c\in X$ one has $|q-1||c|\geq Ray(Y_{A(q,T)}(x,y),c)$,
then Lemma \ref{Y(x,y)Y(y,z)=Y(x,z)} does not apply (cf.
\cite[Section 15]{DV-Dwork}). It may happen (cf. Remark
\ref{k_0-example}) that there exists a (smallest) $k_0\geq 0$ such
that condition \eqref{GOOD condition for Y(x,y)} holds for
$q^{k_0}$ instead of $q$, and for $Y_{A(q^{k_0},T)}(x,y)$ instead
of $Y_{A(q,T)}(x,y)$. There then exists  a Taylor solution $Y_c\in
M_n(\a_\Omega(c,R))$ of the iterated system
$\sigma_{q^{k_0}}(Y_c)=A(q^{k_0},T)Y_c$. In this case, for all
$c\in X(\Omega)$, we can recover a solution $Y^{\mathrm{big}}$ of
the system $\sigma_{q}(Y^{\mathrm{big}})=A(q,T)Y^{\mathrm{big}}$
itself in the algebra of analytic functions over the disjoint
union of disks $\bigcup_{i=0}^{k_0-1}\mathrm{D}^-(q^ic,R)$. Indeed
$\sigma_q$ acts on the algebra
$\prod_{i\in\mathbb{Z}/k_0\mathbb{Z}}M_n(\a_{K}(q^ic,R) )$ by
$\sigma_q((M_{q^ic}(T))_{i\in\mathbb{Z}/k_0\mathbb{Z}})=(M_{q^{i+1}c}(qT)
)_{i\in\mathbb{Z}/k_0\mathbb{Z}}$, and so one has
\begin{equation}\label{tilde Y-T-}
Y^{\mathrm{big}}(T)=(\;Y^{\mathrm{big}}_{q^ic}(T)\;)_i:=(\;
A(q^{i},q^{-i}T)\cdot Y_{c}(q^{-i}T)\;
)_{i\in\mathbb{Z}/k_0\mathbb{Z}}\;.
\end{equation}
In fact  $A(q^{i+1},q^{-i}T)=A(q,T)A(q^i,q^{-i}T)$. This and
related matters are very well explained in \cite{DV-Dwork}.

\subsubsection{} Notice that the relations of Lemma
\ref{Y(x,y)Y(y,z)=Y(x,z)} hold for $Y_A(x,y)$ as a function on
$\mathcal{U}_R$, and not for $Y^{\mathrm{big}}(T)$ (cf.
\eqref{tilde Y-T-}). In other words the expression
$Y^{\mathrm{big}}_A(x,y)$ has no meaning if $|x-y|\geq R$. In
particular the expression \eqref{eq-4-}, which is the main tool of
the Propagation Theorem \ref{main theorem second form}, holds only
if $|x-y|,|z-y|<R$.

\subsubsection{The case of a root of unity.} If $q\in\bs{\mu}(\Q)$ is a root of unity, then even when
a solution $Y\in GL_n(\a_\Omega(c,R))$  exists, \emph{the radius
is not defined} since we may have another solution with different
radius (cf. Example \ref{exemple radius not defined} below). For
this reason, the radius of convergence of the system
\eqref{Y(qT)=AY--} will be not defined if $q\in\bs{\mu}(\Q)$.

\begin{example}\label{exemple radius not defined}
Let $q=\xi$ be a $p-$th root of unity, with $\xi\neq 1$. The
solutions of the unit object at $t^p\in\Omega$ are the functions
$y\in\a_\Omega(t^p,R)$ such that $y(\xi T)=y(T)$.  Every function
in $T^p$ has this property. For example the family of functions
$\{\;y_\alpha:=\exp(\alpha
(T^p-t^p))\;\}_{\alpha\in\Omega}$,
is such that for different values of $\alpha$ one has different
radii.
\end{example}

\subsection{Taylor solutions of $(\sigma_q,\delta_q)-$modules}
In this subsection $q$ may be a root of unity. We preserve the
previous notations. We consider now a system (the notion of
solution of such a system have been defined in section
\ref{section constant solutions}):
\begin{eqnarray}\label{Y(qT)=AY and Y'=GY --}
&\left\{
\begin{array}{rcl}
\sigma_q(Y)&=&A(q,T)\cdot Y\;, \\
&&\\
\delta_q(Y)&=&G(q,T)\cdot Y\;,
\end{array}
\right.& \begin{array}{l}
A(q,T)\in GL_n(\H_K(X))\;,\\
\\
G(q,T)\in \phantom{I}M_n(\H_K(X))\;.
\end{array}
\end{eqnarray}
It can happen that a solution of $\sigma_q^{\M}$ is not a solution
of $\delta_q^{\M}$ as shown by the following example:
\begin{example}\label{example of second obstruction}
Suppose that $q\in\mathrm{D}^-(1,1)$ is not a root of unity. Let
$X:=\mathrm{D}^+(0,|p|^{\frac{1}{p-1}})$,
$A(q,T):=\exp((q-1)T)\in\H_K(X)^\times$, $G(q,T):=0$. Let $c=0$,
and $R<|p|^{\frac{1}{p-1}}$. Then every solution $y(T)\in\a_K(0,R)
$ of the operator $\sigma_q-A(q,T)$ is of the form
$y(T)=\lambda\cdot\exp(T)$, with $\lambda\in K$. If
$\delta_q(y)=0$, then $y=0$. Hence, the
$(\sigma_q,\delta_q)-$module defined by $A(q,T)$ and $G(q,T)$ has
no (non trivial) solutions in $\a_K(0,R)$.
\end{example}

To guarantee the existence of solutions we need a
\emph{compatibility condition} between $\sigma_q$ and $\delta_q$,
which should be written explicitly in terms of matrices of
$\sigma_q^n$ and $\delta_1^n$. This obstruction will not appear in
the sequel of the paper since this condition is automatically
satisfied by analytic $\sigma$-modules (cf. Lemma \ref{alpha
continuous}). This fact will follows from that a solution
$\alpha:\M\to\a_{\Omega}(c,R)$ is continuous (see the proof of
Lemma \ref{alpha continuous}). Observe that Lemma \ref{alpha
continuous} below is not a formal consequence of the previous
theory. Indeed, by Definition \ref{discrete sigma,delta algebra},
the general
$(\sigma,\delta)-$algebra $\C$ used in Definition %
\ref{discrete sigma,delta algebra} has the discrete topology,
hence the morphism $\alpha:\M\to\C$ defining the solution is not
continuous in general.

\begin{lemma}\label{alpha continuous}
Let $U\subseteq \Q(\H_K(X))$ be an open subset, and let $\M$ be an
analytic $(\sigma,\delta)-$module on $U$, representing the family
of equations $\{\sigma_q(Y)=A(q,T)\cdot Y\}_{q\in U}$, with
$A(q,T)\in GL_n(\H_K(X))$, for all $q\in U$. Let $Y_c(T)\in
GL_n(\a_\Omega(c,R))$, $|q-1||c|<R\leq\rho_{c,X}$, be a
\emph{simultaneous} solution of every equation of this family.
Then $Y_c(T)$ is also solution of the equation
\begin{equation}
\delta_q(Y)=G(q,T)\cdot Y\;,
\end{equation}
where $G(q,T):=q\frac{d}{dq}(A(q,T))$ (cf.
\eqref{G(q,T)=qd/dqA(q,T)}). Hence $Y_c(T)$ is solution of the
differential equation defined in section \ref{construction of
delta}:
\begin{equation}\label{Y_c(T) diff eq --}
\delta_1(Y_c(T))=G(1,T)\cdot Y_c(T)\;,
\end{equation}
where $G(1,T)=G(q,q^{-1}T)\cdot A(q,q^{-1}T)^{-1}\in M_n(\H_K(X))$
(cf. \eqref{G(qq',T) = G(q',qT) cdot A(q,T)}).
\end{lemma}
\begin{proof} In terms of modules, the columns of the matrix $Y_c(T)$
correspond to $\H_K(X)-$linear maps $\alpha:\M\to\a_\Omega(c,R)$,
verifying $\sigma_q\circ\alpha=\alpha\circ\sigma_q^{\M}$, for all
$q\in U$ (cf. Section \ref{morphisms as solutions}). We must show
that such an $\alpha$ also commutes with $\delta_q$. This follows
immediately by the continuity of $\alpha$. Indeed, the inclusion
$\H_K(X)\to\a_\Omega(c,R)$ is continuous, and hence every
$\H_K(X)-$linear map $\H_K(X)^n\to\a_\Omega(c,R)$ is continuous.
\end{proof}

\subsection{Twisted Taylor formula for
$(\sigma,\delta)-$modules, and rough estimate of radius}

Let $X$ be a $q$-invariant affinoid. Let
$D_q:=\sigma_q\circ\frac{d}{dT}= \lim_{q'\to
q}\frac{\sigma_{q'}-\sigma_q}{T(q'-q)}=
\frac{1}{qT}\cdot\delta_q$. For all $q\in\Q(X)$ and all
$f(T)\in\H_K(X)$, one has
\begin{eqnarray}
\D_q(f\cdot
g)&=&\sigma_q(f)\cdot\D_q(g)+\D_q(f)\cdot\sigma_q(g)\;,\\
(d/dT\circ\sigma_q)&=&q\cdot(\sigma_q\circ d/dT)\;,\\
\D_q^n&=&q^{n(n-1)/2}\cdot\sigma_q^n\circ(d/dT)^n\;,\\
\|\D_q^n(f(T))\|_{X}&\leq&
\frac{|n!|}{r_X^n}\cdot\|f(T)\|_{X}\qquad\textrm{(cf. Lemma
\ref{|f'| leq |f| / r_A})}\;.
\end{eqnarray}
Hence, for all $c\in K$, $\D_q^n(T-c)^i=\frac{i!}{(i-n)!}\cdot
q^{n(n-1)/2}\cdot(q^{n}T-c)^{i-n}$ if $n\leq i$, and
$\D_q^n(T-c)^i=0$ if $n>i$. This shows that if $f(T):=\sum_{i\geq
0}a_i\cdot\frac{(T-c)^i}{(i!)\cdot q^{i(i-1)/2}}\in
\a_{\Omega}(c,R)$ is a formal series, with
$|q-1||c|<R\leq\rho_{c,X}$, then $a_n=\D_q^n(f)(c/q^n)$, and the
usual Taylor formula can be written as
\begin{equation}\label{taylor expansion of a formal series}
f(T)=\sum_{n\geq 0}\D_q^n(f)(c/q^n)\cdot\frac{(T-c)^n}{(n!)\cdot
q^{n(n-1)/2}}\quad.
\end{equation}

The following proposition gives the analogue of the
classical rough estimate for differential and $q-$difference
equations (cf. \cite[4.1.2]{Ch}, \cite[4.3]{DV-Dwork}).

\begin{proposition}\label{Taylor (sigma,delta)-solution has a big radius}
Let $c\in X(\Omega)$. Assume that the system \eqref{Y(qT)=AY and
Y'=GY --} has a Taylor solution $Y_c\in M_n(\a_\Omega(c,R_c))$,
with $|q-1||c|<R_c\leq\rho_{c,X}$. For all $q$-invariant
sub-affinoid $X'\subseteq X$, containing
$\mathrm{D}^+(c,|q-1||c|)$, one has
\begin{equation}\label{hui}
R_c \geq
\frac{|p|^{\frac{1}{p-1}}}{\max(\;r_{X'}^{-1}\|A(q,T)\|_{X'}\;,
\;\|G(q,T)/qT\|_{X'}\; )}\;.
\end{equation}
In particular if $X'$ is a disk $\mathrm{D}^+(c,\rho)$, with
$|q-1||c|\leq \rho\leq\rho_{c,X}$, then
\begin{equation}
R_c \geq
\frac{|p|^{\frac{1}{p-1}}\cdot\rho}{\max(\;|A(q,T)|_{(c,\rho)}\;,
\;\frac{|G(q,T)|_{(c,\rho)}}{\max(1,|c|/\rho)}\; )}\;.
\end{equation}
\end{proposition}
\begin{proof} The matrix $Y_c(T)$ satisfies
$\sigma_q^n(Y_c(T))=A_{[n]}(q,T)\cdot Y_c(T)$, and
$D_q^n(Y_c(T))=F_{[n]}(q,T)\cdot Y_c(T)$, where
$F_{[0]}=\mathrm{Id}=A_{[0]}$, $A_{[1]}:=A(q,T)$,
$F_{[1]}:=\frac{1}{qT}G(q,T)$, and
\begin{eqnarray}
A_{[n]}&:=&\sigma_q^{n-1}(A_{[1]})\cdots \sigma_q(A_{[1]})\cdot
A_{[1]}\;,\\
F_{[n+1]}&:=&\sigma_q(F_{[n]})\cdot F_{[1]}+D_q(F_{[n]})\cdot
A_{[1]}\;.\label{G_n+1=G_1 sigma_q(G_n)+ A D_q(G_n)}
\end{eqnarray}
Hence one has
\begin{equation}\label{Taylor formal solution of delta_q}
Y_c(T):=\sum_{i\geq 0}F_{[n]}(c/q^n)\frac{(T-c)^n}{(n!)\cdot
q^{n(n-1)/2}}\;\;,
\end{equation}
which is a hybrid between the usual Taylor formula and the
Taylor formula for $q-$difference equations. Inequalities
\ref{hui} then follow from the inequality
\begin{equation}
|F_{[n]}(c/q^n)|_{\Omega}\leq\|F_{[n]}\|_{X'}\leq
\max\left(\|F_{[1]}\|_{X'}\;,\;\frac{1}{r_{X'}}\cdot\|A_{[1]}\|_{X'}\right)^n
\end{equation}
If $X'=\mathrm{D}^+(c,\rho)$, then the last term is equal to
$\frac{1}{\rho^n}\cdot
\max\left(\frac{|G(q,T)|_{(c,\rho)}}{\max(1,|c|/\rho)}\;,\;|A(q,T)|_{(c,\rho)}\right)^n.$
Indeed $r_{\mathrm{D}^+(c,\rho)}=\rho$,
$F_{[1]}=\frac{1}{qT}G(q,T)$, and
$|T|_{(c,\rho)}=|(T-c)+c|_{(c,\rho)}=\max(\rho,|c|)$, hence
$|F_{[1]}|_{(c,\rho)}=\frac{1}{|q|\max(|c|,\rho)}\cdot|G(q,T)|_{(c,\rho)}$,
and $|q|=1$. \end{proof}

\section{Generic radius of convergence and solvability}

\begin{definition}[(Generic radius of convergence)]\label{generic radius of
convergence} Let $q\in\Q(X)$ (resp. $q\in \Q(X)-\bs{\mu}(\Q)$),
let $c\in X(K^{\mathrm{alg}})$, and let $\mathrm{D}^+(c,\rho)$,
$|q-1||c|<\rho\leq\rho_{c,X}$,  be a $q$-invariant disk. Let $\M$
be the $(\sigma_q,\delta_q)-$module (resp. $\sigma_q-$module)
defined by the system \eqref{Y(qT)=AY and Y'=GY --} (resp.
\eqref{Y(qT)=AY--}). Let
$R_{t_{c,\rho}}:=Ray(Y(x,y),t_{c,\rho})=Ray(Y(x,y),|.|_{(c,\rho)})$
be the radius of convergence\footnote{In the case of the
$q$-difference equation \eqref{Y(qT)=AY--}, the radius
$R_{t_{c,\rho}}$ is given by definition \eqref{R_c--}. In the case
of the system \eqref{Y(qT)=AY and Y'=GY --} the radius
$R_{t_{c,\rho}}$ is given indifferently by definition \eqref{R_c
in the differential case} or by definition \eqref{R_c--}, indeed
under our assumptions these two definitions are equal since
$Y_{A(q,T)}(x,y)=Y_{G(1,T)}(x,y)$. However observe that the
definition \eqref{R_c--} exists only if $q \in \Q-\bs{\mu}(\Q)$,
while definition \eqref{R_c in the differential case} preserves
its meaning on the root of unity.} of $Y_{A(q,T)}(T,t_{c,\rho})$.
Assume that
\begin{equation}\label{frefrefre}
|q-1||t_{c,\rho}|<R_{t_{c,\rho}} 
\;.\footnote{Observe that
$\rho_{c,X}=\rho_{t_{c,\rho},X}$, indeed
$\mathrm{D}^+(c,r)=\mathrm{D}^+(t_{c,\rho},r)$, for all $r\geq
\rho$.}
\end{equation}
We define the \emph{$(c,\rho)-$generic radius of convergence of
$\M$} to be the real number
\begin{equation}\label{Ray}
Ray(\M,|.|_{c,\rho})\;:=\;\min\left(\;R_{t_{c,\rho}}\;,\;\rho_{c,X}\;\right)
\; > \; |q-1||c|\;.
\end{equation}
\end{definition}

\subsubsection{} The assumption \eqref{frefrefre} ensures that the
disk of convergence of $Y(x,y)$ at $y=t_{c,\rho}$ is
$q$-invariant. While the bound
$Ray(\M,|.|_{c,\rho})\leq\rho_{c,X}$ ensures that $Y(x,y)$ is
invertible in the disk $\mathrm{D}^-(t_{c,\rho},R)$, for all
$0<R\leq Ray(\M,|.|_{c,\rho}) $ (cf. Lemma \ref{det neq 0}). We
recall that $|t_{c,\rho}|=\min(|c|,\rho)$, and that
$\|\cdot\|_{\D^+(c,\rho)}=\max_{y_0\in
\D^+_{K^{\mathrm{alg}}}(c,\rho)}|\cdot|_{y_0}$. Hence, by the
transfer principle (cf. Section \ref{transfert}), one has:
\begin{equation}
R_{t_{c,\rho}} := Ray(Y(x,y),t_{c,\rho}) =
Ray(Y(x,y),\mathrm{D}^+(c,\rho)) =
\min_{y_0\in\mathrm{D}^+_{K^\mathrm{alg}}(c,\rho)}
Ray(Y(x,y),y_0)\;.
\end{equation}
The number $Ray(\M,|.|_{(c,\rho)})$ is invariant under change of
basis in $\M$, while the number
$R_{t_{c,\rho}}=Ray(Y(x,y),|.|_{(c,\rho)})$ depends on the choice
of basis. Observe that $Ray(\M,|.|_{(c,\rho)})$ depends on the
affinoid $X$, and on the semi-norm $|.|_{(c,\rho)}$ defined by
$t_{c,\rho}$, but not on the particular choice of $t_{c,\rho}$
(cf. Section \ref{remremju}).

\begin{definition}[(Solvability)]
Let $\M$ be a $\sigma_q-$module (resp. a
$(\sigma_q,\delta_q)-$module) on $\H_K(X)$. We will say that $\M$
is solvable at $t_{c,\rho}$ if
\begin{equation}
Ray(\M,|.|_{(c,\rho)})=\rho_{c,X}\;.
\end{equation}
\end{definition}

\subsubsection{Continuity and $\log$-concavity of the Radius.}
Notice that every point $|.|_*$ in the Berkovich space associated
to $X$ is of the form $|.|_{c,\rho}$, for a suitable $\rho\geq 0$,
and for a point $c$ in $X(L)$, where $(L,|.|)/(K,|.|)$ is a
sufficiently large extension of complete valued fields. One may
verify that $|.|_*\mapsto Ray(\M,|.|_{c,\rho})$ is a well defined
function on the Berkovich space (i.e. the Radius does not depends
on the chosen $c$, but only on $|.|_*$). In a recent pre-print
(cf. \cite{DV-Balda}) it have been proved that the function $|.|_*
\mapsto Ray(\M,|.|_*)$ is continuous on the Berkovich Space. We
refer to \cite{DV-Balda} for a very inspiring treatment to this
subject.

We notice that this generalizes a previous statement (cf.
\cite{Ch-Dw}) proving, for all $c\in X(L)$, the continuity of the
function $\rho\mapsto Ray(\M,|.|_{c,\rho})$.

Let now $(L,|.|)/(K,|.|)$ be any extension of complete valued
fields. Let $c\in X(L)$. The function $\rho\mapsto
Ray(\M,|.|_{c,\rho})$ defined on $[0,\rho_{c,X}]$ is
$\log$-concave (cf. Def. \ref{log graphic ..}), and it can be
proved that it is piecewise $\log$-affine. This follows
essentially by the definition of the Radius (cf. \eqref{R_c--}),
and by Lemma \ref{log-convex property}.

\subsection{Solvability over an annulus and over the Robba ring}

Let $\B:=\a_K(I)$, with $I=]r_1,r_2[$, and let $\M$ be a
$\sigma_q-$module (resp. a $(\sigma_q,\delta_q)-$module) on
$\a_K(I)$.
For all $c\in K$, $|c|\in I$, one has
$t_{c,|c|}=t_{0,|c|}$. For all affinoid
$X\subseteq\mathcal{C}(I)$ containing the disk
$\mathrm{D}^-(c,|c|)$ one has $\rho_{c,X}=|c|$. Then the norm
$|.|_{c,|c|}:\a_K(I)\longrightarrow\mathbb{R}_{\geq}$ and the
generic radius $Ray(\M,|.|_{(c,|c|)})$, do not depend on the
choice of $c$  or the affinoid $X$, but only on $|c|$. Hence, for all $\rho\in
I$, we chose an arbitrary $c\in \Omega$, with $|c|=\rho\in I$, and
we set
\begin{equation}\label{ZERTY}
\index{t_rho@$t_{\rho}:=t_{0,\rho}$}
\index{Ray(M,rho)@$Ray(\M,\rho):=Ray_0(\M,\rho)$}
t_\rho:=t_{c,\rho}\;,\quad \mathrm{and}\qquad
Ray(\M,\rho):=Ray(\M,|.|_{(c,\rho)})\;.
\end{equation}

To define the radius we need the assumption
$|q-1||t_{\rho}|<\rho_{t_{\rho},X}=\rho$ (cf. Definition
\ref{generic radius of convergence}). Since $|t_{\rho}|=\rho$,
this assumption is equivalent to
\begin{equation}
|q-1|<1\;.
\end{equation}

\begin{definition}[(solvability at $\rho$)]
\index{sigma_q-Mod(A_K(I))^sol(rho)@$\sigma_q-\mathrm{Mod}(\a_K(I))^{\mathrm{sol(\rho)}}$}
Let $q\in\Q_1-\bs{\mu}(\Q_1)$ (cf. Definition \eqref{Q and Q_1}).
Let $\M$ be a $\sigma_q-$module on $\a_K(I)$. We will say that
$\M$ is solvable at $\rho\in I$ if
\begin{equation}
Ray(\M,\rho)=\rho\;.
\end{equation}
\end{definition}

\subsubsection{Solvability over $\R_K$ or $\Hd_K$.}

Let $q\in\Q_1-\bs{\mu}(\Q_1)$. Let $\M$ be a $\sigma_q-$module 
over $\R_K$. By definition $\M$ comes, by scalar extension, from a
module $\M_{\varepsilon_1}$ defined on an annulus
$\mathcal{C}(]1-\varepsilon_1,1[)$. If $\varepsilon_2>0$, and if
$\M_{\varepsilon_2}$ is another module on
$\mathcal{C}(]1-\varepsilon_2,1[)$ satisfying
$\M_{\varepsilon_2}\otimes_{\a_K(]1-\varepsilon_2,1[)}\R_K\simto\M$,
then there exists a
$\varepsilon_3\leq\min(\varepsilon_1,\varepsilon_2)$ such that
\begin{equation}
\M_{\varepsilon_1}\otimes\a_K(]1-\varepsilon_3,1[)\xrightarrow[]{\;\;\sim\;\;}
\M_{\varepsilon_2}\otimes\a_K(]1-\varepsilon_3,1[)\;.
\end{equation}
Hence the limit $\lim_{\rho\to 1^-}Ray(\M_{\varepsilon},\rho)$ is
independent of the choice of the module $\M_{\varepsilon}$.

\begin{definition}\label{sigma-Mod(Hd)^[r]}
Let $q\in\Q_1-\bs{\mu}(\Q_1)$, and let $|q-1|< r\leq 1$. We define
\begin{equation}\index{sigma_q-Mod(H^dag_K)^[r]@$\sigma_q-\Mod(\Hd_K)^{[r]}$}
\sigma_q-\Mod(\Hd_K)^{[r]}\;,
\end{equation}
as the full sub category of $\sigma_q-\Mod(\Hd_K)$ whose objects
satisfy
\begin{eqnarray}
Ray(\M,1)\geq r \;,\qquad(r > |q-1|)\;,
\end{eqnarray}
as illustrated below in the log-graphic of the function
$\log(\rho)\mapsto \log(Ray(\M,\rho)/\rho)$ (cf. Def. \ref{log
graphic ..}):
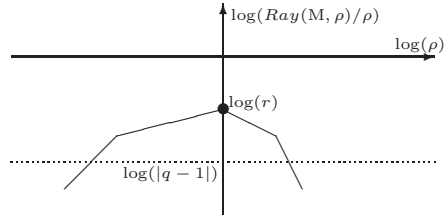
\begin{center}
\begin{picture}(160,80)
\put(80,0){\vector(0,1){80}}%
\put(0,60){\vector(1,0){160}}%
\put(20,10){\line(1,1){20}}%
\put(40,30){\line(4,1){40}}%
\put(80,40){\line(2,-1){20}}%
\put(100,30){\line(1,-2){10}}
\qbezier[80](0,20)(80,20)(160,20)%
\put(77.5,37.5){$\bullet$} %
\put(82,41){\begin{tiny}$\log(r)$\end{tiny}}
\put(42,14){\begin{tiny}$\log(|q-1|)$\end{tiny}}
\put(145,63){\begin{tiny}$\log(\rho)$\end{tiny}}
\put(83,74){\begin{tiny}$\log(Ray(\M,\rho)/\rho)$\end{tiny}}
\end{picture}
\end{center}
Objects in $\sigma_q-\Mod(\Hd_K)^{[1]}$ will be called
\emph{solvable}.
\end{definition}

\begin{definition}\label{sigma-Mod(R)^[r]}
Let $q\in\Q_1-\bs{\mu}(\Q_1)$, and let $|q-1|\leq r\leq 1$. We
define
\begin{equation}\index{sigma_q-Mod(R_K)^[r]@$\sigma_q-\Mod(\R_K)^{[r]}$}
\sigma_q-\Mod(\R_K)^{[r]}\;,
\end{equation}
as the full sub category of $\sigma_q-\Mod(\R_K)$ formed by
objects $\M$ satisfying $\lim_{\rho\to 1^{-}}Ray(\M,\rho)\geq r$,
and there exists $\varepsilon_q>0$ such that $Ray(\M,\rho) >
|q-1|$, for all $\rho\in]1-\varepsilon_q,1[$. 
\if{\begin{enumerate} \item  if $r>|q-1|$,
then $\lim_{\rho\to 1^{-}}Ray(\M,\rho)\geq r$,

\item if $r=|q-1|$, then there exists $\varepsilon_q>0$ such that
$Ray(\M,\rho) > r$, for all $\rho\in]1-\varepsilon_q,1[$.
\end{enumerate}}\fi
\if{This condition can be also expressed by asking
\begin{eqnarray}
\lim_{\rho\to 1^-}Ray(\M,\rho)\geq r \;.
\end{eqnarray}
and moreover, in the particular case in which $\lim_{\rho\to
1^-}Ray(\M,\rho)=r=|q-1|<1$, we ask that there exist
$\varepsilon_q>0$, such that $Ray(\M_{\varepsilon_q},\rho)>
|q-1|\rho$, for all $1-\varepsilon_q<\rho<1$,}\fi
There are two possible cases $r>|q-1|$, and $r=|q-1|$, as
illustrated in the following pictures:

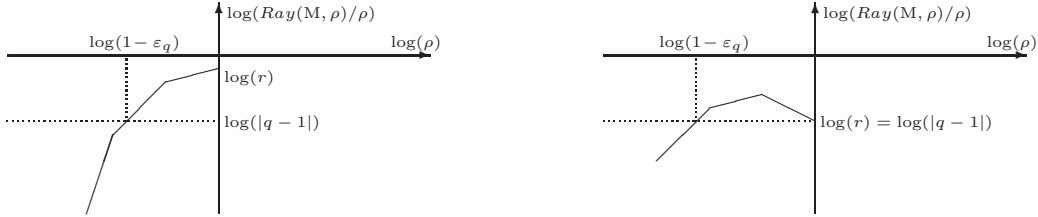
\begin{picture}(220,90)
\put(100,5){\vector(0,1){80}}%
\put(20,65){\vector(1,0){160}}%
\put(50,5){\line(1,3){10}}%
\put(60,35){\line(1,1){20}}%
\put(80,55){\line(4,1){20}}
\qbezier[40](20,40)(60,40)(100,40)%
\qbezier[15](65,40)(65,52.5)(65,65)%
\put(51,68){\begin{tiny}$\log(1\!-\varepsilon_q)$\end{tiny}}
\put(102,38){\begin{tiny}$\log(|q-1|)$\end{tiny}}
\put(102,55){\begin{tiny}$\log(r)$\end{tiny}}
\put(165,68){\begin{tiny}$\log(\rho)$\end{tiny}}
\put(103,79){\begin{tiny}$\log(Ray(\M,\rho)/\rho)$\end{tiny}}
\end{picture}
\qquad
\begin{picture}(160,80)
\put(80,5){\vector(0,1){80}}%
\put(0,65){\vector(1,0){160}}%
\put(20,25){\line(1,1){20}}%
\put(40,45){\line(4,1){20}}%
\put(60,50){\line(2,-1){20}}
\qbezier[40](0,40)(40,40)(80,40)%
\qbezier[15](35,40)(35,52.5)(35,65)%
\put(21,68){\begin{tiny}$\log(1\!-\varepsilon_q)$\end{tiny}}
\put(82,38){\begin{tiny}$\log(r)=\log(|q-1|)$\end{tiny}}
\put(145,68){\begin{tiny}$\log(\rho)$\end{tiny}}
\put(83,79){\begin{tiny}$\log(Ray(\M,\rho)/\rho)$\end{tiny}}
\end{picture}\\
Objects in $\sigma_q-\Mod(\R_K)^{[1]}$ will be called
\emph{solvable}.
\end{definition}

\begin{remark}
Notice that in definition \ref{sigma-Mod(Hd)^[r]} the existence of
$\varepsilon_q>0$ such that $Ray(\M,\rho)
> |q-1|$, for all $\rho\in]1-\varepsilon_q,1+\varepsilon_q[$ is
automatically verified since one assumes $r>|q-1|$.
\end{remark}

\subsubsection{Analogous definitions for
$(\sigma_q,\delta_q)$-modules.}
\index{sigma_q,delta_q-Mod(R_K)^[r]@$(\sigma_q,\delta_q)-\Mod(\R_K)^{[r]}$,
$(\sigma_q,\delta_q)-\Mod(\Hd_K)^{[r]}$,
$(\sigma_q,\delta_q)-\mathrm{Mod}(\a_K(I))^{\mathrm{sol(\rho)}}$}
In the case of $(\sigma_q,\delta_q)$-modules, the generic radius
of convergence is defined even if $q$ is a root of unity. We give
then analogous definitions of
$(\sigma_q,\delta_q)-\Mod(\B)^{[r]}$, for $\B:=\R_K$ or
$\B:=\Hd_K$, without any restrictions on
$q$. 

\subsection{Generic radius for discrete and analytic
objects over $\R_K$ and $\Hd_K$}\label{generic radius for discrete
or analytic objects} In this section $\B=\R_K$ or $\B=\Hd_K$.
\begin{definition}
For all $\varepsilon>0$ let
\begin{equation}\label{I_epsilon}
I_\varepsilon:= \left\{
\begin{array}{lcl}
]\;1-\varepsilon\;,\;1\;[\;,&\textrm{ if }&\B=\R_K\\
&&\\
]\;1-\varepsilon\;,\;1+\varepsilon\;[\;,&\textrm{ if }&\B=\Hd_K\;\quad.\\
\end{array}
\right.
\end{equation}
\end{definition}

\begin{definition} For all subset
$S\subseteq\mathrm{D}^-(1,1)=\Q_1$, for all $0<\tau<1$, we set
\begin{equation}
S_\tau:=S\cap\mathrm{D}^-(1,\tau)\;.
\end{equation}
\end{definition}

\begin{definition}\label{sigma-Mod(R)^[r]_S}
Let $0 < r \leq 1$. Let $S\subseteq\mathrm{D}^-(1,1)$,
$S^\circ\neq \emptyset$. We denote by
\begin{equation}\index{sigma-Mod(R_K)^[r]_S@$\sigma-\Mod(\R_K)^{[r]}_S$,
$\sigma-\Mod(\Hd_K)^{[r]}_S$}
 \sigma-\Mod(\B)^{[r]}_S\;,
\end{equation}
the full subcategory of $\sigma-\Mod(\B)_S$ whose objects $\M$
have the following properties:
\begin{enumerate}
\item The restriction of $\M$ to every $q\in S$ 
belongs to $\sigma_q-\Mod(\B)^{[r]}\;;$ %
\item For all $\tau$ such that $0<\tau<r$, there exists
$\varepsilon_{\tau}>0$ such that the restriction
$\mathrm{Res}^{\ph{S}}_{\ph{S_\tau}}(\M)$ comes, by scalar
extension, from an object
\begin{equation}
\M_{\varepsilon_{\tau}}\in\;\;
\sigma-\Mod(\a_K(I_{\varepsilon_{\tau}}))_{S_\tau}^{\mathrm{disc}}
\end{equation}
such that, for all $\rho\in I_{\varepsilon_\tau}$, and for all
$q,q'\in S_\tau$ one has (cf. \eqref{Y_A(T,c)})
\begin{equation}\label{Yq=Y_q' for all gggg}
Y_{A(q,T)}(T,t_{\rho})=Y_{A(q',T)}(T,t_\rho)\;.
\end{equation}
\end{enumerate}
Objects in $\sigma_q-\Mod(\B)_S^{[1]}$ will be called
\emph{solvable}.
\end{definition}

\begin{example}
This example justifies the condition i) given in the preceding
definition. Let $r:=\omega:=|p|^{\frac{1}{p-1}}$, and let
$S=\mathrm{D}^-(1,\omega)$. Let $\M$ be the discrete
$\sigma$-module over the Robba ring defined by the family of
equations $\{\;\sigma_q-A(q,T)\;\}_{q\in S}$, where
$A(q,T):=\exp((q^{-1}-1)T^{-1})$. Then
$Y(x,y):=\exp(x^{-1}-y^{-1})$ is the simultaneous solution of
every equation of this family. Observe that $A(q,T)\in\R_K$ if and
only if $|q^{-1}-1|<\omega$, but if $|q-1|$ tends to $\omega^-$,
then \emph{the matrices $A(q,T)$ do not all belong to the same
annulus}. Indeed $A(q,T)\in\a_K(I_{\varepsilon})$ if and only if
$|q^{-1}-1|<\omega(1-\varepsilon)$.
\end{example}

\begin{remark}
Condition i) implicitly implies that $S\subseteq\mathrm{D}^-(1,r)$
if $\B=\Hd_K$ (cf. Def. \ref{sigma-Mod(Hd)^[r]}), and
$S\subseteq\mathrm{D}^+(1,r)$ if $\B=\R_K$ (cf Def.
\ref{sigma-Mod(R)^[r]}).
\end{remark}

\subsubsection{Analogous definitions for
$(\sigma_q,\delta_q)$-modules.}\label{truy}
\index{sigma,delta-Mod(R_K)^[r]_S@$(\sigma,\delta)-\Mod(\R_K)^{[r]}_S$,
$(\sigma,\delta)-\Mod(\Hd_K)^{[r]}_S$} One defines analogously
$(\sigma,\delta)-\Mod(\B)^{[r]}_S$, \emph{but without restrictions
on $S\subseteq \mathrm{D}^-(1,r)$}, as the subcategory of
$(\sigma,\delta)-\Mod(\B)_S$, whose objects verify conditions i)
and ii), in which equation \eqref{Yq=Y_q' for all gggg} is
replaced by (cf. Definitions \eqref{Y_G-+-} and \eqref{Y_A(T,c)})
\begin{equation}
Y_{G(1,T)}(T,t_{\rho})=Y_{A(q,T)}(T,t_\rho)\;,
\end{equation}
 for all $\rho\in
I_{\varepsilon_{\tau}}$, and all $q\in S_\tau$.

\section{The Propagation Theorem}

\label{section - main theorem}

\subsection{Taylor admissible modules}

\begin{definition}[(Taylor admissible discrete modules on $S$)]
\label{taylor taylor admissible} \label{def of Taylor
adm}\label{def of Taylor adm}   Let
$X:=\mathrm{D}^+(c_0,R_0)-\cup_{i=1}^n\mathrm{D}^-(c_i,R_i)$ be an
affinoid, and let
$S\subseteq\Q_1(X)$, be a subset with $S^\circ\neq\emptyset$ (cf. \eqref{S'}). %
Let $(\M,\sigma^{\M})$ be a discrete $\sigma$-module defined by
the family of equations
\begin{equation}\label{gugg}
\{\sigma_q-A(q,T)\}_{q\in S}\;,\;\; A(q,T)\in
GL_n(\H_K(X))\;,\;\forall\;q\in S\;.
\end{equation}
We will say that $(\M,\sigma^{\M})$ is \emph{Taylor admissible on
$X$, with generic radius greater than $r$}, if :
\begin{itemize}
\item[$(1)$] One has
$S\subseteq\mathrm{D}^-(1,r/\max(|c_0|,R_0))$; %
\item[$(2)$] there exists a matrix $Y(x,y)$, convergent in
$\mathcal{U}_R$ (cf. \eqref{U_R}), with $R\geq r$ satisfying, for
all $q\in S$, the condition \eqref{FRREISTUGHY}, that is
\begin{equation}\label{condition for the admissibility}
 r \; \leq \; R \; \leq \; r_X\;;
\end{equation}
\item[$(3)$]  $Y(x,y)$ is simultaneous solution of every equation
of the family \eqref{gugg}.
\end{itemize}
The full subcategory of $\sigma-\Mod(\H_K(X))_S^{\mathrm{disc}}$
whose objects are Taylor admissible, with generic radius greater
than $r$, will be denoted by
\begin{equation}\index{sigma,delta-Mod(H_K(X))^r_S@$(\sigma,\delta)-\Mod(\H_K(X))^{[r]}_S$,
$\sigma-\Mod(\H_K(X))^{[r]}_S$} \sigma-\Mod(\H_K(X))_S^{[r]}\;.
\end{equation}
Moreover we
set
\begin{equation}\label{Remark - ffff}
\index{sigma,delta-Mod(H_K(X))^adm_S@$(\sigma,\delta)-\Mod(\H_K(X))^{\mathrm{adm}}_S$,
$\sigma-\Mod(\H_K(X))^{\mathrm{adm}}_S$}
\sigma-\Mod(\H_K(X))_S^{\mathrm{adm}}\;:=\;\bigcup_{r}\sigma-\Mod(\H_K(X))_S^{[r]}\;\;.
\end{equation}
where $r\leq r_X$ runs in the set of real numbers such that
$S\subseteq\mathrm{D}^-(1,r/\max(|c_0|,R_0))$. We define
analogously the categories $(\sigma,\delta)-\Mod(\H_K(X))^{[r]}_S$
and $(\sigma,\delta)-\Mod(\H_K(X))^{\mathrm{adm}}_S$ of
\emph{admissible $(\sigma,\delta)-$modules on $S$}. Namely the
condition $S^\circ\neq\emptyset$ is suppressed, and if
$(\M,\sigma^{\M},\delta_1^{\M})$ is a discrete
$(\sigma,\delta)$-modules on $S$ defined by a system of equations
(cf. \eqref{delta(Y)=YG}), then the Taylor solution
$Y_{G(1,T)}(x,y)$ (cf. \eqref{Y_G-+-}) of the differential
equation defined by $\delta_1^{\M}$ satisfies \eqref{condition for
the admissibility}, and moreover is simultaneously solution of
every equation defined by $\sigma_q^{\M}$, for all $q\in S$.
\end{definition}

\subsubsection{Taylor Admissibility over $\Hd_K(X)$.} We define
\begin{equation}
\sigma-\Mod(\Hd_K(X))^{[r]}_S \;,\quad (\textrm{resp.}\;
(\sigma,\delta)-\Mod(\Hd_K(X))^{[r]}_S )
\end{equation}
as the full subcategory of $\sigma-\Mod(\Hd_K(X))_S$ (resp.
$(\sigma,\delta)-\Mod(\Hd_K(X))_S$) formed by objects whose
restriction belongs to $\sigma-\Mod(\H_K(X))_S^{[r]}$ (resp.
$(\sigma,\delta)-\Mod(\H_K(X))^{[r]}_S$).

\begin{remark}
If $X=\{|T|=1\}$, $\Hd_K(X)=\Hd_K$ (cf. \eqref{Hd_K}), this
definition is equivalent to Def. \ref{sigma-Mod(R)^[r]_S}.
\end{remark}

\subsubsection{Taylor admissibility over $\R_K$.}
We preserve the notations of section \ref{generic radius for
discrete or analytic objects}.
\begin{definition}
We will say that an object is Taylor admissible over an annulus
$\mathcal{C}(I)$ if its restriction to every sub-annulus
$\mathcal{C}(J)$, with $J$ compact, $J\subseteq I$, is Taylor
admissible (cf. Definition \ref{taylor taylor admissible}).
\end{definition}

One defines Taylor admissibility over $\R_K$ by reducing to the
case of modules over a single annulus
$\mathcal{C}(I_\varepsilon)$, for some $\varepsilon>0$
sufficiently close to $0$. One finds in this way exactly the
Definition \ref{sigma-Mod(R)^[r]_S}:

\begin{definition}\label{Taylor adm over R_K and Hd_K .....}
Let $S\subseteq\mathrm{D}^-(1,1)$, with $S^{\circ}\neq \emptyset$.
Let $\tau_S:=\sup_{q\in S}|q-1|$. We set
\begin{equation}
\sigma-\Mod(\R_K)_S^{\mathrm{adm}}\;:=\;
\sigma-\Mod(\R_K)_S^{[\tau_S]}\;.
\end{equation}
We give the same definition for $(\sigma,\delta)$-modules, without
assuming that `` $S^{\circ}\neq \emptyset$ '' :
$(\sigma,\delta)-\Mod(\R_K)_S^{\mathrm{adm}}:=(\sigma,\delta)-\Mod(\R_K)_S^{[\tau_S]}$.
\end{definition}

\subsection{Propagation Theorem}

\begin{remark}\label{Remark explaining the Propagation}
We preserve notations of Definition \ref{taylor taylor
admissible}. If $\M$ is Taylor admissible on $X$, then, in
particular, $\M$ is trivialized by $\a_K(c,R)$, for all $c\in
X(K)$. Hence we can apply $\C$-Deformation and $\C$-Confluence to
$\M$, with $\C=\a_K(c,R)$ (cf. section \ref{constant deformation
section}). It will follows from the proof of Theorem \ref{main
theorem second form}, that this confluence does not depend on the
chosen point $c\in X(K)$.
\end{remark}

\begin{theorem}[(Propagation Theorem first form)]\label{main theorem}
Let $X$ be an \emph{affinoid}. Then, if $q\in
\Q_1(X)-\bs{\mu}(\Q_1(X))$, the natural restriction functor
\begin{equation}
\bigcup_{U}\mathrm{Res}^{U}_q\;:\;
\bigcup_{U}\sigma-\Mod(\H_K(X))_U^{\mathrm{adm}}
\;\xrightarrow[]{\;\;\;\;}\;\sigma_q-\Mod(\H_K(X))^{\mathrm{adm}}
\end{equation}
is an equivalence, where $U$ runs over the set of all open
neighborhood of $q$. The analogous fact is true for
$(\sigma,\delta)$-modules without supposing $q\notin\bs{\mu}(\Q)$.
\end{theorem}
\begin{proof} By Lemma \ref{The forget from S to S' is fully faithful},
$\cup_U\mathrm{Res}^U_{\{q\}}$ is fully faithful. Indeed for all modules
$\M,\N$  over $U$, by admissibility, there exists a number
$R$, with $|q-1|\max(|c_0|,R_0)<R\leq r_X$, such that, for all
$c\in X(K)$, the algebra $\mathrm{C}:=\a_K(c,R)$ trivializes both
$\M$ and $\N$. The essential surjectivity of
$\cup_{U}\mathrm{Res}^U_{\{q\}}$ will follow from Theorem
\ref{main theorem second form} below. \end{proof}

\begin{theorem}[(Propagation Theorem second form)]\label{main theorem second form}
Let $X=\mathrm{D}^+(c_0,R_0)-\cup_{i=1}^n\mathrm{D}^-(c_i,R_i)$.
Let $q\in\Q_1(X)-\bs{\mu}(\Q_1(X))$. Let
\begin{equation}\label{system---}
Y(q\cdot T)=A(T)\cdot Y(T)\;,\quad A(T)\in GL_n(\H_K(X))\;
\end{equation}
be a \emph{Taylor admissible} $q-$difference equation (cf. Def.
\ref{taylor taylor admissible}).
Then there exists a matrix $A(Q,T)$ uniquely determined by the
following properties:
\begin{enumerate}
\item $A(Q,T)$ is analytic and invertible in the domain
\begin{equation}
\mathrm{D}^-\Bigl(\;1\;,\;\frac{R}{\max(|c_0|,R_0)}\;\Bigr)\times
X\;\;\subset\; \mathbb{A}_K^2\;,
\end{equation} %
\item The matrix $A(Q,T)$ specialized at $(q,T)$ is equal to $A(T)$\;, %
\item For all $q'\in \mathrm{D}^-(1,R/\max(|c_0|,R_0))$, the
Taylor solution matrix $Y_A(x,y)$ of the equation
\eqref{system---} (cf. \eqref{Y_A(T,c)}) simultaneously satisfies
\begin{equation}\label{frequency}
Y_A(q'\cdot T,y)\;=\;A(q',T)\cdot Y_A(T,y)\;\;.
\end{equation} %
\end{enumerate}
Moreover the matrix $A(Q,T)$ is independent of the choice of
solution $Y_A(x,y)$.
\end{theorem}
\begin{proof}
By equation \eqref{frequency}, the matrix $A(Q,T)$ must be equal
to
\begin{equation}
A(Q,T)=Y_A(Q\cdot T,y)\cdot Y_A(T,y)^{-1}=Y_A(Q\cdot T,y)\cdot
Y_A(y,T)=Y_A(Q\cdot T,T)\;.\label{independence from c----}
\end{equation}
This makes sense since $Y_{A(q,T)}(x,y)$ is invertible in its
domain of convergence (cf. Lemma \ref{Y(x,y)Y(y,z)=Y(x,z)}). Hence
$A(Q,T)$ converges in the domain of convergence of $Y_A(QT,T)$ and
is invertible in that domain, since $Y_A(x,y)$ is. By
admissibility, there exists $|q-1|\max(|c_0|,R_0)<R\leq r_X$ such
that $Y_A(x,y)$ converges for all $(x,y)\in\mathcal{U}_R$, i.e.
for all $(x,y)$ such that $|x-y|<R$ (cf. \eqref{U_R}). Then
$Y_A(QT,T)$ converges for $|Q-1||T|<R$. Since $|T|\leq \sup_{c\in
A}|c|=\max(|c_0|,R_0)$, it follows that $Y(QT,T)$ converges for
$|Q-1|< R/\max(|c_0|,R_0)$.\end{proof}

\begin{remark}\label{admissible ==> analytic}
By the propagation Theorem, every object of
$\sigma-\Mod(\H_K(X))^{\mathrm{adm}}_{U}$ and of
$(\sigma,\delta)-\Mod(\H_K(X))^{\mathrm{adm}}_{U}$ is
automatically analytic.
\end{remark}

\begin{corollary}\label{compare with}
Let $ \max(|c_0|,R_0) < r \leq  r_X $, and let
$S\subseteq\mathrm{D}^-(1,r/\max(|c_0|,R_0))$, such that
$S^\circ\neq\emptyset$. For all $q\in S^\circ$ one has the
following diagram in which all functors are equivalences by Remark
\ref{shhhshhsssuu}:
\begin{equation}\label{diagram ..R..RR}
\xymatrix{ \sigma-\Mod(\H_K(X))^{[r]}_{S}
\quad\ar@{=}[r]^-{\eqref{sigma-an=(sigma,delta)-an}}
\ar[d]_{\mathrm{Res}^{\;S}_{\{q\}}}^{\wr}
\ar@{}[dr]|{\odot}& %
\quad (\sigma,\delta)-\Mod(\H_K(X))^{[r]}_{S}
\ar[d]_{\wr}^{\mathrm{Res}^{\;S}_{\{q\}}} \\
\sigma_q-\Mod(\H_K(X))^{[r]}\quad
&\quad(\sigma_q,\delta_q)-\Mod(\H_K(X))^{[r]}
\ar[l]^-{\textrm{\emph{Forget} }\delta_q}_-{\sim}\;.}
\end{equation}
By considering the union for all $r$ (cf. Equation \eqref{Remark -
ffff}) one has the following statement. If
$\tau_q:=|q-1|\max(|c_0|,R_0)$, one then has the equivalences:
\begin{equation}
\xymatrix{
\begin{displaystyle}
\bigcup_{r>\tau_q}\end{displaystyle}\sigma-\Mod(\H_K(X))^{\mathrm{adm}}_{\mathrm{D}^-(1,r)}
\quad\ar@{=}[r]^-{\eqref{sigma-an=(sigma,delta)-an}}
\ar[d]_{\bigcup_{r>\tau_q}\mathrm{Res}^{\mathrm{D}^-(1,r)}_{\;\{q\}}}^{\wr}
\ar@{}[dr]|{\odot}& %
\quad\begin{displaystyle}\bigcup_{r>\tau_q}\end{displaystyle}
(\sigma,\delta)-\Mod(\H_K(X))^{\mathrm{adm}}_{\mathrm{D}^-(1,r)}
\ar[d]_{\wr}^{\bigcup_{r > \tau_q}\mathrm{Res}^{\mathrm{D}^-(1,r)}_{\;\{q\}}} \\
\sigma_q-\Mod(\H_K(X))^{\mathrm{adm}}
&\quad(\sigma_q,\delta_q)-\Mod(\H_K(X))^{\mathrm{adm}}
\ar[l]^{\textrm{\emph{Forget} }\delta_q}_{\sim}\;.}
\end{equation}

In particular, if $q,q'\in\mathrm{D}^-(1,1)-\bs{\mu}_{p^{\infty}}$
verify $\max(|q-1|,|q'-1|)\max(|c_0|,R_0) < r$, then, by the
formalism introduced in Section \ref{constant deformation
section}, if $\mathrm{D}:=\mathrm{D}^-(1,r/\max(|c_0|,R_0))$, one
has an equivalence:
\begin{equation}\label{Deformation}
\mathrm{Res}_{q'}^{\mathrm{D}}\circ(\mathrm{Res}^{\mathrm{D}}_{q})^{-1}\;:\;
\sigma_q-\Mod(\H_K(X))^{[r]}\xrightarrow[]{\;\;\sim\;\;}
\sigma_{q'}-\Mod(\H_K(X))^{[r]}\;.
\end{equation}
The same statement holds for $(\sigma,\delta)$-modules without
assuming $q,q'\notin\bs{\mu}_{p^{\infty}}$. \hfill\CVD
\end{corollary}

\begin{definition}
In the notation of Corollary \ref{compare with} (cf. Equation
\eqref{Deformation}), if $q,q'\notin\bs{\mu}_{p^{\infty}}$, we set
\begin{equation}
\mathrm{Def}_{q,q'}^{\mathrm{Tay}}:=\mathrm{Res}_{q'}^{\mathrm{D}}\circ(\mathrm{Res}^{\mathrm{D}}_{q})^{-1}
\;:\;\sigma_q-\Mod(\H_K(X))^{[r]}\xrightarrow[]{\;\;\sim\;\;}
\sigma_{q'}-\Mod(\H_K(X))^{[r]} \;.
\end{equation}
We denote again by $\mathrm{Def}_{q,q'}^{\mathrm{Tay}}$, without
assuming $q,q'\in\bs{\mu}_{p^{\infty}}$, the analogous functor for
$(\sigma,\delta)$-modules. Moreover, if
$q\notin\Q(X)-\bs{\mu}_{p^{\infty}}$, then we set
\begin{equation}
\Conf_q^{\mathrm{Tay}}:=\mathrm{Def}_{q,1}^{\mathrm{Tay}}\circ(\mathrm{Forget}\;\delta_q)^{-1}\;:\;
\sigma_q-\Mod(\H_K(X))^{[r]}\xrightarrow[]{\;\;\sim\;\;}\delta_{1}-\Mod(\H_K(X))^{[r]}\;.
\end{equation}
\end{definition}

By remark \ref{Remark explaining the Propagation}, the functor
$\Conf^{\mathrm{Tay}}_q:(\sigma_q,\delta_q)-\Mod(\H_K(X))^{[r]}
\xrightarrow[]{\;\sim\;} \sigma_q-\Mod(\H_K(X))^{[r]}$ of diagram
\eqref{diagram ..R..RR} coincides with $\Conf^{\C}_q$ (cf.
Definition \ref{Conf_q,q'}), where $\C$ is equal to $\a_K(c,r)$,
where $r$ is as in the corollary \ref{compare with}, and where
$c\in X(K)$ is arbitrarily chosen.

\subsubsection{Root of unity.}
If $q\in\bs{\mu}_{p^{\infty}}$, then the categories
$\sigma_q-\Mod(\H_K(X))_S^{[r]}$ and
$\sigma_q-\Mod(\H_K(X))_S^{\mathrm{adm}}$ are not defined. In this
case we cannot expect any  equivalence between
$(\sigma_q,\delta_q)-\Mod(\H_K(X))^{\mathrm{adm}}$ with a full
subcategory of $\sigma_q-\Mod(\H_K(X))$ because the first category
is $K$-linear and the second is not. In this case we will see in
Proposition \ref{frobenius are trivial in the root of 1} that the
functor ``$\textrm{Forget}\;\delta_q$'' is not very interesting
since it sends every $(\sigma_q,\delta_q)-$module with Frobenius
structure into the trivial $\sigma_q-$module (i.e., a direct sum
of the copies of the unit object).

\subsubsection{}\label{explicit computation of the derivation as limit}
Starting from a Taylor \emph{admissible}
$\sigma_q-$module $\M$ over $\B$, one can \emph{compute} the
differential equation
$\mathrm{Conf}_q^{\mathrm{Tay}}(\M)\in\delta_1-\Mod(\B)$ by the
relation
\begin{equation}
G(1,T)=\lim_{q\to 1}\frac{A(q,T)-\mathrm{Id}}{q-1}=\lim_{n\to
+\infty}\frac{A(q^{p^n},T)-\mathrm{Id}}{q^{p^n}-1}\;,
\end{equation}
where $A(q^{p^n},T)=A(q,q^{p^n-1}T)A(q,q^{p^n-2}T)\cdots A(q,T)$.
The propagation theorem provides the convergence of this limit in
$M_n(\B)$; The reader may have the feeling that this limit should
be easy to compute, but (without introducing the Taylor solution)
the convergence of this limit and its explicit computation are
\emph{highly non trivial facts}. It is surprising to see that the
admissibility condition, which is not a strong assumption,
actually implies such a deep fact.

\begin{remark}
It should be possible to generalize the main theorem to other kind
of operators, different from $\sigma_q$. In other words it should
be possible to ``deform'' differential equations into
``$\sigma-$difference equations'', where $\sigma$ in an
automorphism different from $\sigma_q$, but sufficiently close to
the identity. In a future work we will study the action of a
$p-$adic Lie group on differential equations.
\end{remark}

\subsection{Extending the Confluence Functor to the case $|q-1|=|q|=1$}

Let $q\in\Q(X)-\bs{\mu}(\Q(X))$ be such that $q^{k_0}\in\Q_1(X)$,
for some $k_0\geq 1$.\footnote{For an annulus centered at $0$, the
condition $q^{k_0}\in\Q_1(A)=\mathrm{D}^-(1,1)$ is equivalent to
$\bar{q}\in\mathbb{F}_p^{\mathrm{alg}}$.} By composing with the
evident functor
\begin{equation}
\sigma_q-\Mod(\H_K(X))\xrightarrow[\qquad]{}\sigma_{q^{k_0}}-\Mod(\H_K(X))\;,
\end{equation}
one defines ``\emph{$k_0$-Taylor admissible objects}'' of
$\sigma_q-\Mod(\H_K(X))$ as objects whose image is Taylor
admissible in $\sigma_{q^{k_0}}-\Mod(\H_K(X))$. Since the sequence
$\{q^{k_0p^n}\}_{n\geq 0}$ tends to $1$, then, for $k_0$
sufficiently large, $q^{k_0}$ satisfies the condition of section
\ref{k_0}, in order that $d_{k^0}$ verifies equality
\eqref{d_q(f)leq r_A^-1 f, fifi}. We obtain then a Confluence
Functor:
\begin{equation}
\sigma_q-\Mod(\H_K(X))^{k_0-\mathrm{adm}}\xrightarrow[\qquad]{}
\delta_1-\Mod(\H_K(X))^{\mathrm{adm}}\;.
\end{equation}

The converse of this fact (i.e. the deformation of a differential
equation into a $q-$difference equation with $|q|=1$ and $|q-1|$
large) remains an open problem.

\begin{remark} \label{k_0-example} Notice that there exist equations in $\sigma_q-\Mod(\H_K(X))$ which
are not $k_0$-Taylor admissible, for all $k_0\geq 1$. %
For example consider the rank one equation $\sigma_q - a$, with
$a\in K$, $|a|>1$. Suppose also that $|q-1|<|p|^{\frac{1}{p-1}}$,
in order that $\liminf_n|[n]_q^!|^{1/n}=|p|^{\frac{1}{p-1}}$. Then
the radius is small and one can compute it explicitly by applying
\cite[Prop.4.6]{DV-Dwork}. One has
$Ray((\M,\sigma_q^{\M}),\rho)=|a|^{-1}|p|^{\frac{1}{p-1}}|q-1|\rho<|q-1|\rho$,
and
$Ray((\M,\sigma_{q^{k_0}}^{\M}),\rho)=|a|^{-k_0}|p|^{\frac{1}{p-1}}|q^{k_0}-1|\rho<|q^{k_0}-1|\rho$.
\end{remark}

\subsection{Propagation Theorem over $\Hd_K$ and $\R_K$}\label{Propagation Theorem for other rings}

The Propagation Theorem is true over every base ring $\B$
appearing in this paper, up to a correct definition for the notion
of ``Taylor admissible''. We state here the results for $\Hd_K$
and $\R_K$.

\begin{proposition}\label{solvable extend to all the disk}
Let again $\B:=\Hd_K$, or $\B:=\R_K$, let  $0<r\leq 1$, and let
$S\subseteq\mathrm{D}^-(1,r)$, be a subset, with
$S^{\circ}\neq\emptyset$. Let $\M \in \sigma-\Mod(\B)_S^{[r]}$
(i.e. in particular $\M$ is admissible). Then $\M$ is the
restriction to $S$ of an \emph{analytically $\C$-constant} module
over all the disk $\mathrm{D}^-(1,r)$. Moreover, the restriction
functor is an equivalence:
\begin{equation}
\sigma-\Mod(\B)^{[r]}_{\mathrm{D}^-(1,r)}
\xrightarrow[\sim]{\mathrm{Res}^{\mathrm{D}^-(1,r)}_{S}}
\sigma-\Mod(\B)^{[r]}_S\;.
\end{equation}
In particular solvable modules extend to the whole disk
$\mathrm{D}^{-}(1,1)$. The analogous assertion holds for
$(\sigma,\delta)-$modules, without supposing that $S^{\circ}\neq
\emptyset$:
\begin{equation}(\sigma,\delta)-\Mod(\B)^{[r]}_{\mathrm{D}^-(1,r)}
\xrightarrow[\sim]{\mathrm{Res}^{\mathrm{D}^-(1,r)}_{S}}
(\sigma,\delta)-\Mod(\B)^{[r]}_{S}\;.
\end{equation}
\end{proposition}
\begin{proof} 
By Lemma \ref{The forget from S to S' is fully faithful}, it
suffices to prove the essential surjectivity of
$\mathrm{Res}^{\mathrm{D}^-(1,r)}_S$. The proof is straightforward
and essentially the same as the proof of the Propagation Theorem
\ref{main theorem}.\end{proof}

\begin{corollary}\label{propagation confluence and def over Robba}
Let $q,q'\in\mathrm{D}^-(1,1)-\bs{\mu}_{p^{\infty}}$. Let
$r\in\mathbb{R}$ satisfy
\begin{equation}
\max(|q-1|,|q'-1|)<r\leq 1\;.
\end{equation}
Then one has an equivalence
\begin{equation}
\sigma_q-\Mod(\R_K)^{[r]}\xrightarrow[\sim]{\;\;\Def_{q,q'}^{\mathrm{Tay}}\;\;}\sigma_{q'}-\Mod(\R_K)^{[r]}\;.
\end{equation}
The same equivalence holds between
$(\sigma_q,\delta_q)-\Mod(\R_K)^{[r]}$ and
$(\sigma_{q'},\delta_{q'})-\Mod(\R_K)^{[r]}$, without assuming
$q\notin \bs{\mu}_{p^\infty}$. Moreover, if
$q\notin\bs{\mu}_{p^\infty}$, and if $|q-1|<r$, then we have an
equivalence
\begin{equation}
(\sigma_q,\delta_q)-\Mod(\R_K)^{[r]}
\xrightarrow[\sim]{\textrm{``}\mathrm{Forget
}\;\delta_q\textrm{''}} \sigma_q-\Mod(\R_K)^{[r]}.
\end{equation}
As usual we set
$\Conf_q^{\mathrm{Tay}}:=\mathrm{Def}_{}^{\mathrm{Tay}}\circ(\textrm{Forget
}\delta_q)^{-1}$. The analogous statement holds for
$\Hd_K$.\hfill\CVD
\end{corollary}

\subsubsection{Unipotent equations}
\label{unipotent in detail}
We shall compute the deformation $\Def^{\mathrm{Tay}}_{1,q}$ of
the differential module $U_m$ defined by the equation
\begin{equation}
\delta_1(Y_{U_m})=\left(
\begin{smallmatrix}
0&1&0&\cdots&0\\
0&0&1&\cdots&0\\
&&&\begin{picture}(0,10)\put(-5,0){\begin{tiny}$\ddots$\end{tiny}}\end{picture}&\\
0&0&0&\cdots&1\\
0&0&0&\cdots&0\\
\end{smallmatrix} \right)\cdot Y_{U_m}\;,\qquad
Y_{U_m}(x,y)=\left(\begin{smallmatrix}
1&\ell_1&\cdots&\ell_{m-2}&\ell_{m-1}\\
0&1&\ell_1&\cdots&\ell_{m-2}\\
\begin{picture}(0,10)\put(-1,1){\begin{tiny}$\vdots$\end{tiny}}\end{picture}&&&&
\begin{picture}(0,10)\put(-1,1){\begin{tiny}$\vdots$\end{tiny}}\end{picture}\\
0&\cdots&0&1&\ell_1\\
0&\cdots&0&0&1\\
\end{smallmatrix}
\right) \;,
\end{equation}
where $\ell_n := [\log(x)-\log(y)]^n / n!$. One has
$\sigma_q^x(\ell_n(x,y))=[\log(qx)-\log(y)]^n/n!=
(\log(q)+\log(x)-\log(y))^n/n!=\sum_{i=0}^{n}\frac{\log(q)^{n-k}}{(n-k)!}\cdot\ell_k$.
The matrix of $\sigma_q^{U_m}$ is then
\begin{equation}
A(q,T)=\left(
\begin{smallmatrix}
1&\;\log(q)&\frac{\log(q)^2}{2}
&\cdots&\frac{\log(q)^{m-1}}{(m-1)!}\\
0&1&\log(q)
&\cdots&\frac{\log(q)^{m-2}}{(m-2)!}\\
\vdots&&
&&\vdots\\
0&0&\cdots
&1&\log(q)\\
0&0&\cdots
&\cdots&1\\
\end{smallmatrix} \right)\;\;.
\end{equation}

\subsection{Classification of solvable rank one $q-$difference equations over $\R_{K_\infty}$}
\label{classific Rk1}

In this section we classify rank one solvable $q-$difference
equations over $\R_{K_\infty}$ by applying the deformation
$\mathrm{Def}^{\mathrm{Tay}}_{1,q}$ to the classification of the
differential equations obtained in \cite{Rk1}. We recall the
classification of the rank one solvable differential equations
over $\R_{K_\infty}:=\cup_{s\geq 0}\R_{K_s}$ (see below).

We fix a Lubin-Tate group $\mathfrak{G}_P$ isomorphic to
$\widehat{\mathbb{G}}_m$ over $\mathbb{Z}_p$. We recall that
$\mathfrak{G}_P$ is defined by an uniformizer $\mathrm{w}$ of
$\mathbb{Z}_p$, and by a series $P(X)\in X\mathbb{Z}_p[[X]]$,
satisfying $P(X)\equiv \mathrm{w}\cdot
X\pmod{X^{2}\mathbb{Z}_p[[X]]}$ and $P(X)\equiv
X^p\pmod{p\mathbb{Z}_p[[X]]}$. Such a formal series is called a
\emph{Lubin-Tate series}. We fix now a sequence
$\bs{\pi}:=(\pi_m)_{m\geq 0}$,
$\pi_m\in\mathbb{Q}_p^{\mathrm{alg}}$, such that $P(\pi_0)=0$,
$\pi_0\neq 0$ and $P(\pi_{m+1})=\pi_m$, for all $m\geq 0$. The
element $(\pi_m)_{m\geq 0}$ is a generator of the Tate module of
$\mathfrak{G}_P$ which is a free rank one $\mathbb{Z}_p-$module.
We set $K_s:=K(\pi_s)$ and $K_{\infty}:=\cup_{s\geq 0}K_s$. We
denote by $k_s$ and $k_\infty$ the respective residual fields. The
tower $K\subseteq K_0\subseteq K_1\subseteq\ldots$ does not depend
n the choice of $\bs{\pi}$, nor on
$\mathfrak{G}_P\cong\widehat{\mathbb{G}}_m$. One has
$K_s=K(\xi_{s})$, where $\xi_s$ is a primitive $p^{s+1}$-th root
of unity. For example, one can choose
$\mathfrak{G}_P=\mathbb{G}_m$, hence $P(X)=(X+1)^p-1$, and $\pi_m
= \xi_{m}-1$, where $\xi_m$ is a compatible sequence of primitive
$p^{m+1}-$th root of $1$, i.e. $\xi_0^p=1$ and
$\xi_m^p=\xi_{m-1}$, for all $m\geq 0$. One has the following
facts:
\begin{enumerate}
\item Every rank one \emph{solvable} differential module over
$\R_{K}$ has a basis in which the associated operator is
\begin{equation}\label{L(a_0,f^-(T))}
\L(a_0,\bs{f}^-(T))\;:=\;\delta_1-\Bigl(a_0 -
\sum_{j=0}^{s}\pi_{s-j}\sum_{i=0}^{j}f_{i}^-(T)^{p^{j-i}}
\partial_{T,\log}(f_{i}^-(T)) \Bigr)\;,
\end{equation}
where $a_0\in\mathbb{Z}_p$, and
$\bs{f}^-(T):=(f_0^-(T),\ldots,f_s^-(T))$ is a Witt vector in
$\W_s(T^{-1}\O_{K_s}[T^{-1}])$, with $K_s:=K(\pi_s)$. Notice that
even if $\pi_j$ does not belong to $K$, the resulting polynomial
$\sum_{j=0}^{s}\pi_{s-j}\sum_{i=0}^{j}f_{i}^-(T)^{p^{j-i}}
\partial_{T,\log}(f_{i}^-(T))$ has, by assumption, coefficients in $K$. %

\item The Taylor solution at $\infty$ of the differential module
in this basis is given by the so called $\bs{\pi}$-exponential
attached to $\bs{f}^-(T)$:
\begin{equation}\label{phi^-(T)}
T^{a_0}\cdot\et_{p^m}(\bs{f}^-(T),1)\;:=\;
T^{a_0}\cdot\exp\Bigl(\sum_{j=0}^s\pi_{s-j}\frac{\phi_j^-(T)}{p^j}\Bigr)\;,
\end{equation}
where
$\ph{\phi_0^-(T),\ldots,\phi_s^-(T)}\in(T^{-1}\O_{K_s}[T^{-1}])^{s+1}$
is the phantom vector of $\bs{f}^-(T)$, namely one has
$\phi_j^-(T)=\sum_{i=0}^jp^if_i^-(T)^{p^{j-i}}$.

\item \label{isomorphism class} The correspondence
$\bs{f}^-(T)\mapsto\et_{p^s}(\bs{f}^-(T),1)$ is a group morphism
\begin{equation}
\W_s(T^{-1}\O_{K_s}[T^{-1}])\xrightarrow[]{\et_{p^s}(-,1)}
1+\pi_sT^{-1}\O_{K_s}[[T^{-1}]]\;.\end{equation} %
Notice that if $L(0,\bs{f}^-(T))$ has its coefficients in $\R_K$
($\subset\R_{K_s}$) then also $\et_{p^s}(\bs{f}^-(T),1)$ lies in
$1+T^{-1}\O_{K}[[T^{-1}]]$ (because it is its Taylor solution at
$\infty$).

\item Conversely, $\L(a_0,\bs{f}^-(T))$ is solvable for all pairs
$(a_0,\bs{f}^-(T))\;\in\;
\mathbb{Z}_p\times \W_s(T^{-1}\O_{K_s}[T^{-1}])$.%

\item The operator $\L(a_0,\bs{f}^-(T))$ has a (strong) Frobenius
structure (cf. Def. \ref{def of frob structure - order}) if and
only if
$a_0\;\in\;\mathbb{Z}_{(p)}:=\mathbb{Z}_{p}\cap\mathbb{Q}$.

\item The operators $\L(a_0,\bs{f}^-_1(T))$ and
$\L(b_0,\bs{f}_2^-(T))$ (with coefficients in
$\R_K$($\subset\R_{K_s}$)) define isomorphic differential modules
(over $\R_K$) if and only if $a_0-b_0\in\mathbb{Z}$ and the
Artin-Schreier equation
\begin{equation}\label{asW}
\Fb(\overline{\bs{g}^-(T)})-\overline{\bs{g}^-(T)}=\overline{\bs{f}^-_1(T)-\bs{f}^-_2(T)}
\end{equation}
has a solution $\overline{\bs{g}^-(T)}$ in
$\W_s(k^{\mathrm{alg}}(\!(t)\!))$, where $t$ is the reduction of
$T$, and $\Fb$ is the Frobenius of
$\W_s(k^{\mathrm{alg}}(\!(t)\!))$ (sending
$(\bar{g}_0,\ldots,\bar{g}_s)$ into
$(\bar{g}_0^p,\ldots,\bar{g}_s^p)$). This happens if and only if
the equation $L(0,\bs{f}^-_1(T)-\bs{f}^-_2(T))$ is trivial over
$\R_K$, and also if and only if
$\et_{p^s}(\bs{f}^-_1(T)-\bs{f}^-_2(T),1)$ is
overconvergent.\footnote{Indeed the overconvergence of
$\et_{p^s}(\bs{f}^-_1(T)-\bs{f}^-_2(T),1)$ is independent on the
residual field, for this reason we can look for solution of the
Artin-Schreier-Witt equation \eqref{asW} with coefficients in the
more general field $k^{\mathrm{alg}}$ instead of $k$.}

\end{enumerate}

By deformation, every solvable $q$-difference equation, with
$|q-1|<1$, has a solution at $\infty$ of the form $T^{a_0} \cdot
\et_{p^s}(\bs{f}^-(T),1)$. Its matrix in this basis is then
$$A(q,T)=\et_{p^s}(\bs{f}^-(qT),1)/\et_{p^s}(\bs{f}^-(T),1)
=\et_{p^s}(\bs{f}^-(qT)-\bs{f}^-(T),1)\;.$$ %
The deformation guarantees that $A(q,T)\in\R_K$. This is confirmed
by the fact that $\bs{f}^-(qT)$ and $\bs{f}^-(T)$ have the same
reduction in $\W_s(k^{\mathrm{alg}}(\!(t)\!))$, and hence
$\et_{p^s}(\bs{f}^-(qT)-\bs{f}^-(T),1)\in\R_K$ by  point vi) of
the previous classification.

\section{Quasi unipotence and $p-$adic local monodromy theorem}

\label{section - quasi unipotence and plmt}

In this section we show how to deduce the $q$-analogue of the
$p$-adic local monodromy theorem (cf. \cite{An}, \cite{Ked},
\cite{Me}) by deformation.

Let $K$ be a complete discrete valued field with perfect residue
field (this hypothesis is necessary to have the $p$-adic local
monodromy theorem). Let $\Ed_K\subset\R_K$ be the so called
\emph{bounded Robba ring}, $\Ed_K:=\{\sum_{i\in\mathbb{Z}}
a_iT^i\in\R_K\;|\; \sup|a_i|<+\infty,\lim_{i\to -\infty}
|a_i|=0\}$. Then, since $K$ is discrete valued,
$(\Ed_K,|\cdot|_{(0,1)})$ is a \emph{Henselian} valued field, with
residue field $k(\!(t)\!)$. It has two topologies arising from
$|\cdot|_{(0,1)}$, and from the inclusion in $\R_K$. It is not
complete with respect to none of these two topologies, but $\Ed_K$
is dense in $\R_K$. One has the inclusions
\begin{equation}
\Hd_K\;\;\subset\;\;\Ed_K\;\;\subset\;\;\R_K\;.
\end{equation}
\vspace{-1cm}
\subsection{Frobenius Functor and Frobenius Structure}
Let $\varphi:K\to K$ be an absolute Frobenius (i.e. a ring
morphism lifting of the $p-$th power map of $k$). Since $\R_K$ is
not a local ring, and does not have a residue ring, we need a
special definition:
\begin{definition}
An absolute \emph{Frobenius} on $\R_K$ (resp. $\Hd_K$, $\Ed_K$) is
a continuous ring morphism, again denoted by
$\varphi:\R_K\to\R_K$, extending $\varphi$ on $K$ and such that
$\varphi(\sum a_iT^i)=\sum \varphi(a_i)\varphi(T)^i$, where
$\varphi(T)=\sum_{i\in\mathbb{Z}}b_iT^i\in\R_K$ (resp.
$\varphi(T)\in\Hd_K$, $\varphi(T)\in\Ed_K$) verifies $|b_i|<1$,
for all $i\neq p$, and $|b_p-1|<1$.
\end{definition}
\begin{definition}
We denote by $\phi$ the particular absolute Frobenius on $\R_K$
given by the choice
\begin{equation}\label{phi(T)=T^p}\index{phi(T)@$\phi(T):=T^p$}
\phi(T):=T^p\;,\qquad \phi(f(T)):=f^{\varphi}(T^p)\;.
\end{equation}
where $f^\varphi(T)$ is the series obtained from $f(T)$ by
applying $\varphi:K\to K$ on the coefficients.
\end{definition}
Let $\B$ be one of the rings $\Hd_K$, $\Ed_K$, or $\R_K$. For all
$q\in\mathrm{D}^-(1,1)$, the following diagrams are commutative
\begin{equation}
\xymatrix{
\B\ar@{}[rd]|{\odot}\ar[r]^{\phi}\ar[d]_{\sigma_{q^p}}&\B
\ar[d]^{\sigma_q}\ar@{}[rrd]|{;}&
&\B\ar@{}[rd]|{\odot}\ar[d]_{p\cdot\delta_1}\ar[r]^{\phi}&\B\ar[d]^{\delta_1}&\\
\B\ar[r]_{\phi}&\B &&\B\ar[r]_{\phi}&\B & .}
\end{equation}

\begin{definition}[(Frobenius functor)]\label{Frobenius Functor} Let
$S\subseteq\mathrm{D}^-(1,r)$, $0<r\leq 1$. Let
\begin{equation}\label{r'=min( r^1/p , r|p| )}
r':=\min(\;r^{1/p}\;,\;r\cdot|p|^{-1}\;)\;.
\end{equation}
The Frobenius functor (cf. def. \ref{sigma-Mod(R)^[r]_S})
\begin{eqnarray}\index{phi^*@$\phi^*$}
\phi^*:&(\sigma,\delta)-\Mod(\B)^{[r]}_{S}\xrightarrow[]{\;\;\quad\;\;}
(\sigma,\delta)-\Mod(\B)^{[r']}_S\;,& \\
(\textrm{resp.
}\;\phi^*:&\phantom{(\delta,)}\sigma-\Mod(\B)^{[r]}_{S}\xrightarrow[]{\;\;\quad\;\;}
\phantom{(\delta,)}\sigma-\Mod(\B)^{[r']}_S\;)&
\end{eqnarray}
is defined as $\phi^*(\M,\sigma^{\M},\delta_1^{\M})=
(\;\phi^*(\M)\;,\;\sigma^{\phi^*(\M)}\;,\;\delta_1^{\phi^*(\M)}\;)$,
where
\begin{enumerate}
\item $\phi^*(\M):=\M\otimes_{\B,\phi}\B$ is the scalar extension
of $\M$ via $\phi$, %
\item the morphism $\sigma^{\phi^*(\M)}$ is given by
$\sigma_q^{\phi^*(\M)}=\sigma_{q^p}^{\M}\otimes\sigma_q^{\B}$:
\begin{equation}
q\longmapsto\sigma_{q^p}^{\M}\otimes\sigma_q\;\;:\;\;S\xrightarrow[
]{\;\;\sigma^{\phi^*(\M)}\;\;}\Autc(\phi^*(\M))\;,
\end{equation}
\item the derivation is given by
\begin{equation}
\delta_1^{\phi^*(\M)}=(p\cdot\delta_1^{\M})\otimes\mathrm{Id_{\B}}+\mathrm{Id}_{\M}\otimes\delta_1^{\B}\;,
\end{equation}
\item a morphism $\alpha:\M\to\N$ is sent into $\alpha\otimes
1:\phi^*(\M)\to\phi^*(\N)$.
\end{enumerate}
\end{definition}
\begin{remark}\label{Remark on the definition of the
pull-back-***} The fact that the functor $\phi^*$ sends
$(\sigma,\delta)-\Mod(\B)^{[r]}_{S}$ into
$(\sigma,\delta)-\Mod(\B)^{[r']}_{S}$ with this particular value
of $r'$ (cf. Equation \eqref{r'=min( r^1/p , r|p| )}) results from
the fact that this result is true for \emph{differential}
equations (cf. \cite[Appendix]{Pu}, and \cite[Prop.7.2]{Ch-Me}),
and from the confluence.
\end{remark}

\subsubsection{}\label{Remark on the definition of the pull-back by
Frobenius} We observe that the pull-back $\varphi^*(\M)$ is
actually a $\sigma$-module over $S^{1/p}:=\{q\in K \;|\; q^p\in S
\}$. Indeed $\phi^*(\M)$ is canonically endowed with the action of
$\sigma_{q^{1/p}}^{\phi^*(\M)}:=\sigma_{q}^{\M}\otimes\sigma_{q^{1/p}}:\phi^*(\M)\to\phi^*(\M)$,
for all roots $q^{1/p}$ of $q$. This fact was used in \cite{An-DV}
to define the so called confluent weak Frobenius structure (cf.
Definition \ref{confluent weak Frob Str}).

If $\M\in (\sigma,\delta)-\Mod(\Hd_K)^{[r]}_S$, then we can
consider its Taylor solution at $1$: $Y(T,1)=\sum_{i\geq
0}Y_i(T-1)^i\in GL_n(\a_K(1,1))$, $Y_i\in M_n(K)$. Then the Taylor
solution of $\phi^*(\M)$ is given by
\begin{equation}
Y^{\phi}(T^p,1)\;\;:=\;\;\sum_{i\geq 0}\varphi(Y_i)(T^p-1)^i\;.
\end{equation}

The matrices of $\phi^*(\sigma_q)$ and $\phi^*(\delta_1)$ are the
following. Let $\e=\{\mathrm{e}_1,\ldots,\mathrm{e}_n\}$ be a
basis of $\M$. Let $\sigma_q-A(q,T)$ and $\delta_1-G(1,T)$ be the
operators associated to $\sigma_q^\M$ and $\delta_1^{\M}$ in this
basis. Then the operators associated to $\phi^*(\M)$ in the basis
$\e\otimes 1$ are
\begin{equation}
\sigma_q - A^\varphi(q^p,T^p) \;,\qquad \delta_1-p\cdot
G^\varphi(1,T^p)\;,
\end{equation}
where, according with \eqref{A(q q',T)=A(q',qT) A(q,T)}, one has
$A(q^p,T)=A(q,q^{p-1}T)\cdots A(q,qT)A(q,T)$.

\subsubsection{Frobenius Structure.} The functor
$\phi^*:\delta_1-\Mod(\R_K)^{[1]}\xrightarrow[]{\;\;\sim\;\;}\delta_1-\Mod(\R_K)^{[1]}$
is an equivalence (cf. \cite[Cor.8.14]{Ch-Me}). By deformation
$\phi^*$ is hence an auto-equivalence of
$\sigma-\Mod(\R_K)^{[1]}_{S}$ (if $S^{\circ}\neq \emptyset$) and
$(\sigma,\delta)-\Mod(\R_K)^{[1]}_{S}$ (without assuming
$S^{\circ}\neq\emptyset$).

\begin{definition}[(Frobenius structure)]\label{def of frob structure -
order} Let $\B$ be one of the rings $\Hd_K$, $\Ed_K$, or $\R_K$.
Let $S\subseteq\mathrm{D}^-(1,1)$ be a subset. Let $\M$ be a
discrete $\sigma-$module (resp. $(\sigma,\delta)-$module) over
$S$. We will say that $\M$ has a \emph{Frobenius structure of
order $h\geq 1$}, if there exists an $\B$-isomorphism
$(\phi^*)^{h}(\M)\xrightarrow[]{\;\;\sim\;\;}\M$ of
$\sigma$-modules over $S$ (cf. Section \ref{Remark on the
definition of the pull-back by Frobenius}), where
$(\phi^*)^{h}:=\phi^*\circ\cdots\circ\phi^*$, $h-$times. We denote
by
\begin{equation}\index{sigma-Mod(B)_S^phi@$\sigma-\mathrm{Mod}(\B)_S^{(\phi)}$,
$(\sigma,\delta)-\mathrm{Mod}(\B)_S^{(\phi)}$}
\sigma-\mathrm{Mod}(\B)_S^{(\phi)}\;,\qquad(\textrm{resp.}\quad(\sigma,\delta)-\mathrm{Mod}(\B)_S^{(\phi)}\;)
\end{equation}
the full subcategory of $\sigma-\Mod(\B)_S^{[1]}$ (resp.
$(\sigma,\delta)-\mathrm{Mod}(\B)_{S}^{[1]}$) whose objects have a
Frobenius structure of some order.
\end{definition}

If $\M$ has a Frobenius structure, then $r=r'$ (cf. \eqref{r'=min(
r^1/p , r|p| )}) and hence $\M$ is \emph{solvable}:
\begin{equation}
\sigma-\Mod(\B)^{(\phi)}_S\;\subset\;\sigma-\Mod(\B)^{[1]}_S\;.
\end{equation}
Hence objects in $\sigma-\Mod(\B)^{(\phi)}_S$ and
$(\sigma,\delta)-\Mod(\B)^{(\phi)}_S$ are, in particular,
\emph{admissible}.

If $Y(T,1)$ is the Taylor solution of
$\M\in(\sigma,\delta)-\Mod(\Hd_K)_{S}^{[1]}$ at $1$, then the fact
that $\M$ has a Frobenius structure of some order $h\geq 1$, is
equivalent to the existence of a matrix $H(T)\in GL_n(\Hd_K)$ such
that
\begin{equation}\label{Y(T^p) = theta(T) Y(T)}
Y^{\varphi^h}(T^{p^h},1)=H(T)\cdot Y(T,1)\;.
\end{equation}
Indeed $\a_K(1,1)$ is a $\Hd_K$-discrete $\sigma$-algebra over
$\mathrm{D}^-(1,1)$ trivializing $\M$ (cf. Def.\ref{discrete
sigma,delta algebra}). In particular the equivalences
$\mathrm{Def}^{\mathrm{Tay}}_{q,q'}$ and
$\Conf^{\mathrm{Tay}}_{q}$ send objects with Frobenius structure
into objects with Frobenius structure.

\begin{proposition}\label{frobenius are trivial in the root of 1}
Let $\xi$ be a $p^n$-th root of unity, and let
$q\in\Q_1-\bs{\mu}(\Q_1)$. Let
$\M\in\sigma_q-\Mod(\Hd_K)^{(\phi)}$. Then
$\mathrm{Def}_{q,\xi}^{\mathrm{Tay}}(\M)$ is trivial (i.e.
isomorphic to a direct sum of copies of the unit object).
\end{proposition}
\begin{proof} Let $Y(T,1)\in GL_n(\Hd_K)$ be the Taylor solution at $1$ of
$\M$ in some basis $\e$. Then, by \eqref{Y(T^p) = theta(T) Y(T)},
there exists $H(T)$ such that $Y^{\varphi^h}(T^{p^h},1)=H(T)\cdot
Y(T,1)$. Hence, one also has
$Y^{\varphi^{nh}}(T^{p^{nh}},1)=H_n(T)\cdot Y(T,1)$, for some
$H_n(T)\in GL_n(\Hd_K)$. Since
$\sigma_\xi(Y^{\varphi^{nh}}(T^{p^{nh}}))=Y^{\varphi^{nh}}(T^{p^{nh}})$,
it follows that in the basis $H_n(T)\cdot\e$ the matrix of
$\sigma_{\xi}$ is trivial: $A(\xi,T)=\mathrm{Id}$ (cf. Section
\ref{morphisms as solutions}). \end{proof}

\subsection{Special coverings of $\Hd_K$}
\label{Katz-Matsuda equivalence section}

We recall briefly the notions of \emph{special coverings}. The
residue field of $\Ed_K$ is $k(\!(t)\!)$ (with respect to the norm
$|.|_{(0,1)}$). On the other hand, the residue ring of $\Hd_K$
(with respect to the Gauss norm $|.|_{(0,1)}$) is $k[t,t^{-1}]$.
One has
\begin{equation}
\xymatrix{\O_{\Hd_K}\ar@{}[dr]|{\odot}\ar@{}[r]|{\subseteq}\ar[d]&\O_{\Ed_K}\ar[d]\\
k[t,t^{-1}]\ar@{}[r]|{\subseteq}& k(\!(t)\!).}
\end{equation}
We denote by $\O_K[T,T^{-1}]^{\dag}$ \label{O_K T,T^-1 dag} the
weak completion of $\O_K[T,T^{-1}]$, in the sense of Monsky and
Washnitzer (cf. \cite{Monsky-Washnitzer-F1}). One has
\begin{equation}
\Hd_K=\O_K[T,T^{-1}]^{\dag}\otimes_{\O_K} K\;.
\end{equation}
%
%
%
%
%
%
Let us look at the residual situation. The morphism
\begin{equation}
\widehat{\eta}:=\mathrm{Spec}(k(\!(t)\!))\;\;\hookrightarrow\;\;
\mathbb{G}_{m,k}= \mathrm{Spec}(k[t,t^{-1}])
\end{equation}
gives rise, by pull-back, to a map
\begin{equation}
\left\{
\begin{array}{l}
\textrm{Finite \'Etale}\\
\textrm{coverings of }\widehat{\eta}
\end{array}
\right\} \xleftarrow[]{\textrm{Pull-back}}\left\{
\begin{array}{l}
\textrm{Finite \'Etales}\\
\textrm{coverings of }\mathbb{G}_{m,k}
\end{array}
\right\}\;.
\end{equation}
It is known (cf. \cite[2.4.9]{Katz-local-to-global}) that this map
is surjective, and moreover that there exists a full sub-category
of the right hand category, called \emph{special coverings of
$\mathbb{G}_{m,k}$}, which is equivalent, via pull-back, to the
category on the left hand side. Special coverings are defined by
the property that they are tamely ramified at $\infty$, and that
their geometric Galois group has a unique $p-$Sylow subgroup (cf.
\cite[1.3.1]{Katz-local-to-global}).

On the other hand, if $\pi\in\O_K$ is a uniformizing element, then
both $(\O_{\Ed_K},(\pi))$ and $(\O_K[T,T^{-1}]^{\dag},(\pi))$ are
Henselian couples in the sense of \cite[Ch.II]{Ray} (cf.
\cite[5.1]{Matsuda-unipotent}). One can show that the preceding
situation lifts to characteristic $0$. One has the following
equivalences:
\begin{equation}
\xymatrix{ \protect{\left\{
\begin{array}{l}
\textrm{Special}\\
\textrm{extensions of }\Hd_K
\end{array}
\right\}
\ar[r]_-{\sim}^-{-\otimes\Ed_K}\ar@{}[dr]|-{\odot} }&
\protect{\left\{\begin{array}{l}
\textrm{Finite unramified}\\
\textrm{extensions of }\Ed_K
\end{array}
\right\}}\ar[r]_-{\sim}^-{-\otimes\R_K}&\protect{\left\{\begin{array}{l}
\textrm{Special}\\
\textrm{extensions of }\R_K
\end{array}
\right\}}\\
\protect{\left\{
\begin{array}{l}
\textrm{Special extensions}\\
\textrm{of }\O_K[T,T^{-1}]^\dag
\end{array}
\right\} \ar[r]_-{\sim}^-{-\otimes\O_{\Ed_K}}\ar@{}[dr]|-{\odot}
\ar[u]_{\wr}^{-\otimes K}
\ar[d]_{-\otimes k}^{\wr}  }&%
\protect{\left\{\begin{array}{l}
\textrm{Finite onramified}\\
\textrm{extensions of }\O_{\Ed_K}
\end{array}
\right\}\ar[d]^{-\otimes k}_{\wr}\ar[u]_{-\otimes K}^{\wr} } &\\
\protect{\left\{
\begin{array}{l}
\textrm{Special}\\
\textrm{coverings of }\mathbb{G}_{m,k}
\end{array}
\right\} \ar[r]_-{\sim}^-{\textrm{Pull-back}} }& \protect{\left\{
\begin{array}{l}
\textrm{Finite \'etale}\\
\textrm{coverings of }\hat{\eta}
\end{array}
\right\}}& }
\end{equation}
where, by special extension of $\O_K[T,T^{-1}]^{\dag}$ (resp.
$\Hd_K$, $\R_K$) we mean a finite \'etale Galois extension of
$\O_K[T,T^{-1}]^{\dag}$ (resp. $\Hd_K$, $\R_K$) coming, by
Henselianity, from a special cover of $\mathbb{G}_{m,k}$.

\begin{lemma}\label{(S^dag(F))^G=Hd_K}
Let $\F/k(\!(t)\!)$ be a finite Galois extension with Galois group
$\G$. Let $\mathcal{S}^{\dag}(\F)/\Hd_K$ be the corresponding
Special extension of $\Hd_K$. Then
$(\mathcal{S}^{\dag}(\F))^{\G}=\Hd_K$.
\end{lemma}
\begin{proof}
By \cite[Exposé V, Cor.3.4]{SGA-1},
$(\mathcal{S}^{\dag}(\F))^{\G}/\Hd_K$ is a Special extension. The
assertion is then easy since, by the above equivalence there is
bijection between Special sub-algebras of $\mathcal{S}^{\dag}(\F)$
over $\Hd_K$, and sub-extensions of $\F/k(\!(t)\!)$.
\end{proof}

\subsubsection{Extension of $\sigma_q$ to Special extensions.}\label{unique extension of sigma to
Ed'}
\begin{lemma}[\protect{\cite[Section 11.3]{An-DV}}] Let $F/k(\!(t)\!)$ be a
finite separable extension. Let
$\mathcal{F}^{\dag}/\O_K[T,T^{-1}]^{\dag}$ be the corresponding
special extensions. The automorphism $\sigma_q$ of
$\O_K[T,T^{-1}]^{\dag}$ extends to an automorphism
$\mathcal{F}^{\dag}$. The extension is unique up to
$\O_K[T,T^{-1}]^{\dag}$-automorphisms of $\mathcal{F}^{\dag}$. The
same statement holds for the extensions $(\Hd_K)'/\Hd_K$,
$(\Ed_K)'/\Ed_K$, $(\R_K)'/\R_K$ corresponding to $F/k(\!(t)\!)$.
In particular there exists a unique extension of $\sigma_q$ to
$\mathcal{F}^{\dag}$, $(\Hd_K)'$, $(\Ed_K)'$, $(\R_K)'$ inducing
the identity on $F$. 
\end{lemma}
\begin{proof} The proof results from the formal properties of Henselian
couples (cf. \cite{Ray}).
\end{proof}
By uniqueness the extension of $\sigma_q$ commutes with the action
of $\mathrm{Gal}(k(\!(t)\!)^{\mathrm{sep}}/k(\!(t)\!))$.

\begin{remark}
Every finite extensions of $\mathbb{C}(\!(T)\!)$ is of the
form $\mathbb{C}(\!(T^{m/n})\!)$. Up to change the variable we
have an isomorphism $\mathbb{C}(\!(T^{m/n})\!) \cong
\mathbb{C}(\!(Z)\!)$. Analogously it can be seen that a finite
unramified extension of $\Ed_{K}$ is (non canonically) isomorphic
to $\Ed_{K'}$ for some finite $K'/K$. In this case the link
between the variable $Z$ and the variable $T$ is rather complicate
and essentially unknown. \emph{The great problem of the theory is
that the extended automorphism does not send $Z$ into $qZ$.} The
general ``Confluence'' theory introduced in section \ref{theory of
deformation} will be crucial in solving this problem.
\end{remark}
\if{\begin{remark} One can show that every unramified extension
$(\Ed_K)'$ of $\Ed_K$ (resp. special extension $\R_K'$ of $\R_K$)
is \emph{non canonically} isomorphic to $\Ed_{K'}$ (resp.
$\R_{K'}$), for some finite Galois unramified extension $K'/K$.
This is analogous to the classical situation, in which every
extension of $\mathbb{C}(\!(T)\!)$ is of the form
$\mathbb{C}(\!(T^{m/n})\!)$, for some integers $m,n\geq 0$, and so
it is isomorphic to $\mathbb{C}(\!(Z)\!)$, with $T=f(Z):=Z^{n/m}$.
In the case of special extensions of $\Ed_K$ or $\R_K$, the
analogous fact is true: $(\Ed_{K,T})'\cong \Ed_{K',Z}$, for some
new variable $Z$, and a finite unramified extension $K'/K$. But
the relation between the new variable and the old one is highly
non trivial, and essentially unknown. If $k'/k$ is the residue
field of $K'$, and if $t=\overline{f}(z)\in k'(\!(z)\!)$ is the
relation between $t$ and $z$ in characteristic $p$, then the
relation between $T$ and $Z$ in characteristic $0$ is given by
$T=f(Z)\in \O_{\Ed_K}$, where $f(Z)$ is an arbitrary Laurent
series, obtained from $\overline{f}(z)$ by lifting coefficient by
coefficient.
\end{remark}
}\fi

\subsection{Quasi unipotence of differential equations and canonical extension}\label{Special
extensions}

In this section we recall some known facts on p-adic
\emph{differential} equations.

\begin{definition}
We denote by $\widetilde{\Hd_K}$ (resp. $\widetilde{\Ed_K}$,
$\widetilde{\R_K}$) the union of all finite special (resp.
unramified, special) extensions of $\Hd_K$ (resp. $\Ed_K$, $\R_K$)
in an algebraically closure of the field of fractions of $\R_K$.
\end{definition}

\begin{definition}
Let $S\subseteq\mathrm{D}^-(1,1)$ be a subset (resp.
$S\subseteq\mathrm{D}^-(1,1)$, with $S^{\circ}\neq\emptyset$). A
discrete $(\sigma,\delta)-$module on $S$ (resp. discrete
$\sigma-$module on $S$) is called \emph{quasi-unipotent} if it is
trivialized by the discrete $(\sigma,\delta)-$algebra
\begin{equation}\label{R'[log(T)]}
\widetilde{\Hd_K}[\log(T)]\qquad \textrm{(resp. }\;
\widetilde{\Ed_K}[\log(T)]\;,\;\;
\widetilde{\R_K}[\log(T)]\;\;\textrm{)}\;\;.
\end{equation}
\end{definition}

Let $\B:=\Hd_K$, or $\Ed_K$, $\R_K$. We observe that $\M$ is
trivialized by $\widetilde{\B}[\log(T)]$, if and only if $\M$ is
trivialized by $\B'[\log(T)]$, where $\B'$ is a (\emph{finite})
special extension of $\B$. Indeed the entries of a fundamental
matrix of solutions of $\M$ in $\widetilde{\B}[\log(T)]$ all lie
in a finite extension.

\begin{theorem}[($p-$adic local monodromy theorem, cf.
\cite{An},\cite{Ked},\cite{Me})]\label{p-adic
local monodromy theorem}%
 Objects in $\delta_1-\Mod(\R_K)^{(\phi)}$ become
quasi-unipotent possibly after a suitable extension of the field
of constants $K$.  In other words, if
$\M\in\delta_1-\Mod(\R_K)^{(\phi)}$, then there exists a finite
extension $K'/K$ such that $\M\otimes_KK'$ is quasi unipotent
(i.e. trivialized by $\widetilde{\Hd_{K'}}[\log(T)]$).\hfill\CVD
\end{theorem}

\begin{theorem}[(\protect{\cite[7.10,7.15]{Matsuda-unipotent}})]\label{canonical extension}
If a differential equation $\M\in\delta_1-\Mod(\R_K)$ is
quasi-unipotent, then it has a Frobenius structure. Moreover, the
scalar extension functor
\begin{equation}\label{jkjkkkkk}
-\otimes\R_K\;:\;\delta_1-\Mod(\Hd_K)^{(\phi)}\;\;\xrightarrow[]{\qquad}\;\;\delta_1-\Mod(\R_K)^{(\phi)}
\end{equation}
is essentially surjective. \hfill\CVD
\end{theorem}

\begin{theorem}[(\protect{\cite[7.15]{Matsuda-unipotent}})]\label{remk-special objects}
There exists a full sub-category of
$\delta_1-\Mod(\Hd_K)^{(\phi)}$, denoted by
$\delta_1-\Mod(\Hd_K)^{\mathrm{Sp}}$, which is equivalent to
$\delta_1-\Mod(\R_K)^{(\phi)}$ via the scalar extension functor
\eqref{jkjkkkkk}. Objects in $\delta_1-\Mod(\Hd_K)^{\mathrm{Sp}}$
category are trivialized by
$\widetilde{\Hd_K}[\log(T)]$.\hfill\CVD
\end{theorem}

\begin{definition}[(Canonical extension)]\label{CAN-+---}\index{Can@$\mathrm{Can}$}
Objects in $\delta_1-\Mod(\Hd_K)^{\mathrm{Sp}}$ will be called
\emph{special objects}. We will denote by
\begin{equation}\label{Can}
\delta_1-\Mod(\R_K)^{(\phi)}
\xrightarrow[\sim]{\;\;\mathrm{Can}\;\;}\delta_1-\Mod(\Hd_K)^{\mathrm{Sp}}
\;\;\subset \;\;\delta_1-\Mod(\Hd_K)^{(\phi)}
\end{equation}
the section of the functor \eqref{jkjkkkkk}, whose image is the
category of special objects (cf. Theorem \ref{remk-special
objects}). We will call it the \emph{canonical extension functor}.
\end{definition}

\begin{corollary}[(\protect{\cite[7.1.6]{An}})]\label{N otimes U_m}
Let $\M\in\delta_1-\Mod(\R_K)^{(\phi)}$, then, up to replacing $K$
by a finite extension $K'/K$, 
$\M$ decomposes in a direct sum of submodules of the form
$\N\otimes U_m$, where $\N$ is a module trivialized by a special
extension of $\R_K$, and $U_m$ is the $m-$dimensional object
defined by the operator (cf. section \ref{unipotent in detail})
\begin{equation}
\delta_1-\left( \begin{smallmatrix}
0&1&0&\cdots&0\\
0&0&1&\cdots&0\\
&&&\ddots&\\
0&0&0&\cdots&1\\
0&0&0&\cdots&0\\\end{smallmatrix} \right)\;.\qquad\qquad\Box
\end{equation}
\end{corollary}

\begin{remark}
The $\log(T)$ appearing in \eqref{R'[log(T)]}, is added uniquely
to trivialize the module of the form $U_m$, for $m\geq 2$ (cf.
section \ref{unipotent in detail}).
\end{remark}

\begin{lemma}\label{Ed' generated by solutions}
Let $\N\in\delta_1-\Mod(\Hd_K)^{\mathrm{Sp}}$ be a special object
trivialized by $\widetilde{\Hd_{K}}$. Let
$\widetilde{Y}=(\widetilde{y}_{i,j})\in GL_n(\widetilde{\Hd_K})$
be a fundamental matrix solution of $\N$. Let $(\Ed)'$ (resp.
$\R'$) be the smallest special extension of $\Ed_K$ (resp.
$\R_{K}$), such that $\N\otimes\Ed_K$ is trivialized by $(\Ed)'$
($\N\otimes\R_K$ is trivialized by $\R'$). Then one has
\begin{equation}
(\Ed)'=\Ed_{K}[\{\widetilde{y}_{i,j}\}_{i,j}]\;,\qquad
\R'=\R_{K}[\{\widetilde{y}_{i,j}\}_{i,j}]\;.
\end{equation}
In other words, the smallest special extension of $\Ed_K$ (resp.
$\R_K$) trivializing $\M$ is generated by the solutions of $\M$.
\end{lemma}
\begin{proof} Since $\M$ is trivialized by $(\Ed)'$, one has
$\Ed_K[\{\widetilde{y}_{i,j}\}_{i,j}]\subseteq (\Ed)'$. Hence the
differential field $\Ed_{K'}[\{\widetilde{y}_{i,j}\}_{i,j}]$ is an
unramified extension, and is then a special extension. Since
$(\Ed)'$ is minimal,
$\Ed_K[\{\widetilde{y}_{i,j}\}_{i,j}]=(\Ed)'$. The case over
$\R_K$ follows from the case over $\Ed_K$. \end{proof}

\begin{corollary}\label{unique extension of sigma_q to Ed[y_i,j]}
We preserve the notations of Lemma \ref{Ed' generated by
solutions}. There exists a unique $K-$linear ring automorphism
$\sigma_q$ of $\Ed_K[\{\widetilde{y}_{i,j}\}_{i,j}]$, which
induces the identity on the residue field.
\end{corollary}
\begin{proof} By Lemma \ref{Ed' generated by solutions},
$\Ed_K[\{\widetilde{y}_{i,j}\}_{i,j}]$ is a special extension
(i.e. Henselian). Hence, by Section \ref{unique extension of sigma
to Ed'}, the extension of $\sigma_q$ to
$\Ed_K[\{\widetilde{y}_{i,j}\}_{i,j}]$ is unique.\end{proof}

\begin{corollary}\label{canonical extension for sigma and (sigma,delta) modules}
Let $S\subseteq\mathrm{D}_K^-(1,1)$ 
The scalar extension functor
\begin{equation}
-\otimes\R_K\;:\;(\sigma,\delta)-\Mod(\Hd_K)_{S}^{(\phi)}\;\;
\longrightarrow\;\; (\sigma,\delta)-\Mod(\R_K)_{S}^{(\phi)}
\end{equation}
is essentially surjective. Moreover there exists a full
sub-category of
$(\sigma,\delta)-\Mod(\Hd_K)_{S}^{(\phi)}$, 
which we call
$(\sigma,\delta)-\Mod(\Hd_K)_{S}^{\mathrm{Sp}}$, 
equivalent via $-\otimes\R_K$ to
$(\sigma,\delta)-\Mod(\R_K)_{S}^{(\phi)}$. 
The same statement is true for $\sigma$-modules under the
assumption $S^{\circ}\neq \emptyset$.
\end{corollary}
\begin{proof} By Proposition \ref{solvable extend to all the disk}, we can assume that $S=\mathrm{D}^-(1,1)$. By Theorem
\ref{canonical extension} there exists a basis of $\M$ in which
the matrix $G(1,T)$ of $\delta_1^{\M}$ lies in $M_n(\Hd_K)$.
Moreover, $\mathrm{Can}(\M,\delta_1^{\M})$ is Taylor admissible,
since all solvable differential equations are Taylor admissible.
By Proposition \ref{solvable extend to all the disk}, for all
$q\in\mathrm{D}^-(1,1)$, the matrix $A(q,T):=Y_G(qT,T)$ belongs
also to $GL_n(\Hd_K)$. This proves the essential surjectivity. The
fully faithfulness follows by deformation of Thm.
\ref{remk-special objects} (cf. Cor. \ref{compare with}).
\end{proof}

\subsubsection{} It is not clear to us if the smallest
special extension of $\Hd_K$ trivializing a given
$\M\in\delta_1-\Mod(\Hd_K)^{\mathrm{Sp}}$ is generated (over
$\Hd_K$) by the entries of a fundamental matrix of solution of
$\M$. So we are obliged to give the following definition.

\begin{definition}\label{def of C^et_K}
We denote by $\widetilde{\C_K^{\textrm{\'et}}}$ the sub-algebra of
$\widetilde{\Hd_K}$ generated, over $\Hd_K$, by the entries of
every fundamental solutions matrix of each object in
$\delta_1-\Mod(\Hd_K)^{\mathrm{Sp}}$ which is trivialized by
$\widetilde{\Hd_K}$.
\end{definition}

With the notation of Corollary \ref{canonical extension for sigma
and (sigma,delta) modules}, the inclusions
$(\sigma,\delta)-\Mod(\Hd_K)_S^{\mathrm{Sp}}\subset
(\sigma,\delta)-\Mod(\Hd_K)_S^{(\phi)}\subset
(\sigma,\delta)-\Mod(\Hd_K)_S^{[1]}$ are strict (the same holds
for $(\sigma,\delta)$-modules). For example the equation
$\delta_1(y)=a_0y$, with $a_0\in\mathbb{Z}_p-\mathbb{Z}_{(p)}$, is
solvable, but without Frobenius structure (cf. section
\ref{classific Rk1}). On the other hand an object with Frobenius
structure could have non zero $p$-adic slope at $1^+$ (hence
irregular at $\infty$), hence it is not special. Unfortunately we
have no examples of \emph{non special} equations with Frobenius
structure, but trivialized by
$\widetilde{\C^{\textrm{ét}}_K}[\log(T)]$.

\subsection{Quasi unipotence of $\sigma-$modules and $(\sigma,\delta)-$modules with Frobenius structure}
\label{quasi unipotence of q-diff - section}

This section is devoted to prove the following

\begin{theorem}[($p-$adic local monodromy theorem (generalized
form))]\label{p-adic local monodromy theorem (generalized form)}
Let $S\subset\mathrm{D}^-(1,1)$, be a subset (resp.
$S^{\circ}\neq\emptyset$). Then every object
$\M\in(\sigma,\delta)-\Mod(\R_K)^{(\phi)}_S$ (resp.
$\M\in\sigma-\Mod(\R_K)^{(\phi)}_S$) is quasi unipotent, after
replacing $K$, if necessary, by a finite extension $K'/K$
depending on $\M$.
\end{theorem}

This result simplifies, and generalizes the analogous result of
\cite{An-DV}. The proof is obtained by deformation of the $p$-adic
local monodromy theorem of \emph{differential} equation (cf.
Theorem \ref{p-adic local monodromy theorem}).

The proof is essentially the following. Assume that $S=\{q\}$,
with $q\notin\bs{\mu}_{p^{\infty}}$. By canonical extension (cf.
Cor. \ref{canonical extension for sigma and (sigma,delta)
modules}) $\M$ is trivialized by $\widetilde{\R_K}[\log(T)]$ if
and only if $\mathrm{Can}(\M)$ is trivialized by
$\widetilde{\Hd_K}[\log(T)]$ (or equivalently by
$\widetilde{\C^{\textrm{ét}}_K}[\log(T)]$). Hence we can assume
that $\M\in\sigma_q-\Mod(\Hd_K)^{\mathrm{Sp}}$. Firstly apply the
confluence functor to obtain a differential equation
$\Conf_q^{\mathrm{Tay}}(\M,\sigma_q^{\M})$. We prove then in Lemma
\ref{la chiave} below that
$\Conf_q^{\mathrm{Tay}}(\M,\sigma_q^{\M})$ is
$\widetilde{\C^{\textrm{ét}}_K}[\log(T)]$-extensible to
$\mathrm{D}^-(1,1)$ (cf. Def. \ref{definition of strongly
confluent at q}). Hence we obtain, by deformation, another
$q$-difference module
$\Def_{1,q}^{\widetilde{\C^{\textrm{ét}}_K}[\log(T)]}(\Conf_q^{\mathrm{Tay}}(\M,\sigma_q^{\M}))$
over $\Hd_K$ (cf. section \ref{theory of deformation}). This
$q$-difference module is quasi unipotent since, by definition, it
has the same solutions in
$\widetilde{\C^{\textrm{ét}}_K}[\log(T)]$ of the quasi unipotent
differential equation $\Conf_q^{\mathrm{Tay}}(\M,\sigma_q^{\M})$.
We shows then that there is an embedding
$\widetilde{\C^{\textrm{ét}}_K}[\log(T)]\subseteq
\a_{K^{\mathrm{alg}}}(1,1)$ commuting with $\delta_1$, $\varphi$,
and with $\sigma_q$, for all $q\in\mathrm{D}^-(1,1)$ (cf. Lemma
\ref{embedding in a_K(1,1)}). This proves that the restriction of
$\Def_{1,q}^{\mathrm{Tay}}$ to the category of objects trivialized
by $\widetilde{\C^{\textrm{ét}}_K}[\log(T)]$ coincides with
$\Def_{1,q}^{\widetilde{\C^{\textrm{ét}}_K}[\log(T)]}$ (cf.
Sections \ref{C_1 subset C_2 ==> Then
Def^C_1=Def^C_2},\ref{Confluence depends on C}), because
$\Def_{1,q}^{\mathrm{Tay}}=\Def_{1,q}^{\mathrm{C}}$, with
$\C=\a_K(1,1)$ (cf. Remark \ref{Remark explaining the
Propagation}) or equivalently $\C=\a_{K^{\mathrm{alg}}}(1,1)$ (cf.
Corollary \ref{Conf Tay = Conf a_K(1,1) = Conf widetilde R_K}).
Hence
\begin{equation}
\Def_{1,q}^{\widetilde{\C^{\textrm{ét}}_K}[\log(T)]}(\Conf_q^{\mathrm{Tay}}(\M,\sigma_q^{\M}))\;=\;
\Def_{1,q}^{\mathrm{Tay}}(\Conf_q^{\mathrm{Tay}}(\M,\sigma_q^{\M}))\;=\;(\M,\sigma_q^{\M})\;.
\end{equation}
In particular $(\M,\sigma_q^{\M})$ is trivialized by
$\widetilde{\C^{\textrm{ét}}_K}[\log(T)]$ and is hence quasi
unipotent.

\begin{lemma}\label{la chiave}
Let $\M\in\delta_1-\Mod(\Hd_K)^{\mathrm{Sp}}$. Assume that $K$ is
sufficiently large so that $\M$ is quasi unipotent. Let $(\Hd_K)'$
be the smallest special extension of $\Hd_{K}$ such that $\M$ is
trivialized by $(\Hd_K)'[\log(T)]$. Let $\widetilde{Y}\in
GL_n(\widetilde{\Hd_{K}}[\log(T)])$ be a fundamental matrix
solution of the differential equation $\M$. Then there exists a
finite extension $K'/K$ such that the matrix
\begin{equation}
\widetilde{A}(q,T):=\sigma_q(\widetilde{Y})\cdot
\widetilde{Y}^{-1}
\end{equation}
belongs to $GL_n(\Hd_{K'})$, for all $q\in\mathrm{D}^-_{K'}(1,1)$.
In particular the operator $\sigma_q$ acting on
$\widetilde{\Ed_K}$ stabilizes both $\widetilde{\Hd_K}$, and
$\widetilde{\C^{\textrm{ét}}_K}$, and hence $\M$ is
$\widetilde{\C_K^{\textrm{ét}}}[\log(T)]$-extensible to the whole
disc $\mathrm{D}^-(1,1)$ (cf. Def. \ref{definition of strongly
confluent at q}).
\end{lemma}
\begin{proof} We can suppose that $K=K'$.
By Corollary \ref{N otimes U_m}, and by canonical extension (cf.
definition \ref{CAN-+---}), one can assume that $\M=N$, or
$\M=U_m$, where $\N$ is trivialized by an Galois étale extension
$(\Hd_K)'$ of $\Hd_K$, and where $U_m$ is defined over $\Hd_K$ as
in Corollary \ref{N otimes U_m}. The case ``$\M=U_m$'' is trivial,
since both the matrices of $\delta_1^{U_m}$ and of
$\sigma_q^{U_m}$ can be described explicitly as in section
\ref{unipotent in detail}. Let now $\M=\N$ (i.e., $\M$ is
trivialized by $\widetilde{\Hd_K}$). In this case the solution
matrix $\widetilde{Y}$ lies in $GL_n((\Hd_K)')$. The special
extension $(\Hd_K)'/\Hd_K$ corresponds via the equivalence of
section \ref{Katz-Matsuda equivalence section} to a finite Galois
extension $\F/k(\!(t)\!)$. Let $\G:=\mathrm{Gal}(\F/k(\!(t)\!))$,
then $\G$ acts on $(\Hd_K)'$ by $\Hd_K$-automorphisms, and
moreover the fixed points under this action are exactly the
elements of $\Hd_K$ (cf. Lemma \ref{(S^dag(F))^G=Hd_K}). After
enlarging $K$, if necessary, for all $\gamma\in\G$, one has
\begin{equation}
\gamma(\widetilde{Y})=\widetilde{Y}\cdot
H_\gamma\;,\quad\textrm{with}\quad H_\gamma\in GL_n(K)\;.
\end{equation}
Indeed by Lemma \ref{Ed' generated by solutions} the corresponding
Galois extension $(\Ed_K)'/\Ed_K$ is generated by the entries of
$\widetilde{Y}$. Hence $(\Ed_K)'/\Ed_K$ is a Picard-Vessiot
extension of $\Ed_K$ with \emph{differential} Galois group $\G$.
It follows then by Picard-Vessiot theory that $H_\gamma \in
GL_n(K)$ (cf. \cite[Obs.1.26]{VS}). Since $\sigma_q$ commutes with
every $\gamma\in\G$ (cf. Section \ref{unique extension of sigma to
Ed'}), one finds
\begin{equation}
\gamma(\widetilde{A}(q,T)) = \gamma(\sigma_q(\widetilde{Y})\cdot
\widetilde{Y}^{-1}) = \sigma_q(\widetilde{Y})\cdot H_\gamma
\cdot(\widetilde{Y}\cdot H_\gamma)^{-1}=\widetilde{A}(q,T)\;.
\end{equation}
Hence $\widetilde{A}(q,T)$ belongs to $\Hd_K$, for all $|q-1|<1$.
\end{proof}

\begin{lemma}\label{embedding in a_K(1,1)}
Let
$\a_{K^{\mathrm{alg}}}(1,1):=\bigcup_{K'/K=\textrm{finite}}\a_{K'}(1,1)$.
There exists an embedding
$\widetilde{\C^{\textrm{ét}}_K}[\log(T)]\subseteq
\a_{K^{\mathrm{alg}}}(1,1)$ commuting with the actions of
$\delta_1$, of $\varphi$, and  of $\sigma_q$, for all
$q\in\mathrm{D}^-_{K^{\mathrm{alg}}}(1,1)$. In other words
$\a_{K^{\mathrm{alg}}}(1,1)$ is a
$\widetilde{\C^{\textrm{ét}}_K}[\log(T)]$-$(\sigma,\delta)$-algebra
over the disc $\mathrm{D}^-_{K^{\mathrm{alg}}}(1,1)$, and one has
the following diagram of discrete $\Hd_K-(\sigma,\delta)-$algebras
over $\mathrm{D}^-(1,1)$: 
\begin{equation}
\xymatrix{\a_{K^{\mathrm{alg}}}(1,1)&\ar@{}[l]|-{\supset}\widetilde{\C^{\textrm{ét}}_K}\ar@{}[r]|-{\subset}&
\widetilde{\Hd_K}\ar@{}[r]|-{\subset}&\widetilde{\Ed_K}\ar@{}[r]|-{\subset}&\widetilde{\R_K}\\
\a_K(1,1)\ar@{}[u]|-{\cup}&\ar@{}[l]|-{\supset}\Hd_K\ar@{}[u]|-{\cup}&\subset&\Ed_K\ar@{}[u]|-{\cup}\ar@{}[r]|-{\subset}&\R_K\ar@{}[u]|-{\cup}\;.}
\end{equation}
\end{lemma}
\begin{proof}
In the following we assume $K$ to be sufficiently large in order
that every Special objects appearing in the proof is quasi
unipotent. Let $\M\in\delta_1-\Mod(\Hd_K)^{\mathrm{Sp}}$, be a
Special differential equation trivialized by $\widetilde{\Hd_K}$.
Let $(\C^{\textrm{ét}}_K)'$ be the smallest sub-$\Hd_K$-algebra of 
$\widetilde{\Hd_K}$ trivializing $\M$. By definition
$(\C^{\textrm{ét}}_K)'$ is generated over $\Hd_K$ by the entries
$\{\widetilde{y}_{i,j}\}_{i,j}$ of a matrix solution
$\widetilde{Y}$ of $\M$ in $\widetilde{\Hd_K}$. We consider
$\Hd_K[\log(T)]$, $(\C^{\textrm{ét}}_K)'[\log(T)]$,
$\widetilde{\C^{\textrm{ét}}_K}[\log(T)]$, $\a_K(1,1)$ as
differential algebras (we forget the actions of $\sigma_q$ in a
first time). We have an embedding $\Hd_K[\log(T)]\subset\a_K(1,1)$
commuting with $\delta_1$ sending the symbol $\log(T)$ into the
power series $\sum_{n\geq 1}(-1)^{n-1} (T-1)^n/n\in\a_K(1,1)$. We
extends this embedding to $(\C_K^{\textrm{ét}})'[\log(T)]$ as
follows. Since the differential equation $\M$ has its coefficients
in $\Hd_K$ we can consider its Taylor solutions $Y(T,1)$ at the
point $1$. Since $\M$ is solvable, then $Y(T,1)\in
GL_n(\a_K(1,1))$. Let now $\mathcal{F}_K:=\mathrm{Frac}(\Hd_K)$ be
the field of fractions of $\Hd_K$. Since $\mathcal{F}_K$ is a
field, then (up to enlarge $K$) we can apply the Picard-Vessiot
theory to obtain an isomorphism
$\mathcal{F}_K[\{\widetilde{y}_{i,j}\}_{i,j}]
\xrightarrow[]{\;\;\sim\;\;}
\mathcal{F}_K[\{y_{i,j}(T,1)\}_{i,j}]$, sending
$\widetilde{y}_{i,j}$ into $y_{i,j}(T,1)$, and commuting with
$\delta_1$. Clearly this isomorphism identifies
$(\C^{\textrm{ét}}_K)'=\Hd_K[\{\widetilde{y}_{i,j}\}_{i,j}]$ with
$\Hd_K[\{y_{i,j}(T,1)\}_{i,j}]\subset\a_K(1,1)$. If $\M'$ is
another differential equation, and if
$\Hd_K[\{\widetilde{y}'_{i,j}\}_{i,j}]$ is the corresponding
Picard-Vessiot extension identified with
$\Hd_K[\{y_{i,j}'(T,1)\}_{i,j}]\subset\a_K(1,1)$, then the
embedding corresponding to $\M\oplus\M'$ extends these two
embedding since the entries of a solutions of $\M\oplus\M'$ are
the families
$\{\widetilde{y}_{i,j},\widetilde{y}'_{h,k}\}_{i,j,h,k}$ and
$\{y_{i,j}(T,1),y'_{h,k}(T,1)\}_{i,j,h,k}$ respectively. It is
hence clear that this family of embedding are compatible, so that
we obtain an embedding
$\widetilde{\C^{\textrm{ét}}_K}\subseteq\a_K(1,1)$ commuting with
$\delta_1$, and consequently
$\widetilde{\C^{\textrm{ét}}_K}[\log(T)]\subseteq\a_K(1,1)$ also
commutes with $\delta_1$. Notice that $\log(T)$ is algebraically
free over $\Hd_K$ and hence over $\widetilde{\C^{\textrm{ét}}_K}$
which is union of finite algebras over $\Hd_K$ (cf. Lemma
\ref{(S^dag(F))^G=Hd_K}). We can now check that this embedding
commutes with $\sigma_q$ (resp. $\varphi$), by looking to its
action on the entries $\{\widetilde{y}_{i,j}\}_{i,j}$, and
$\{y_{i,j}(T,1)\}_{i,j}$. Hence it is enough to prove that the
isomorphism $\Ed_K[\{\widetilde{y}_{i,j}\}_{i,j}]
\xrightarrow[]{\;\;\sim\;\;} \Ed_K[\{y_{i,j}(T,1)\}_{i,j}]$
commutes with $\sigma_q$ (resp. $\varphi$). Observe that, if we
fix an embedding of $\mathcal{F}_K[\{y_{i,j}(T,1)\}_{i,j}]$ in an
fixed algebraically closure of $\Ed_K$, then
$\Ed_K[\{y_{i,j}(T,1)\}_{i,j}]$ is, by definition, the smallest
field containing $\Ed_K$ and $\{y_{i,j}(T,1)\}_{i,j}$. The action
of $\delta_1,\sigma_q,\varphi$ are defined on
$\Ed_K[\{y_{i,j}(T,1)\}_{i,j}]$ as the extensions of
$\delta_1,\sigma_q,\varphi$ on
$\mathcal{F}_K[\{y_{i,j}(T,1)\}_{i,j}]$. This extension exists
since $\Hd_K[\{y_{i,j}(T,1)\}_{i,j}]\cap\Ed_K=
\Hd_K[\{\widetilde{y}_{i,j}\}_{i,j}]\cap\Ed_K= \Hd_K$ (cf. Lemma
\ref{(S^dag(F))^G=Hd_K}), and since $\delta_1,\sigma_q,\varphi$
act on $Y(T,1)$ by multiplication by matrices with coefficients in
$\Hd_K$, (cf. \cite[Sect.6,no.1,Prop.1]{Bou-Alg-4-5}). By Lemma
\ref{unique extension of sigma_q to Ed[y_i,j]}, there exists a
unique extension of $\sigma_q$ to
$\Ed_K[\{\widetilde{y}_{i,j}\}_{i,j}]$, and of course a unique
extension of $\varphi$ since
$\Ed_K[\{\widetilde{y}_{i,j}\}_{i,j}]/\Ed_K$ is unramified. Hence
the isomorphism $\Ed_K[\{\widetilde{y}_{i,j}\}_{i,j}]
\xrightarrow[]{\sim} \Ed_K[\{y_{i,j}(T,1)\}_{i,j}]$ commutes with
$\sigma_q$ and  $\varphi$.
\end{proof}

\begin{remark}
The same statement holds for $\a_K(c,1)$ instead of $\a_K(1,1)$,
providing that $|c|=1$, $c\in K$, and $\varphi(c)=c$.
\end{remark}

\begin{proof}[Proof of Theorem \ref{p-adic local monodromy theorem (generalized form)}]
By Proposition \ref{solvable extend to all the disk}, one has
$(\sigma,\delta)-\Mod(\R_K)^{(\phi)}_S =
(\sigma,\delta)-\Mod(\R_K)^{(\phi)}_{\mathrm{D}^-(1,1)}$, (resp.
$\sigma-\Mod(\R_K)^{(\phi)}_S =
\sigma-\Mod(\R_K)^{(\phi)}_{\mathrm{D}^-(1,1)}$). On the other
hand, $(\sigma,\delta)-\Mod(\R_K)^{(\phi)}_{\mathrm{D}^-(1,1)}=
\sigma-\Mod(\R_K)^{(\phi)}_{\mathrm{D}^-(1,1)}$ (cf.
\eqref{sigma-an=(sigma,delta)-an}). Moreover, if
$q\in\mathrm{D}^-(1,1)-\bs{\mu}_{p^{\infty}}$, then
$(\sigma,\delta)-\Mod(\R_K)^{(\phi)}_{\mathrm{D}^-(1,1)}=(\sigma_q,\delta_q)-\Mod(\R_K)^{(\phi)}$
(resp. $\sigma-\Mod(\R_K)^{(\phi)}_{\mathrm{D}^-(1,1)} =
\sigma_q-\Mod(\R_K)^{(\phi)}$). Hence, without loss of generality,
we can assume that $\M$ is a Taylor admissible
$(\sigma_q,\delta_q)-$module, with Frobenius structure. 
\if{By
Remark \ref{Remark explaining the Propagation}, Section \ref{C_1
subset C_2 ==> Then Def^C_1=Def^C_2}, and Lemmas \ref{la
chiave},\ref{embedding in a_K(1,1)}, }\fi
The proof follows now by the discussion after Theorem \ref{p-adic
local monodromy theorem (generalized form)}.
\end{proof}

\begin{corollary}\label{Conf Tay = Conf a_K(1,1) = Conf widetilde R_K}
Let $S\subset\mathrm{D}^-(1,1)-\bs{\mu}_{p^{\infty}}$ (resp.
$S\subset\mathrm{D}^-(1,1)$). We have the following equalities
\begin{equation}
\Conf^{\mathrm{Tay}}_q\;\;\stackrel{\textrm{Rem.} \ref{Remark
explaining the
Propagation}}{=}\;\;\Conf^{\a_K(1,1)}_q\;\;=\;\;\Conf_q^{\a_{K^{\mathrm{alg}}}(1,1)}\;\;\stackrel{\textrm{Lemma
} \ref{embedding in
a_K(1,1)}}{=}\;\;\Conf^{\widetilde{\C^{\textrm{ét}}_K}[\log(T)]}_q\;,
\end{equation}
where the first three equalities holds for these functors on
$\sigma-\Mod(\Hd_K)^{[1]}_{S}$ (resp.
$(\sigma,\delta)-\Mod(\Hd_K)^{[1]}_{S}$), while, in the last
equality, one considers the restrictions of these functors to the
full subcategory of $\sigma-\Mod(\Hd_K)^{[1]}_{S}$ (resp.
$(\sigma,\delta)-\Mod(\Hd_K)^{[1]}_{S}$) of objects trivialized by
$\widetilde{\C^{\textrm{ét}}_{K}}[\log(T)]$. In particular the
last equality holds on $\sigma-\Mod(\Hd_K)_S^{\mathrm{Sp}}$ (resp.
$(\sigma,\delta)-\Mod(\Hd_K)_S^{\mathrm{Sp}}$). The same relation
holds for deformation functors.
\end{corollary}
\begin{proof}
By Remark \ref{Remark explaining the Propagation} the restriction
of $\Conf^{\mathrm{Tay}}_q$ to the category of \emph{solvable}
objects coincides with $\Conf^{\a_K(1,1)}_q$. A \emph{solvable}
object over $\Hd_K$ is trivialized by $\a_{K^{\mathrm{alg}}}(1,1)$
if and only if it is trivialized by $\a_K(1,1)$. Indeed both these
conditions are verified if and only if its Taylor solution at $1$
converges on $\mathrm{D}^-(1,1)$. Hence $\Conf^{\a_K(1,1)}_q =
\Conf_q^{\a_{K^{\mathrm{alg}}}(1,1)}$ on solvable objects. Now, by
Theorem \ref{p-adic local monodromy theorem (generalized form)},
Special objects are trivialized by
$\widetilde{\C^{\textrm{ét}}_K}[\log(T)]$ hence by
$\a_{K^{\mathrm{alg}}}(1,1)$ (cf. Lemma \ref{embedding in
a_K(1,1)}).
\end{proof}

\subsection{The confluence of Andr\'e-Di Vizio}\label{Andre e Di Vizio}
In this last section
we prove that the restriction of $\Conf_q^{\mathrm{Tay}}$ to
$\sigma_q-\Mod(\R_K)^{(\phi)}$ is isomorphic to the functor
``$\Conf$'' defined in \cite[Section 15.1]{An-DV}. In all this
last section $q\in\mathrm{D}^-(1,1)-\bs{\mu}_{p^{\infty}}$.

We recall that an \emph{antecedent} of a $\sigma_q-$module $\M$
over $\R_K$ is a $\sigma_{q^p}-$module $\M_{1}$ such that
$\phi^*(\M_1)$ is isomorphic to $\M$ as $\sigma_q$-module. The
antecedent is unique up to isomorphisms, because this fact is true
for differential equations (cf. Remark \ref{Remark on the
definition of the pull-back-***}). In order to preserve the
notations of \cite{An-DV}, we fix a $s\geq 1$, and we call $\M_1$
the $s$-th. antecedent of $\M$, i.e.
$\Phi:(\phi^*)^{s}(\M_{1})\xrightarrow[]{\;\sim\;}\M$.

The following definition was given in \cite{An-DV} under the
assumption $|q-1|<|p|^{1/(p-1)}$. The same definition holds for
$q\in\mathrm{D}^-(1,1)-\bs{\mu}_{p^{\infty}}$.

\begin{definition}[(\protect{\cite[12.11]{An-DV}})]\label{confluent weak Frob Str}
Let $s\in\mathbb{N}_{>0}$. A \emph{Confluent Weak Frobenius
Structure} (CWFS) on a $\sigma_q$-module
$\M_0:=(\M_0,\sigma_q^{\M_0})\in\sigma_q-\Mod(\R_K)$ is a sequence
$\{\sigma_{q^{p^{sm}}}^{\M_{m}}\}_{m\geq 0}$ of
$q^{p^{sm}}$-difference operators on $\M_0$, together with a
family of isomorphisms
\begin{equation}\label{isomorphismsssdes}
\Phi_m:(\;(\phi^{*})^{s}(\M_{0})\;,\;(\phi^*)^s(\sigma_{q^{p^{s(m+1)}}}^{\M_{m+1}})\;)
\xrightarrow[]{\;\;\sim\;\;}
(\M_{0},\sigma_{q^{p^{sm}}}^{\M_m})\;,
\end{equation}
of $q^{p^{sm}}-$difference modules (identifying
$(M_0,\sigma_{q^{p^{sm}}}^{\M_m})$ to the $s$-th antecedent of
$(\M_0,\sigma_q^{\M_0})$), such that:
\begin{enumerate}
\item The operators
$\Delta_{q^{p^{sm}}}^{\M_{m}}:=(\sigma_{q^{p^{sm}}}^{\M_{m}}-\mathrm{Id}^{\M_0})/(q^{p^{sm}}-1)$
converge to a derivation $\Delta^{\M_\infty}$ on $\M_0$.  %
\item  If $\M_\infty:=(\M_0,\Delta^{\M_\infty})$ is this
differential module, then the sequence of isomorphisms
\eqref{isomorphismsssdes} converges to a Frobenius isomorphism
$\Phi_\infty:\phi^*(\M_\infty)\xrightarrow[]{\;\;\sim\;\;}\M_\infty$. %
\end{enumerate}
We denote by
\begin{equation}\label{sigma_q-Mod(R_K)^(conf(phi))}
\sigma_q-\Mod(\R_K)^{\mathrm{conf}(\phi)}
\end{equation}
 the category whose objects
are families of operators
$(\M_0,\{\sigma_{q^{p^{sm}}}^{\M_m}\}_{m\geq 0})$ on $\M_0$
admitting the existence of a family $\{\Phi_m\}_{m\geq 0}$ making
it on a confluent weak Frobenius structure on
$(\M_0,\sigma_q^{\M_0})$. A morphism
$\alpha:(\M_0,\{\sigma_{q^{p^m}}^{\M_m}\}_{m\geq
0})\longrightarrow (\N_0,\{\sigma_{q^{p^m}}^{\N_m}\}_{m\geq 0})$
is a $\R_K$-linear morphisms $\alpha:\M_0\to\N_0$ verifying
simultaneously $\alpha\circ\sigma_{q^{p^{sm}}}^{\M_m}=
\sigma_{q^{p^{sm}}}^{\N_m}\circ\alpha$, for all $m\geq 0$.
\end{definition}

\subsubsection{Construction of CWFSs.}\label{infinitely many WCFS} A
$q$-difference module $(\M_0,\sigma_q^{\M_0})$ admits infinitely
many Confluent Weak Frobenius Structures (CWFS), even if
$(\M_0,\sigma_q^{\M_0})$ admits a (strong) Frobenius structure.
Indeed if a CWFS $(\M_0,\{\sigma_{q^{p^{sm}}}^{\M_m}\}_{m\geq
0},\{\Phi_m\}_{m\geq 0})$ on $(\M_0,\sigma_q^{\M_0})$ is given, we
give now an algorithm to produce infinitely many CWFS on
$(\M_0,\sigma_q^{\M_0})$. Let $\{\psi_m:\M_0\simto\M_0\}_{m\geq
0}$ be a sequence of $\R_K$-linear automorphisms of $\M_0$ such
that $\lim_m\psi_m=\mathrm{Id}^{\M_0}$. Define
\begin{equation}
\sigma_{q^{p^{sm}}}^{\M_m'}:=\psi_m\circ\sigma_{q^{p^{sm}}}^{\M_m}\circ\psi_m^{-1}\;,\qquad
\Phi_m':=\psi_m\circ\Phi_m\circ[\phi^*(\psi_{m+1})]^{-1}\;.
\end{equation}
One easily checks that
$(\M_0,\{\sigma_{q^{p^{sm}}}^{\M_m'}\}_m,\{\Phi_m'\}_m)$ is again
a CWFS on $(\M_0,\sigma_{q}^{\M_0})$. Notice that this new CWFS is
not always isomorphic to the first one (even if
$\psi_0=\mathrm{Id}^{\mathrm{M_0}}$). Indeed, by definition, an
isomorphism is a single arrow $\alpha:\M_0\to\M_0$ satisfying
\emph{simultaneously} $\alpha\circ\sigma_{q^{p^{sm}}}^{\M_m}=
\sigma_{q^{p^{sm}}}^{\M_m'}\circ\alpha$, for all $m\geq 0$.
Nevertheless, since $\lim_{m}\psi_m=\mathrm{Id}^{\M_0}$, the limit
differential equation is the same for all CWFS defined in this way
(cf. Remark \ref{Delta_infty depends on all the family M_m}). We
observe moreover that $\psi_m$ defines an isomorphism of
$q^{p^{sm}}$-difference modules between
$(\M_0,\sigma_{q^{p^{sm}}}^{\M_m})$ and
$(\M_0,\sigma_{q^{p^{sm}}}^{\M_m'})$, this agrees with the
uniqueness of the antecedent by Frobenius. If $\widetilde{Y}_m$ is
the solution of $(\M_0,\sigma_{q^{p^{sm}}}^{\M_m})$ in
$GL_n(\widetilde{\R_K}[\log(T)])$, and if $B_m(T)\in GL_n(\R_K)$
is the matrix of $\psi_m$, then the solution of
$(\M_0,\sigma_{q^{p^{sm}}}^{\M_m'})$ is given by
$B_m(T)\widetilde{Y}_m$.

\begin{remark} Assume that $\M_0$ admits a (strong) Frobenius structure.
The constancy of the solution does not follows from the preview
definition. Indeed a solution of $(\M_0,\sigma_q^{\M_0})$, with
values in $\C$ is a morphism $\alpha:\M_0\to\C$ satisfying
$\alpha\circ\sigma_q^{\M_0}=\sigma^{\C}_q\circ\alpha$ (cf. section
\ref{section - solution formal def}). The fact that $\alpha$ is a
solution of $(\M_0,\sigma_q^{\M_0})$ does not implies that
$\alpha$ commutes also with $\sigma_{q^{p^{sm}}}^{\M_m}$. Indeed
the data $(\M_0,\{\sigma_{q^{p^{sm}}}^{\M_m}\}_{m\geq 0})$ is
\emph{not necessarily} a discrete $\sigma$-module over
$S=\{q^{p^{sm}}\}_{m\geq 0}$, because the map
$S\to\mathrm{Aut}^{\mathrm{cont}}(\M_0)$ sending $q$ into
$\sigma_{q^{p^{sm}}}^{\M_m}$ is not supposed to have any coherency
(cf. Remark \ref{first remark},(iii)). To obtain the constancy of
the solutions we need to ``\emph{rigidify}'' these constructions
by introducing the notion of \emph{$\C$-constant} $\sigma$-module
(cf. Remarks \ref{intro -constancy of solutions}, \ref{Morphisms
between sigma modules}, and Example \ref{Morphisms between sigma
modules -examples}).
\end{remark}

\subsubsection{} \label{chi(M)} We have an evident
fully faithful functor
\begin{equation}
\chi_q^{(\phi)}\;:\;\sigma_q-\Mod(\R_K)^{(\phi)}\;\;\longrightarrow\;\;\sigma_q-\Mod(\R_K)^{\mathrm{conf(\phi)}}
\end{equation}
defined by
\begin{equation}
\chi_q^{(\phi)}(\M_0,\sigma_q^{\M_0})\;:=\;(\M_0,\{(\sigma_q^{\M_0})^{p^{sm}}\}_{m\geq
0})\;,
\end{equation}
where $s$ is sufficiently large to have an isomorphism
$\Phi:(\phi^*)^s(\M_0,\sigma_q^{\M_0})\xrightarrow[]{\;\sim\;}(\M_0,\sigma_q^{\M_0})$,
and $\Phi_m:=\Phi$, for all $m\geq 0$. On the other hand we have
another functor (cf. \cite[Sect.12.3]{An-DV})
\begin{equation}
Lim_\infty^{(\phi)}\;:\;\sigma_q-\Mod(\R_K)^{\mathrm{conf(\phi)}}\;\;\longrightarrow\;\;\delta_1-\Mod(\R_K)^{(\phi)}
\end{equation}
sending $(\M_0,\{\sigma_{q^{p^{sm}}}^{\M_m}\}_{m\geq
0},\{\Phi_m\}_{m\geq 0})$ into its limit differential equation
$(\M_0,\Delta^{\M_\infty})$. We have actually
\begin{equation}\label{frtvfbdnx,}
Lim_\infty^{(\phi)} \circ
\chi_q^{(\phi)}\;=\;\Conf^{\mathrm{Tay}}_q\;.
\end{equation}
Indeed if $(\M_0,\sigma_q^{\M_0})$ has a (strong) Frobenius
Structure, then $(\M_0,\{(\sigma_q^{\M_0})^{p^{sm}}\}_{m\geq 0})$
is a (solvable) Taylor admissible $\sigma$-module over
$S:=\{q^{p^{sm}}\}_{m\geq 0}$ (cf. Def. \ref{Taylor adm over R_K
and Hd_K .....}). Hence, by Section \ref{explicit computation of
the derivation as limit}, the differential equation
$\Conf_q^{\mathrm{Tay}}(\M_0,\sigma_q^{\M_0})$ is given by the
limit
$\Delta^{\M_\infty}:=\lim_{m\to\infty}\Delta_{q^{p^{sm}}}^{\M_{m}}$
of definition \ref{confluent weak Frob Str}. Moreover since the
operator $\sigma_{q^{p^{sm}}}^{\M_m}$ is determined by the
knowledge of the solutions of $(\M_0,\sigma_{q^{p^{sm}}}^{\M_m})$
in $GL_n(\widetilde{\R_K}[\log(T)])$, then
$\chi^{(\phi)}_q(\M_0,\sigma_q^{\M_0})$ is the \emph{unique} CWFS
on $(\M_0,\sigma_q^{\M_0})$ such that the fundamental matrix
solution of $(\M_0,\sigma_q^{\M_0})$ in
$GL_n(\widetilde{\R_K}[\log(T)])$ (or equivalently its Taylor
solution in $\a_K(1,1)$, cf. Lemma \ref{embedding in a_K(1,1)}) is
simultaneously solution of every
$(\M_0,\sigma_{q^{p^{sm}}}^{\M_m})$.

\begin{remark}\label{Delta_infty depends on all the family M_m}
It is not clear whether the limit differential equation
$(\M_0,\Delta^{\M_\infty})$ depends on the particular CWFS on
$(\M_0,\sigma_q^{\M_0})$ or, analogously, if there exists two non
isomorphic $q$-difference modules endowed with CWFS giving the
same limit differential equation. Indeed both these phenomena
arise in the category $\sigma_q-\Mod(\R_K)^{\mathrm{conf}}$
defined below.
\end{remark}

\begin{lemma}\label{chi = D^conf circ V^(phi)}
If $K$ is algebraically closed, then the functor $\chi^{(\phi)}_q$
is isomorphic to the functor
$D_{\sigma_q}^{\mathrm{conf}(\phi)}\circ V_{\sigma_q}^{(\phi)}$ of
\cite[Cor. 14.8]{An-DV}. Hence the functor
$\Conf^{\mathrm{Tay}}_q(\stackrel{\eqref{frtvfbdnx,}}{=}Lim_\infty^{(\phi)}\circ\chi_q^{(\phi)})$
is isomorphic to the confluence functor
$\Conf:=Lim_\infty^{(\phi)}\circ
D_{\sigma_q}^{\mathrm{conf}(\phi)}\circ V_{\sigma_q}^{(\phi)}$ as
it was defined in \cite[Section 15.1]{An-DV}.
\end{lemma}
\begin{proof}
As explained in the introduction,
$V_{\sigma_q}^{(\phi)}(\M,\sigma_q^{\M})$ (resp.
$V_{d}^{(\phi)}(\M,\delta_1^{\M})$) is the (dual of the) space of
solutions of $(\M,\sigma_q^{\M})$ (resp. $(\M,\delta_1^{\M})$) in
$\widetilde{\R_K}[\log(T)]$. By definition $D_d^{(\phi)}\circ
V_d^{(\phi)}\cong\mathrm{Id}$, and $D_{\sigma_q}^{(\phi)}\circ
V_{\sigma_q}^{(\phi)}\cong\mathrm{Id}$. Then $D_{d}^{(\phi)}\circ
V_{\sigma_q}^{(\phi)}=\Conf_q^{\widetilde{\R_K}[\log(T)]}$ is the
functor sending $(\M,\sigma_q^{\M})$ into the differential
equation having the same solutions in $\widetilde{\R_K}[\log(T)]$.
By definition (cf. \cite[Prop. 12.17]{An-DV}) one has
$D_{\sigma_q}^{\mathrm{conf}(\phi)} \cong \chi^{(\phi)}_q \circ
D_{\sigma_q}^{(\phi)}$. This proves that the functor
$\mathrm{Conf}$ of \cite{An-DV} is equal to
$\Conf_q^{\widetilde{\R_K}[\log(T)]}$. By Corollary \ref{Conf Tay
= Conf a_K(1,1) = Conf widetilde R_K} we conclude.
\end{proof}

\subsubsection{}\label{frfrfrfr,,,,frfr,fr} Lemma \ref{chi = D^conf
circ V^(phi)} clarifies the nature of the functor $\Conf$ of
\cite{An-DV} (cf. Corollary \ref{Conf Tay = Conf a_K(1,1) = Conf
widetilde R_K}). Indeed $\Conf$ is equal to
$\Conf_q^{\mathrm{Tay}}$, and sends a $q$-difference equation into
the differential equation having the same Taylor solutions (or
equivalently having the same ``étale'' solutions in
$\widetilde{\R_K}[\log(T)]$, cf. Lemma \ref{embedding in a_K(1,1)}
and Corollary \ref{Conf Tay = Conf a_K(1,1) = Conf widetilde
R_K}). This functor actually does not depend on the existence of a
Frobenius Structure and exists in the more general context of
\emph{admissible} modules. This generalizes the constructions of
\cite{An-DV} to all $q\in\mathrm{D}^-(1,1)-\bs{\mu}_{p^{\infty}}$,
removing also the assumption ``$K=K^{\mathrm{alg}}$''. Notice that
the equivalence provided by the Propagation Theorem requires only
the definition and the formal properties of the Taylor solution
$Y(x,y)$. For this reason the equivalences
$\Conf^{\mathrm{Tay}}_q$ and $\mathrm{Def}_{q,q'}^{\mathrm{Tay}}$
are not a consequence of the heretofore developed theory.
Conversely our Confluence implies the main results of \cite{An-DV}
and also of \cite{DV-Dwork}.

\subsubsection{A conjecture of \cite{An-DV}.}
Section \ref{infinitely many WCFS} proves that the fully faithful
functor $\chi_q^{(\phi)}$ is not an equivalence. This answer to a
question asked in \cite[Corollary 14.8, and after]{An-DV}.
Nevertheless observe that the existence of a CWFS on
$(\M_0,\sigma_q^{\M_0})$ is equivalent to the existence of a
\emph{strong} Frobenius structure on it. This have been firstly
proved for rank one equations (cf. \cite[Prop.7.3]{An-DV}, the
case with rational coefficient follows actually from section
\ref{classific Rk1}, indeed every rank one equation
with rational exponent has a (strong) Frobenius Structure). %
The general case is proved as follows.

\begin{definition}
We define
\begin{equation}\label{sigma_q-Mod(R_K)^conf}
\sigma_q-\Mod(\R_K)^{\mathrm{conf}}
\end{equation}
as the category whose objects are $\R_K$-modules $\M$ together
with a family of $\sigma_q$-semi-linear automorphisms
$\{\sigma_{q^{p^{sm}}}^{\M}:\M\simto\M\}_{m\geq 0}$ (without any
condition of compatibility) such that the limit
\begin{equation}
\delta_1^{\M}:=\lim_{m\to\infty}\frac{\sigma_{q^{p^{sm}}}^{\M}-\mathrm{Id}}{q^{p^{sm}}-1}
\end{equation}
converges to a connection $\delta_1^{\M}$ on $\M$. Morphisms
between $(\M,\{\sigma_{q^{p^{sm}}}^{\M}\}_{m\geq 0})$ and
$(\N,\{\sigma_{q^{p^{sm}}}^{\N}\}_{m\geq 0})$ are $\R_K$-linear
morphisms $\alpha:\M\to\N$ satisfying simultaneously $\alpha\circ
\sigma_{q^{p^{sm}}}^{\M} = \sigma_{q^{p^{sm}}}^{\N}\circ\alpha$,
for all $m\geq 0$.
\end{definition}

\begin{remark}
We have a functor
\begin{equation}
Lim_\infty\;:\;\sigma_q-\Mod(\R_K)\;\;
\longrightarrow\;\; 
\delta_1-\Mod(\R_K)
\end{equation}
sending $(\M,\{\sigma_{q^{p^{sm}}}^{\M}\}_{m\geq
0})$ into its limit differential equation. Indeed if
$\alpha:\M\to\N$ satisfies simultaneously $\alpha\circ
\sigma_{q^{p^{sm}}}^{\M} = \sigma_{q^{p^{sm}}}^{\N}\circ\alpha$,
for all $m\geq 0$, then, by passing to the limit, one has
$\alpha\circ\delta_1^{\M}=\delta_1^{\N}\circ\alpha$. We have then
the following commutative diagram of categories:
\begin{equation}
\xymatrix{
\sigma_q-\Mod(\R_K)^{[r]}\ar[r]^{\chi_q}\ar@{}[dr]|-{\odot}&\ar@{}[dr]|-{\odot}\sigma_q-\Mod(\R_K)^{\mathrm{conf}}\ar[r]^-{Lim_\infty}&\delta_1-\Mod(\R_K)\\
\sigma_q-\Mod(\R_K)^{(\phi)}\ar[r]_{\chi^{(\phi)}_q}\ar@{}[u]|{\cup}&\sigma_q-\Mod(\R_K)^{\mathrm{conf}(\phi)}
\ar[r]_-{Lim_\infty^{(\phi)}}\ar@{}[u]|{\cup}&\delta_1-\Mod(\R_K)^{(\phi)}\ar@{}[u]|{\cup},}
\end{equation}
where $r\geq |q-1|$, and where $\chi_q$ sends $(\M,\sigma_q^{\M})$
into $(\M,\{(\sigma_q^{\M})^{p^{sm}}\}_{m\geq 0})$. By Section
\ref{explicit computation of the derivation as limit}, as above
\begin{equation}
Lim_\infty\circ\chi_q=\Conf^{\mathrm{Tay}}_q\;:\;\sigma_q-\Mod(\R_K)^{[r]}\;
\xrightarrow[]{\;\;\sim\;\;}\;\delta_1-\Mod(\R_K)^{[r]}\;\subset\;\delta_1-\Mod(\R_K)\;.
\end{equation}
\end{remark}

\begin{corollary}
Let $q\in\mathrm{D}^-(1,1)-\bs{\mu}_{p^{\infty}}$. Let
$(\M,\sigma_q^{\M})\in\sigma_q-\Mod(\R_K)^{[r]}$, with $r\geq
|q-1|$. Then $(\M,\sigma_q^{\M})$ admits a CWFS if and only if it
admits a (strong) Frobenius structure.
\end{corollary}
\begin{proof}
Assume that $Lim_{\infty}\circ\chi_q(\M,\sigma_q^{\M})$ lies in
$\delta_1-\Mod(\R_K)^{(\phi)}$. By \eqref{Y(T^p) = theta(T) Y(T)},
$\mathrm{Def}^{\mathrm{Tay}}_{1,q}\circ
Lim_{\infty}\circ\chi_q(\M,\sigma_q^{\M})$ lies in
$\sigma_q-\Mod(\R_K)^{(\phi)}$. Now since
$\mathrm{Def}^{\mathrm{Tay}}_{1,q} \circ Lim_{\infty} \circ
\chi_q=
 \mathrm{Def}^{\mathrm{Tay}}_{1,q} \circ \Conf_q^{\mathrm{Tay}}
= \mathrm{Id}$, then $\mathrm{Def}^{\mathrm{Tay}}_{1,q} \circ
Lim_{\infty} \circ \chi_q(\M,\sigma_q^{\M})$ is isomorphic to
$(\M,\sigma_q^{\M})$, and has hence (strong) Frobenius structure.
\end{proof}

\subsection{The theory of slopes} In a sequence of papers,
G.Christol and Z.Mebkhout developed a theory of slopes for
$p$-adic differential equations over the Robba ring. We summarize
the main properties in the following theorem.
\begin{theorem}[(cf. \cite{Astx})]\label{break decomposition}
Let $\M$ be a solvable differential module over $\R_K$. There
exists a unique decomposition of $\M$, called \emph{break
decomposition}
\begin{equation}
\M=\oplus_{x\in\mathbb{R}_{\geq 0}}\M(x)\;,
\end{equation}
satisfying the following properties. Let $t_\rho$ be a generic
point for the norm $|\cdot|_\rho$ (cf. \eqref{ZERTY}), then there
exists $\varepsilon
>0$ such that
\begin{enumerate}

\item For all $\rho<]1-\varepsilon,1[$, $\M(x)$ is (the biggest
submodule of $\M$) trivialized by
$\a_K(t_{\rho},\rho^{x+1})$, 

\item For all $\rho<]1-\varepsilon,1[$, and for all $y < x$,
$\M(x)$ has no solutions in $\a_K(t_{\rho},\rho^{y+1})$.
\end{enumerate}
The number $\mathrm{Irr}(\M):= \sum_{x\geq
0}x\cdot\textrm{rank}_{\R_K}(\M(x))$ is called \emph{$p$-adic
irregularity} of $\M$, and it lies in $\mathbb{N}$.
\end{theorem}

The fact that $\mathrm{Irr}(\M)$ is integer is known as the
\emph{Hasse-Arf property}. This theorem has an analogous in the
theory of representations of the Galois group of a local field:

\begin{proposition}[(cf. \protect{\cite{Katz}})]
Let $\mathcal{I}$, $\mathcal{P}$ be the inertia and the wild
inertia subgroups of
$\mathrm{G}:=\mathrm{Gal}(k(\!(t)\!)^{\mathrm{sep}}/k(\!(t)\!))$.
Denote by $\{\mathcal{I}^{(x)}\}_{x\geq 0}$ the ``upper numbering
filtration'' of $\mathcal{I}$. Let $\V$ be a
$\mathbb{Z}[1/p]$-representation of $\G$, such that $\mathcal{P}$
acts through a finite discrete quotient. Then $\V$ admits a
\emph{break decomposition} $\V=\oplus_{x\geq 0}\V(x)$ of
$\G$-submodules $\V(x)$ such that $\V(0)=\V^{\mathcal{P}}$, and
for all $x>0$:
\begin{enumerate}
\item $(\V(x))^{\mathcal{I}^{(x)}}=0$;

\item For all $y>x$, $(\V(x))^{\mathcal{I}^{(y)}}=\V(x)$.
\end{enumerate}
The number $\textrm{Swan}(\V):=\sum_{x\geq
0}\mathrm{rank}_{\mathbb{Z}[1/p]}\V(x)$ is called \emph{Swan
conductor} of $\V$, and it lies in $\mathbb{N}$.
\end{proposition}

For an very inspiring overview about this analogy we refer to
\cite{An-Ramis-2}.

Different authors (cf. \cite{Tsu-swan}, \cite{Matsuda-unipotent},
\cite{Crew-Can-ext}) proved that the equivalence functor
introduced by J.-M.Fontaine (cf. \cite{Fo}, \cite{Ts}),
associating to a finite representation of $\G$, a
$(\varphi,\nabla)$-module over $\Ed_K$ (and hence a differential
module over $\R_K$) preserves the break decompositions. \emph{The
Swan conductor of a representation equals the Irregularity of the
corresponding differential equation}. 

In \cite{An} the author state a family of axiomatic conditions in
a general Tannakian category in order to have a ``\emph{theory of
slopes}''. The previous two cases respect the formalism of
\cite{An}. 

In a second time he conjectured (cf. \cite[Conjecture
4.2]{An-Ramis-2}) that a similar theory of slopes, should exists
also for $\sigma_q-\Mod(\R_K)^{(\phi)}$ and ask if this ``new''
theory of slopes is compatible with that of Christol-Mebkhout on
$\delta_1-\Mod(\R_K)^{(\phi)}$ (via the confluence), and hence
with the ramification theory on
$\underline{\mathrm{Rep}}_{K^{\mathrm{alg}}}(\mathcal{I}_{k^{\mathrm{alg}}(\!(t)\!)}\times\mathbb{G}_a)$
(via the Fontaine's functor $T_1$ of the introduction). He
suggested to proceed in analogy with the theory of
Christol-Mebkhout (cf. \cite{Astx}), reproducing their proofs in
the context of $q$-difference equations in order to obtain a
statement analogous to theorem \ref{break decomposition}. Finally
he asked whether this ``new'' theory of slopes on
$\sigma_q-\Mod(\R_{K^{\mathrm{alg}}})^{(\phi)}$ is compatible or
not with the theory of slopes of Christol-Mebkhout in
$\delta_1-\Mod(\R_{K^{\mathrm{alg}}})^{(\phi)}$ via the
equivalence $\Conf$ that he obtained in \cite{An-DV} for
$|q-1|<|p|^{\frac{1}{p-1}}$.

Afterwards, at the end of 2005, he actually obtained such a theory
of slopes for $\sigma_q-\Mod(\R_{K^{\mathrm{alg}}})^{(\phi)}$,
with $|q-1|<|p|^{\frac{1}{p-1}}$, and established the two
corollaries below in this case. These verifications will be
included in a fort coming paper of Y.André. This part have been
exposed by Y.André at the \emph{$24$-th Nordic and $1$-st
Franco-Nordic Congress of Mathematicians} ($6$ to $9$ January
$2006$, Rejkyavik, Iceland).

The next corollaries prove the above conjecture in the more
general context of $\sigma$-modules. We prove it for all
$|q-1|<1$, without any assumptions about the Frobenius structure,
and without assuming $K=K^{\mathrm{alg}}$. The equivalence
established by Corollary \ref{propagation confluence and def over
Robba} gives in fact the following analogous of Theorem \ref{break
decomposition} for $\sigma-$modules and $(\sigma,\delta)$-modules.
Thank to Proposition \ref{solvable extend to all the disk},
without loss of generality, we can reduce this statement to the
case $S=\{q\}$:

\begin{corollary}\label{q-break decomposition}
Let $|q-1|<1$, $q\in K$, (resp. $q\notin\bs{\mu}_{p^{\infty}}$).
Let $\M\in(\sigma_q,\delta_q)-\Mod(\R_K)^{[1]}$ (resp.
$\M\in\sigma_q-\Mod(\R_K)^{[1]}$). Then $\M$ admits a \emph{break
decomposition} $\M=\oplus_{x\geq 0}\M(x)$, where $\M(x)$ is
characterized by the following properties (analogues to i) and ii)
of Theorem \ref{break decomposition}). There exists $\varepsilon
>0$ such that
\begin{enumerate}

\item For all $\rho<]1-\varepsilon,1[$, $\M(x)$ is (the biggest
submodule of $\M$) trivialized by
$\a_K(t_{\rho},\rho^{x+1})$, 

\item For all $\rho<]1-\varepsilon,1[$, and for all $y < x$,
$\M(x)$ has no solutions in $\a_K(t_{\rho},\rho^{y+1})$.
\end{enumerate}
This decomposition is compatible with the confluence i.e.
$\M(x)=\Def_{1,q}^{\mathrm{Tay}}(\Conf_q^{\mathrm{Tay}}(\M)(x))$.
In particular the irregularity $\mathrm{Irr}_{\sigma_q}(\M) :=
\sum_{x\geq 0} x \cdot \textrm{rank}_{\R_K}\M(x)$ is a natural
number.
\end{corollary}
\begin{proof}
The ``\emph{slopes}'' and the ``\emph{Irregularity}'' are defined,
by Christol and Mebkhout (cf. \cite{Ch-Me}), by means of the
generic \emph{radius of the Taylor solutions}. The $K$-linear
equivalences $\Conf^{\mathrm{Tay}}_q$ and
$\mathrm{Def}^{\mathrm{Tay}}_{1,q}$ preserve, by definition, the
generic Taylor solution. It follows immediately that the
$q$-difference equation inherits then, via the equivalence
$\Conf_q^{\mathrm{Tay}}$, the \emph{slopes} of the attached
differential equation, together with their formal properties
(break decomposition, Hasse-Arf property, ...).
\end{proof}

\begin{corollary}
With the notations of \cite{An-DV} and \cite{An-Ramis-2}, if
$K=K^{\mathrm{alg}}$ is algebraically closed, the functor
$D_{\sigma_q}^{(\phi)}:\sigma_q-\Mod(\R_{K^{\mathrm{alg}}})^{(\phi)}\longrightarrow
\underline{\mathrm{Rep}}_{K^{\mathrm{alg}}}(\mathcal{I}_{k^{\mathrm{alg}}(\!(t)\!)}\times\mathbb{G}_a)$
preserves the slopes (by corollary \ref{q-break decomposition} on
the left hand side, and by the Swan conductor on the right hand
side).
\end{corollary}
\begin{proof}
One has
$D_{\sigma_q}^{(\phi)}=D_d^{(\phi)}\circ\Conf^{\mathrm{Tay}}_q$.
Since $D_d^{(\phi)}$ and $\Conf^{\mathrm{Tay}}$ preserve the
slopes, so does $D_{\sigma_q}^{(\phi)}$.
\end{proof}

\vspace{-0.7cm}

\bibliographystyle{amsalpha}
\bibliography{corrections-to-CM1996}
\vspace{-0.3cm}


\end{document}